\documentclass[11pt,a4paper,
fleqn,reqno]{amsart}                 
\usepackage[T1]{fontenc}                   
\usepackage[utf8]{inputenc}                 
\usepackage[english]{babel}    
 \usepackage{csquotes}
\usepackage{graphicx}                      %
\usepackage[font=small]{quoting}    
\usepackage{float}
\usepackage{caption}

\usepackage[top=1.27in,bottom=1.4in,left=1.0in,right=1.0in]{geometry}                  
\usepackage{verbatim}
\usepackage{color}
\usepackage{xcolor}
\usepackage{picture}
\usepackage{amsmath,amsthm,amsfonts,amssymb}
\usepackage{comment}
\usepackage[draft]{hyperref}
\hypersetup{colorlinks,linkcolor={blue}, citecolor={blue},urlcolor={black}}   
\usepackage{enumitem}
\usepackage{mathtools}
\usepackage{wasysym}
 
\usepackage[url=false,doi=false,isbn=false, giveninits=true, 
style=
numeric,
maxbibnames=99]{biblatex}
\AtBeginBibliography{%
  \setlength{\emergencystretch}{3em}%
}

\usepackage{ragged2e}

\newtheorem{thm}{Theorem}[section] 
\newtheorem{lem}[thm]{Lemma} 
\newtheorem{cor}[thm]{Corollary} 
\newtheorem{prop}[thm]{Proposition} 
 
\newtheorem{defn}[thm]{Definition}

\theoremstyle{definition} 
\newtheorem{rem}{Remark}[section]

\theoremstyle{definition}

\usepackage{esint}

\usepackage{amsmath}

\def\O{\Omega}

\def\S{\Sigma} 
\def\n{\nabla}

\def\p{\partial}

\def\a{\alpha}
\def\b{\beta}
\def\n{\nabla}

\def\o{\omega}
\def\O{\Omega}
\def\p{\partial}

\def\a{\alpha}
\def\b{\beta}
\def\g{\gamma}
\def\d{\delta}
\def\k{\kappa}
\def\l{\lambda}
\def\s{\sigma}

\def\ov{\overline}

\def\n{\nabla}
\def\<{\langle}
\def\>{\rangle}

\def\n{\nabla}

\def\RR{\mathbb{R}}

\def\SS{\mathbb{S}}

\def\o{\omega}
\def\O{\Omega}
\def\p{\partial}

\def\a{\alpha}
\def\b{\beta}
\def\g{\gamma}

\def\d{\delta}
\def\l{\lambda}
\def\K{\mathcal{K}}
\def\s{\sigma}

\def\ov{\overline}

\def\wh{\widehat}

\def\C{\mathcal{C}}

\def\t{\tilde}
\def\ol{\overline}
\def\wt{\widetilde}

{\left\{\begin{array}{@{}l@{}}}{\end{array}\right.}
\patchcmd{\abstract}{\scshape\abstractname}{\textbf{\abstractname}}{}{}
\makeatletter 
\def\@makefnmark{} 
\makeatother

\numberwithin{equation}{section}
\numberwithin{exa}{section}

\begin{document}
	\title[Prescribed $L_p$ curvature problem for convex capillary hypersurface]
	{Prescribed $L_p$ curvature problem for convex capillary hypersurface}

\author[X. Mei]{Xinqun Mei}
\address[X. Mei]{Key Laboratory of Pure and Applied Mathematics, School of Mathematical Sciences, Peking University,  Beijing, 100871, P.R.China}

\email{\href{qunmath@pku.edu.cn}{qunmath@pku.edu.cn}}

\author[G. Wang]{Guofang Wang}
\address[G. Wang]{Mathematisches Institut, Albert-Ludwigs-Universit\"{a}t Freiburg, Freiburg im Breisgau, 79104, Germany}
\email{\href{guofang.wang@math.uni-freiburg.de}{guofang.wang@math.uni-freiburg.de}}

\author[L. Weng]{Liangjun Weng}
\address[L. Weng]{Centro di Ricerca Matematica Ennio De Giorgi, Scuola Normale Superiore, Pisa, 56126, Italy}
\email{\href{mailto:liangjun.weng@sns.it}  {liangjun.weng@sns.it}}



\begin{abstract}
We address the prescribed $L_p$ curvature problem for convex capillary
hypersurfaces in the Euclidean half-space. 
Using the support function, we formulate the problem as a Hessian
quotient equation on a spherical cap with a Robin boundary condition.
For $p\geq1$, we establish existence results under suitable assumptions
on the prescribed data, together with uniqueness in the corresponding
parameter ranges. At the critical exponent, we obtain an associated
eigenvalue problem with uniqueness up to dilation. In the subcritical
range, existence is established for capillary even data. As applications, we obtain
existence results for the capillary $L_p$ Christoffel--Minkowski
problem and prescribed $L_p$ quotient curvature problems, providing a
framework for studying these equations in the presence of a
capillary boundary.
\end{abstract}


\maketitle 
\raggedbottom
 

\section{Introduction}
The study of prescribed curvature problems lies at the core of modern convex geometry and geometric analysis. Over the past decades, a wide array of tools---ranging from the theory of fully nonlinear elliptic equations to geometric curvature flows---has been brought to bear on classical questions such as the Minkowski problem, the Christoffel--Minkowski problem, and their $L_p$ or duality extensions. These developments have revealed a striking unifying principle: many geometric inequalities, variational structures, and rigidity phenomena can be encoded in nonlinear curvature quantities, whose analytic behavior captures surprisingly refined geometric information. In contrast to the extensively studied closed setting, the prescribed curvature problem for hypersurfaces with boundary, particularly in capillary geometries, remains far less understood. The presence of a supporting hyperplane and a prescribed Young's contact angle creates boundary interactions that resist direct application of techniques successful in the closed case. This intrinsic boundary coupling introduces substantial analytic difficulties and necessitates the development of genuinely new ideas and techniques.

Building on this perspective, a natural direction is to investigate \textit{whether the modern theory of prescribed curvature problems admits an analogue in the capillary setting}, particularly for convex capillary hypersurfaces in the Euclidean half-space. Of special interest are $L_p$-type curvature quantities, which play a fundamental role in the variational formulation and existence theory of the classical $L_p$ Christoffel--Minkowski problems. In the presence of a boundary contact angle condition, however, the associated fully nonlinear equations are accompanied by nontrivial boundary interactions, leading to Hessian (or Hessian quotient) type equations coupled with Robin (or Neumann) boundary conditions. Understanding how convexity, capillary geometry, and analytic admissibility interact in this framework is a central theme of the present work. To set the stage, we now turn to the basic definition and example of capillary hypersurfaces in the Euclidean half-space, which will be studied throughout the paper.
 
A {\it capillary hypersurface} in the Euclidean half-space $\overline{\mathbb{R}^{n+1}_{+}}$ $(n \geq 2)$ is a properly embedded, smooth, compact hypersurface $\Sigma \subset \overline{\mathbb{R}^{n+1}_{+}}$ with nonempty boundary $\partial\Sigma \subset \partial\mathbb{R}^{n+1}_{+}$, whose boundary meets the supporting hyperplane at a constant contact angle $\theta \in (0,\pi)$. That is,
 \begin{eqnarray*}
     \<\nu, e\>=-\cos\theta,~~~~\text{ along } ~ \p \S,
 \end{eqnarray*}
where $\nu$ and $e=-E_{n+1}$ are the unit outward normal of $\S\subset\ol{\RR^{n+1}_+}$ and $\p\RR^{n+1}_+\subset \ol{\RR^{n+1}_+}$ respectively.  A basic example of a capillary hypersurface is the unit spherical cap $$\C_{\theta} \coloneqq  \left \{\xi\in \ov{\mathbb{R}^{n+1}_+} ~ \big| ~|\xi-\cos\theta e|= 1 \right\}.$$ Let $\k(X)\coloneqq (\k_1,\cdots, \k_n)$ denote the principal curvatures of $\S$ at $X\in \S$, and let $\s_k(\k(X))$ denote the $k$-th mean curvature of $\S$, see Section \ref{sec-2} for details. Let $\sigma$ be the round metric on $\mathcal C_\theta$, and let
$\mu$ denote the outward unit conormal to $\partial\mathcal C_\theta$.
The support function of the unit spherical cap is
\begin{eqnarray}\label{eqn-ell}
  \ell(\xi)=\langle\xi,\xi-\cos\theta\,e\rangle
  =\sin^2\theta+\cos\theta\langle\xi,e\rangle.
\end{eqnarray}
It satisfies $\nabla^2\ell+\ell\sigma=\sigma$ and
$\nabla_\mu\ell=\cot\theta\,\ell$ on $\partial\mathcal C_\theta$.

In this paper,  we study the prescribed $L_p$ curvature problem for convex capillary hypersurfaces in $\ov{\RR^{n+1}_{+}}$. More precisely,

\

\textit{Fix $p\ge 1$ and $0\leq l<k\leq n$. Given a positive smooth function $f$ on the spherical cap $\C_{\theta}$,  can one find a strictly convex capillary hypersurface $\S\subset \ov{\RR^{n+1}_{+}}$ such that
\begin{eqnarray}\label{lp-curvatur problem-c}
    \frac{\sigma_{n-l}(\kappa(X))}{\sigma_{n-k}(\kappa(X))}=f^{-1}(\wt \nu(X))\<X, \nu(X)\>^{1-p}~?
\end{eqnarray}}

\
 
 Here  
\begin{eqnarray*}
    \wt \nu(X)\coloneqq \nu(X)+\cos\theta e,
\end{eqnarray*}
is the \textit{capillary Gauss map} of $\S$, which is a diffeomorphism from $\S$ to $\C_{\theta}$ if $\S$ is strictly convex. See \cite[Section 2]{MWWX} for more discussion on $\wt \nu $.  

Before proceeding, we briefly recall some progress on the prescribed $L_p$ curvature problem in the closed setting. Let $M$ be a closed, strictly convex hypersurface in $\mathbb{R}^{n+1}$ enclosing the origin and satisfying
\begin{eqnarray}\label{lp-curvatur problem}
    \frac{\sigma_{n-l}(\kappa(X))}{\sigma_{n-k}(\kappa(X))}=g^{-1}({\nu}(X))\<X, \nu(X)\>^{1-p},~~~X\in M,
\end{eqnarray}
for some smooth positive function $g$ on $\SS^{n}$.  
When $p\neq1$, Eq.~\eqref{lp-curvatur problem}
characterizes homothetically evolving profiles of the
anisotropic curvature flow
\begin{eqnarray}\label{a-flow}
    \frac{\partial Y(\cdot, t)}{\partial t}={\rm{sign}}(p-1) \cdot\left(g( \nu(Y))\frac{\sigma_{n-l}(\kappa(Y))}{\sigma_{n-k}(\kappa(Y))}\right)^{\frac{1}{1-p}} \nu.
\end{eqnarray}
 In particular, when $g=1$, the classification of solutions to \eqref{lp-curvatur problem} is crucial for understanding the asymptotic behavior of the normalized flow \eqref{a-flow}, see e.g., \cite{BCD, CD, GLW-crell, 
huisken-sinestrari-1999, James} and references therein.
When $p=1$,  and $k=n$,  Eq. \eqref{lp-curvatur problem} corresponds to the prescribed $(n-l)$-th Weingarten curvature problem in classical differential geometry. For $l=0$, it is the classical Minkowski problem,  which has been completely solved through the work of Minkowski
\cite{Min}, Alexandrov \cite{Alex}, Lewy \cite{Lewy}, Nirenberg \cite{Nire}, Pogorelov \cite{Pog52}, Cheng--Yau \cite{CY76}. For
$1\leq l\leq n-1$, Guan--Guan \cite{GG02} proved the existence of a strictly convex hypersurface that satisfies Eq.  \eqref{lp-curvatur problem} when $g$ is invariant under an automorphic group $G$ without fixed points (for example, if $g$ is an even function on $\SS^n$) and Guan--Ma--Zhou \cite{gmz2006} solved the prescribed quotient curvature problem. When $k=n$ and $1\leq l\leq n-1$, Hu--Ivaki \cite{HI} established the existence of a strictly convex, even solution to Eq. \eqref{lp-curvatur problem} for $1<p<n-l+1$. If $p>n-l+1$, the existence of a strictly convex solution to Eq. \eqref{lp-curvatur problem} was established in \cite{GRW} by an elliptic method and in \cite{Bryan-ivkai-julian-2021} via a flow method.  For the case $p=n-l+1$, Eq. \eqref{lp-curvatur problem} is invariant under dilation, and we refer to it as the eigenvalue problem for the prescribed curvature equation. In  \cite{Lee}, Lee obtained the strictly convex, even solution to the eigenvalue problem through a compactness argument relying on the result in \cite{GRW}.  Recently, in  \cite{MWW-Quotient},  we removed the evenness assumption in \cite{Lee} and resolved the prescribed quotient curvature problem for a closed, strictly convex hypersurface. See also \cite{CNS-4, GLM06, STW}, among many others, for more related results.

Another important motivation for studying the prescribed $L_p$ curvature problem arises from the prescribed area measure problem in convex geometry. Given a closed, strictly convex hypersurface $M\subset\RR^{n+1}$, it induces a $k$-th area measure on $\SS^n$ by $$dS_k\coloneqq \s_k(\n^2 h+h\s) d\s,$$ where $h$ is the support function of the convex body enclosed by $M$ and $d\s$ is the standard volume form on $\SS^n$. Indeed, \eqref{lp-curvatur problem} is equivalent to finding a strictly convex function $h$ on $\mathbb{S}^{n}$ satisfying
\begin{eqnarray}\label{support-equ}
    \frac{\sigma_{k}(\n^{2}h+h\sigma)}{\sigma_{l}(\n^{2}h+h\s)}=gh^{p-1}, \quad {\rm{on}}~~\SS^{n}. 
\end{eqnarray}
In the case $k=n$, $l=0$, and $p \neq 1$, \eqref{support-equ} reduces to the celebrated $L_p$ Minkowski problem initiated by Lutwak \cite{Lut93}, which encompasses both the logarithmic Minkowski problem ($p=0$, see B\"or\"oczky--Lutwak--Yang--Zhang \cite{BLYZ}) and the centro-affine Minkowski problem ($p=-1-n$, see Chou--Wang \cite{CW06}).

For $1\leq k\leq n-1$ and $l=0$,
Eq.~\eqref{support-equ} is the $L_p$
Christoffel--Minkowski equation.
Under suitable structural assumptions on the prescribed
function, Guan--Ma \cite{GM} established existence for
$p=1$, and Hu--Ma--Shen \cite{HMS2004} treated
$p>k+1$ and the associated eigenvalue problem at $p=k+1$.
For $1<p<k+1$, Guan--Xia \cite[Theorem~1.1]{GX}
proved existence for positive even prescribed data satisfying
the spherical convexity condition in that theorem.
Many significant results have been achieved in the past three decades; it is nearly impossible to provide a comprehensive account here. We refer the interested reader to the classical monograph by Schneider \cite{Sch} and a more recent one by B\"or\"oczky--Figalli--Ramos \cite{BRF}.

In \cite{MWW-AIM}, we studied the prescribed Gauss--Kronecker curvature problem for convex capillary hypersurfaces in $\ov{\RR^{n+1}_{+}}$, which we refer to as the capillary Minkowski problem. A necessary and sufficient condition for the solution to the capillary Minkowski problem was obtained therein.  We further investigated the prescribed $k$-th Weingarten curvature problem for convex capillary hypersurfaces in \cite{MWW-Weingarten}, i.e., the problem of finding a strictly convex capillary hypersurface that satisfies 
\begin{eqnarray}\label{c-curvature}
    \sigma_{k}(\kappa(X))=f^{-1}(\wt \nu(X)).
\end{eqnarray} 
Under the condition that $f$ is a capillary even function on $\C_{\theta}$,  we obtained the existence of a strictly convex, capillary even hypersurface solving Eq. \eqref{c-curvature}.  A function $f\colon\C_\theta\to\RR$ is said to be  \textit{capillary even} on $\C_\theta$, if for any point $\xi \coloneqq  (\xi_1, \dots, \xi_n, \xi_{n+1}) \in \C_\theta$,  $$f(\xi_1,\cdots, \xi_n,\xi_{n+1})=f(-\xi_1,\cdots, -\xi_n, \xi_{n+1}).$$
See, e.g. \cite[Definition 1.1]{MWW-Weingarten}. We say that a capillary hypersurface is {capillary even}  if its support function is a {capillary even} function on $\C_\theta$. On the other hand, the capillary Minkowski problem can also be viewed as a prescribed surface area measure problem for the capillary hypersurface, see \cite[Section 1]{MWW-AIM}.  Subsequently,  in \cite{MWW-Lp}, we proposed the capillary $L_p$-Minkowski problem for capillary convex bodies with $p\in \RR$, which seeks to find a capillary convex body with a prescribed capillary $L_p$-surface area measure in the Euclidean half-space. This problem corresponds to the solvability of the following Monge-Amp\`ere equation with a Robin boundary condition:  
\begin{eqnarray*}
\left\{
		\begin{array}{rcll}\vspace{2mm}\displaystyle
			\det(\n^2 h+h\s )&=&  f h^{p-1},& \quad    \hbox{ in } ~\C_{\theta},\\
			\n_\mu h &=& \cot\theta  h, &\quad  \hbox{ on } ~\p \C_\theta,
		\end{array}\right.
	\end{eqnarray*}
    where  $\n h$ and $\n^2h$ are the gradient and the Hessian of $h$ on $\C_\theta$ with respect to the standard spherical metric $\s$ on $\C_{\theta}$, respectively,  $\mu$ is the unit outward normal of $\p \C_\theta \subset \C_\theta$.
    {Existence results for $p\geq1$ were established by elliptic methods in \cite{MWW-AIM, MWW-Lp}, under the assumptions stated
there and with an eigenvalue formulation at the critical exponent.
The capillary even $L_p$ Minkowski problem was subsequently
studied for $-n-1<p<1$ in \cite{H-I} and, by a flow
approach for $p>-n-1$ in \cite{H-H-I}.}

Motivated by the aforementioned developments for closed convex hypersurfaces and the capillary 
\(L_p\) Minkowski problem, it is natural to consider the capillary \(L_p\) Christoffel--Minkowski problem, and, more generally, prescribed $L_p$ curvature problem for convex capillary hypersurfaces. 
In view of \cite[Proposition~2.4]{MWW-AIM} (see also \cite[Lemma~2.4 and Proposition~2.6]{MWWX}), the problem \eqref{lp-curvatur problem-c} can be reformulated as finding a strictly convex function $h$ satisfying the following Hessian quotient equation with a Robin boundary condition:
\begin{eqnarray}\label{quo-equ}
\left\{	\begin{array}{rcll}\vspace{2mm}\displaystyle
			\frac{\sigma_{k}(\n^{2}h+h\sigma)}{\sigma_{l}(\n^{2}h+h\sigma)}&=&f h^{p-1} &\quad    \hbox{ in } ~\C_{\theta},\\
			\n_\mu h &=& \cot\theta  h, &\quad  \hbox{ on }~ \p \C_\theta.
		\end{array}\right.
	\end{eqnarray} 

Throughout, $W[h]\coloneqq\nabla^2h+h\sigma$.
Admissibility means $\lambda(W[h])\in\Gamma_k$; convexity and
strict convexity mean $W[h]\ge0$ and $W[h]>0$, respectively
(see Definition~\ref{def}). 
We focus on the acute-angle range $\theta\in(0,\pi/2)$.
At $\theta=\pi/2$, the Robin boundary condition becomes a Neumann boundary 
condition. For data and solutions admitting smooth even
reflection across the equator, this gives a corresponding
reflection-symmetric problem on $\SS^n$.

Our first main result gives existence and uniqueness for
\eqref{quo-equ} when $p>k-l+1$, and an associated eigenvalue
problem when $p=k-l+1$.
Under additional structural assumptions on $f$,
the admissible solutions are strictly convex.
\begin{thm}\label{thm-admissible-existence}
Let $0\leq l<k\leq n$ and $\theta\in (0, \pi/2)$. Then, for any positive, smooth function $f$ on $\C_{\theta}$, the following results hold.
\begin{enumerate}
    \item  If $p>k-l+1$, then there exists a unique,  smooth, admissible solution $h$ solving Eq. \eqref{quo-equ}.
    \item If $p=k-l+1$, then there exists a  positive, smooth and admissible solution $h$  and a positive constant $\tau$ solving 
\begin{eqnarray}\label{homo-equ}
\left\{
		\begin{array}{rcll}\vspace{2mm}\displaystyle
			\frac{\sigma_{k}(\n^{2}h+h\s)}{\sigma_{l}(\n^{2}h+h\s)}&=&  \tau f h^{k-l}, & \quad    \hbox{in }~~ \C_{\theta},\\
			\n_\mu h&=&\cot\theta h, &\quad  \hbox{ on }~~ \p \C_\theta.
		\end{array}\right.
  \end{eqnarray} 
 Moreover, $\tau$ is unique and the positive admissible solution $h$ is unique up to dilations.
  \end{enumerate}
Furthermore, when $k\leq n-1$,  assume that $f$ satisfies
 \begin{eqnarray}\label{f-condition}
\left\{
		\begin{array}{rcll}\vspace{2mm}\displaystyle   
      \n^{2}f^{-\frac{1}{p+k-l-1}}+f^{-\frac{1}{p+k-l-1}}\sigma&\geq& 0, ~~{\text{in}}~~~\C_{\theta},\\
      \n_{\mu} f+\cot\theta (p+k-l-1)f&\geq& 0,~~{\text{on}}~~~\partial\C_{\theta}.
 	\end{array}\right.
	\end{eqnarray} 
    Then the admissible solutions in (1) and (2)  are strictly convex. 
\end{thm}
For the case $1<p<k-l+1$,  Eq. \eqref{quo-equ} is degenerate due to the absence of a uniform positive lower bound for the solution. For the classical $L_p$ Christoffel--Minkowski problem, Guan--Xia \cite[Section 5]{GX} constructed examples that vanish at some point but satisfy Eq. \eqref{support-equ} with $l=0$, some $p\in (1,k+1)$, and $f$ being a positive, smooth function.   These examples indicate that a certain geometric obstruction on the prescribed function is needed for the existence of a solution to the classical $L_p$ Christoffel--Minkowski problem.  In this paper, we focus on the existence of capillary even, strictly convex solutions to Eq. \eqref{quo-equ}.

\begin{thm}\label{thm-p-small}
 Let $0\leq l<k\leq n$ and $\theta\in (0, \pi/2)$. Let $1\leq p<k-l+1$ and  $f$ be a positive, smooth, capillary even function on $\C_{\theta}$. If $k\leq n-1$,  assume that $f$ satisfies conditions \eqref{f-condition}.  Then there exists a positive, smooth, capillary even, and strictly convex solution $h$ to  Eq. \eqref{quo-equ}.  When $l=0$ and $1<p<k+1$, the solution is unique. 
\end{thm}

We remark that two cases: $l=0$, $1\leq k \leq n-1$, $1< p <k+1$ and $ k=n$, 
$1\leq l\leq n-1$, $1< p<n-l+1$ were studied independently by Hu--Ivaki in \cite{hu-ivakilpchristoffelminkowskiproblem} and \cite{hi-lpcurvature} respectively. For more  details see Remark \ref{rem-update} below.

As a direct consequence of Theorem~\ref{thm-admissible-existence} and Theorem \ref{thm-p-small}, we solve the prescribed capillary $k$-th $L_{p}$-surface area measure problem for capillary convex bodies. More precisely, for a capillary convex body \(\widehat{\Sigma}\subset \ol{\RR^{n+1}_+}\), we define the corresponding capillary \(k\)-th \(L_{p}\)-surface area measure by
\begin{eqnarray}\label{eqn:intro-Lp-kth-measure}
dS_{\widehat{\Sigma},k,p}^{c} \coloneqq \ell^{p} h_{\widehat{\Sigma}}^{1-p} \, dS_{\widehat{\Sigma}, k}= \ell^{p} h_{\widehat{\Sigma}}^{1-p} \s_k \left(W[h_{\wh\S} ]\right) d\s.
\end{eqnarray}
We refer to the corresponding prescribed measure problem as the \emph{capillary \(L_p\) Christoffel--Minkowski problem}, see Section~\ref{sec-2.3} for its variational formulation and more discussion. In terms of the support function \(h\coloneqq h_{\wh\S}\), it is equivalent to the following Hessian equation with a Robin boundary condition:
\begin{eqnarray}\label{eqn-intro-Lp-CM-eq}
\left\{ \begin{array}{rcll}\vspace{1mm} \s_k(\nabla^{2} h + h \, \sigma) &=& f \, h^{\,p-1}, & \quad \text{ in } ~\C_{\theta},\\ \nabla_\mu h &=& \cot \theta \, h, & \quad \text{ on } ~\partial \C_{\theta}. \end{array} \right. \end{eqnarray}
In particular, when $p=1$ and $1\leq k\leq n-1$, the problem reduces to the capillary Christoffel--Minkowski problem. Existence results were obtained independently, under their respective hypotheses, by Hu--Ivaki--Scheuer \cite{HIS2025} and by the authors \cite{MWW-CM}; in the latter work, an additional assumption is the existence of a suitable homotopy path. Thus, the capillary \(L_p\) Christoffel--Minkowski problem naturally interpolates between the capillary $L_p$ Minkowski problem and the capillary Christoffel--Minkowski problem.
We now turn to the general case $p>1$ and $1\le k\le n-1$ for \eqref{eqn:intro-Lp-kth-measure} (or \eqref{eqn-intro-Lp-CM-eq}) and state the main result in the following. 
 \begin{thm}\label{thm-cap-Lp-CM-thm}
   Let $1\leq k\leq n-1$ and  $\theta \in (0, \pi/2)$.  Let  $f$ be a positive smooth function on $\C_{\theta}$. Assume that $f$ satisfies the condition \eqref{f-condition} with $l=0$. The following results hold:
     \begin{enumerate}
         \item When $p>k+1$,  there exists a unique capillary convex body such that its capillary $k$-th $L_p$-surface area measure is equal to $ f\ell^{p}d\sigma$.
         \item When $p=k+1$,  there exists a unique positive constant $\tau$ and a unique capillary convex body (up to dilation)  such that its capillary $k$-th  $L_p$-surface area measure is equal to $\tau  f\ell^{p}d\sigma$.
         \item When $1< p<k+1$ and $f$ is capillary even,   there exists a unique capillary even convex body such that its capillary $k$-th $L_p$-surface area measure is equal to $f\ell^{p}d\sigma$.
     \end{enumerate}
 \end{thm} 
 \begin{rem}In the case \(\theta=\pi/2\), the capillary $L_p$ Christoffel--Minkowski problem reduces, by reflection across the boundary hyperplane, to the classical \(L_p\) Christoffel--Minkowski problem for closed convex hypersurfaces. Consequently, Theorem~\ref{thm-cap-Lp-CM-thm} extends several results in the classical setting. More precisely, parts (1) and (2) correspond to the result of Hu--Ma--Shen \cite[Theorem~1. (i), (ii)]{HMS2004}, while part (3) recovers the even case established by Guan--Xia \cite[Theorem~1.1]{GX}. Therefore, Theorem~\ref{thm-cap-Lp-CM-thm} may be viewed as a natural extension of the classical \(L_p\) Christoffel--Minkowski theory from closed convex hypersurfaces to capillary convex hypersurfaces with contact angle \(\theta\in(0,{\pi}/{2})\). On the other hand, the corresponding problem for obtuse contact angle \(\theta\in({\pi}/{2},\pi)\) remains widely open. The difficulty is rooted in both analytical and geometric aspects. Analytically, the sign structure that underlies the maximum principle in the acute-angle regime breaks down when \(\theta>{\pi}/{2}\), leading to a fundamental obstruction in deriving the a priori estimates, in particular the global $C^2$ estimate or curvature estimate. Geometrically, convex capillary hypersurfaces with obtuse contact angle are no longer convex in the sense of Busemann \cite{Busemann} for hypersurfaces with boundary. These difficulties suggest that new ideas will be required to develop an \(L_p\) Christoffel--Minkowski theory in the obtuse capillary setting. \end{rem}
 
Another application of Theorem~\ref{thm-admissible-existence} and Theorem \ref{thm-p-small} is the solvability of the prescribed $L_p$ curvature problem and related eigenvalue problem for convex capillary hypersurfaces. These results can be viewed as a natural continuation of, and a culmination of, our previous series of works \cite{MWW-Weingarten, MWW-CM, MWW-AIM, MWW-Lp}.
    \begin{thm}\label{thm:1.4}
    Let $0\leq l<k\leq  n$ and $\theta \in (0, \pi/2)$. Assume that $f$ is a positive, smooth function on $\C_{\theta}$. If $k\leq n-1$, suppose that $f$ satisfies condition \eqref{f-condition}. Then the following results hold.
   \begin{enumerate}
       \item If $p>k-l+1$, then there exists a unique, smooth, strictly convex capillary hypersurface $\S\subset\overline{\mathbb{R}^{n+1}_{+}}$ satisfying Eq. \eqref{lp-curvatur problem-c}.
       \item If $p=k-l+1$,  then there exists a  unique (up to a dilation) smooth, strictly convex capillary hypersurface $\S\subset\overline{\mathbb{R}^{n+1}_{+}}$ and  a positive constant $\tau$ satisfying the following eigenvalue problem: 
        \begin{eqnarray*} 
\<X, \nu(X)\>^{k-l}\frac{\sigma_{n-l}(\kappa(X))}{\sigma_{n-k}(\kappa(X))}=\tau^{-1} f^{-1}(\wt \nu(X)),\quad \forall X\in \S.
        \end{eqnarray*}  
    \item If $1\leq p<k-l+1$. Assume that $f$ is capillary even, then there exists a  smooth, capillary even, and strictly convex capillary hypersurface $\S\subset\overline{\mathbb{R}^{n+1}_{+}}$ satisfying Eq. \eqref{lp-curvatur problem-c}.       
    \end{enumerate}
    \end{thm}

We conclude the introduction by outlining the main ideas of the proofs. To establish the existence results in Theorems~\ref{thm-admissible-existence} and \ref{thm-p-small}, we employ the continuity method and the topological degree method, respectively. In both cases, the central ingredient is the derivation of uniform a priori estimates for admissible solutions of \eqref{quo-equ}. Once these estimates are available, the existence theory follows from standard continuity and degree arguments.

\ 

For the range $p>k-l+1$, the $C^0$ estimate follows directly from the maximum principle. In contrast to the case of convex hypersurfaces, however, this estimate no longer yields a gradient bound when $ k\le n-1$. Indeed, admissible solutions need not be convex in this setting, and the lack of convexity destroys the classical implication from $C^0$ to $C^1$ estimates (see, for example, \cite{CY76} and \cite[Lemma 3.2]{MWW-AIM}). Consequently, the gradient estimate requires a separate argument.

The most delicate part of the analysis is the derivation of global $C^2$ estimates, especially the control of the double normal derivative on the boundary. The difficulty originates from the boundary terms appearing in the linearized equation. For convex solutions, the unfavorable term
$$-\sum_{i,j,k} F^{ij} W_{ik} {d}_{kj}$$
(see, e.g., Eq. \eqref{eqn:FijQij}) has a favorable sign and can be controlled by $C f^{\frac1{k-l}}$, where $C>0$ depends only on $n$. Without convexity, however, this mechanism breaks down completely, and no direct estimate is available. To overcome this difficulty, we adapt a boundary barrier argument inspired by the works of Lions--Trudinger--Urbas \cite[Section~3]{LTU} and Ma--Qiu \cite[Section~4]{MQ}, and subsequent refinement of Chen--Zhang \cite[Section~4]{CZ2021}. More precisely, we construct suitable auxiliary functions involving the normal derivative $\n_\mu h$ in a thin neighborhood of $\p\C_\theta$. Combined with a careful analysis of the linearized operator, these barrier functions enable us to control the problematic boundary terms and establish the required double normal derivative estimate. This provides the crucial ingredient for the global $C^2$ estimate and, together with the Evans--Krylov theorem and Schauder theory, yields the higher order estimates. We refer the reader to Section~\ref{sec3.2}, especially Lemma~\ref{c2 est-admissible soluion-big case}, for the detailed argument.

\ 

  When $p=k-l+1$, Eq.~\eqref{quo-equ} becomes dilation invariant, and therefore no \emph{a priori} $C^0$ estimate can be expected for its positive solutions. Following the strategy developed in Guan--Lin \cite[Section~3]{GL2000} and Hu--Ma--Shen \cite[Section~4]{HMS2004}, we instead seek a scale-invariant estimate. More precisely, the key step is to establish the logarithmic gradient estimate
\begin{eqnarray*}
    \|\n \log h\|_{C^0(\C_\theta)}\leq C.
\end{eqnarray*}Once this estimate is available, the lack of normalization caused by the dilation invariance can be overcome.     Compared with the classical closed hypersurface case, the presence of the Robin boundary condition complicates the analysis. In particular, the boundary contributions produced by the maximum principle cannot be controlled by the standard auxiliary functions used in the literature. To overcome this difficulty, motivated by the techniques developed in Ma--Xu \cite[Section~3]{MX} and Ma--Qiu \cite[Section~3]{MQ} for Neumann boundary value problems, we introduce the auxiliary function
\begin{eqnarray*}
    \omega\coloneqq \log h-\cot\theta \tilde{\ell},
\end{eqnarray*}
where $\wt\ell$ is the positive function defined in \eqref{t-l}. The additional correction term compensates for the boundary contribution arising from the Robin boundary condition and allows the maximum principle to be applied effectively. This leads to the desired logarithmic gradient estimate, see Lemma~\ref{lem kC1} for details.

\ 

For $1\leq p<k-l+1$, we establish positive lower and upper bounds
for the support function under the capillary evenness assumption.
Our argument combines a quantitative second derivative estimate
with the geometric radius inequality in Lemma~\ref{chou-wang-lemma}.
The endpoint $p=1$ and the range $1<p<k-l+1$ are treated separately.

When $p=1$, the terms arising from differentiating
$h^{(p-1)/(k-l)}$ vanish. The interior calculation and the boundary
barrier argument yield
\begin{eqnarray*}
\max\limits_{\C_\theta}|\n^2h+h\s|
\leq C\left(1+\|h\|_{C^0(\C_\theta)}\right).
\end{eqnarray*}
See Lemma~\ref{lem-p-one-C2}. This argument follows the boundary
estimate used in \cite[Lemmas~3.4 and~3.5]{MWW-Weingarten}.

For $1<p<k-l+1$, the weighted gradient estimate and the boundary
double normal derivative estimate give
\begin{eqnarray*}
\max\limits_{\C_\theta}|\n^2h+h\s|
&\leq& C\left[1+\|h\|_{C^0(\C_\theta)}
\left(1+(\min\limits_{\C_\theta}h)^{\frac{p-1}{k-l}-1}\right)
+(\max\limits_{\C_\theta}h)^{2-\gamma}
(\min\limits_{\C_\theta}h)^{\gamma-1}\right],~~~~
\end{eqnarray*}
where $\gamma>0$ is chosen as in Lemma~\ref{lem weighted gra}.
The constant $C$ depends only on $n,k,l,p,\theta,\gamma$,
$\min\limits_{\C_\theta}f$, and $\|f\|_{C^2(\C_\theta)}$.

In each case, we combine the corresponding second derivative
estimate with the elementary bounds
\begin{eqnarray*}
\min\limits_{\C_\theta}h\leq C,
\qquad \max\limits_{\C_\theta}h\geq c>0,
\end{eqnarray*}
Lemma~\ref{chou-wang-lemma}, and
Proposition~\ref{prop-radius-control}.
The capillary evenness assumption then allows us to obtain the required
two-sided bounds for $h$ and the global $C^2$ estimate.
The higher order estimates provide the compactness needed
for the degree argument in Section~\ref{sec-5}.

The strategy of combining second derivative estimates with
geometric radius inequalities to derive a priori bounds was
used by Chou--Wang
\cite[Lemmas~2.2, 2.4 and~2.5]{CW00}
in their study of a logarithmic Gauss curvature flow. Related weighted gradient estimates for the $L_p$
Christoffel--Minkowski problem on $\mathbb{S}^n$ were
established by Guan--Xia \cite{GX} and further developed
by Guan \cite{Guan2023}.
For the closed prescribed $L_p$ curvature problem, Hu--Ivaki
\cite[Lemma~2.2 and Proposition~2.3]{HI}
proved a weighted gradient estimate and combined it with
the geometric inequality of Chou--Wang to obtain two-sided bounds and the global $C^2$ estimate for even solutions.
In our case, for $1<p<k-l+1$, we use the weighted gradient estimate
argument of \cite[Lemma 3.1]{hi-lpcurvature}
, together with the capillary analogue of the Chou--Wang
inequality established in
\cite[Lemma~2.3]{MWW-Weingarten} and our 
estimates
for the capillary Hessian quotient equation, including
the boundary second derivative estimates above for the
Robin boundary condition.

For $1<p<k-l+1$, the range of contact angles covered by
this approach depends on the weighted gradient estimate.
We therefore clarify below the relation between the results
of the first version, the related works of Hu--Ivaki, and the estimates used in the present revision.

 \begin{rem}\label{rem-update}
The first arXiv version of this paper
(arXiv:2512.16686v1) established existence for the general
Hessian quotient problem \eqref{quo-equ} in the range
$1<p<k-l+1$, under the angle restriction
$\theta\in \left(\arccos \left((p-1)/(k-l) \right),  \pi /2 \right)$. The results for $p=1$ and $p\geq k-l+1$ are proved for every $\theta \in (0,\pi/2)$. The additional angle restriction in the subcritical range $1<p<k-l+1$ was imposed solely in
order to obtain the weighted gradient estimate \eqref{weighed gredient est} in Lemma \ref{lem weighted gra}.

Independently, Hu--Ivaki
\cite[Theorem~1.1]{hu-ivakilpchristoffelminkowskiproblem} proved Theorem \ref{thm-p-small} 
for  $l=0$, $1\leq k\leq n-1$ and $1<p<k+1$ and every $\theta\in(0,\pi/2)$.
Their proof uses a non-collapsing argument.

Later, in the case $k=n$, $1\leq l \leq n-1$ and $1<p<n-l+1$,   Hu--Ivaki \cite[Theorem~1.1]{hi-lpcurvature} removed our additional angle restriction and obtained the existence result for every $\theta\in(0,\pi/2)$. 
 The key new estimate in \cite[Lemma~3.1]{hi-lpcurvature} uses a weighted gradient estimate based
on the auxiliary function
$\ell^{2-\gamma}|\n u|^2/u^\gamma$, with $u\coloneqq h/\ell$;
as observed in that work, the weighted gradient argument also applies to more
general homogeneous quotient operators.

In the present revised version, we adopt their auxiliary function for the weighted gradient-estimate
argument for the Hessian quotient $\sigma_k/\sigma_l$ in \eqref{quo-equ}.
Together with our boundary second derivative estimates, this yields the
existence statement of Theorem~\ref{thm-p-small} for all $0\leq l<k\leq n$,
$1<p<k-l+1$ and for every $\theta\in(0,\pi/2)$. The corresponding geometric applications of Theorems \ref{thm-cap-Lp-CM-thm} (3) and  \ref{thm:1.4} (3) are updated accordingly. 

 \end{rem}

\

{\textbf{Organization of the paper.}}
In Section~\ref{sec-2}, we recall basic properties of elementary symmetric functions and review the geometry of capillary convex hypersurfaces (or bodies). From a variational viewpoint, we introduce the capillary $k$-th $L_p$-surface area measure and reformulate the capillary $L_p$ Christoffel--Minkowski problem. Section~\ref{sec-3} is devoted to a priori estimates for admissible (and strictly convex) solutions of Eq.~\eqref{quo-equ}. In Section~\ref{sec-4}, under additional assumptions on $f$, we establish the strict convexity of solutions in the case $k \le n-1$. Finally, Section~\ref{sec-5} contains the existence and uniqueness results for admissible (strictly convex) solutions to Eq.~\eqref{quo-equ}.
 
\

\section{Preliminaries}\label{sec-2}
In this section, we first recall several well-known properties of the elementary symmetric functions. We then review some geometric aspects of capillary convex bodies and introduce the notion of the capillary $k$-th $L_p$ surface area measure associated with such bodies.

\subsection{Elementary symmetric functions}
In this subsection, we recall the definition and some basic properties of the elementary symmetric polynomial functions.  
\begin{defn}
For $k = 1, 2,\ldots, n,$ the $k$-th elementary symmetric function $\sigma_k$ is defined by
\begin{eqnarray*} 
\sigma_k(\lambda) = \sum _{1 \le i_1 < i_2 <\cdots<i_k\leq n}\lambda_{i_1}\lambda_{i_2}\cdots\lambda_{i_k},
 \qquad \text {for} \quad\lambda\coloneqq (\lambda_1,\ldots,\lambda_n)\in \mathbb{R}^{n}.
\end{eqnarray*}
\end{defn}
We use the convention that $\sigma_0=1$ and $\sigma_k =0$ for $k>n$. Let $H_k(\lambda)$ be the normalization of 
$\sigma_{k}(\lambda)$ given by \begin{eqnarray*}
    H_k(\lambda)\coloneqq\frac{1}{\binom{n}{k}}\s_k(\lambda).
\end{eqnarray*} 
Denote  $\sigma _k (\lambda \left| i \right.)$ the symmetric
	function with $\lambda_i = 0$ and $\sigma _k (\lambda \left| ij \right.)$ the symmetric function with $\lambda_i =\lambda_j = 0$.  Recall that the  G{\aa}rding cone is defined as
\begin{eqnarray*} 
\Gamma_k \coloneqq \left\{ \lambda  \in \mathbb{R}^n \mid \sigma _i (\lambda ) > 0,~~\forall 1 \le i \le k \right\}.
\end{eqnarray*} 
\begin{defn}
    Let $A=\{A_{ij}\}$ be an $n\times n$ symmetric matrix. For $k=1,2,\cdots, n$, we define
    \begin{eqnarray*}
        \sigma_{k}(A)\coloneqq \sigma_{k}(\lambda(A))=\sum\limits_{1\leq i_{1}<i_{2}\cdots< i_{k}\leq n}\lambda_{i_{1}}(A)\lambda_{i_{2}}(A)\cdots \lambda_{i_{k}}(A),
    \end{eqnarray*}
where $\lambda(A)\coloneqq (\lambda_1(A),\lambda_2(A),\ldots,\lambda_n(A))$ denotes the vector of eigenvalues of $A$. 
Equivalently, $\sigma_k(A)$ is the sum of all $k\times k$ principal minors of $A$.
\end{defn}
We also denote by $\sigma _m (A \left|
i \right.)$ the symmetric function with $A$ deleting the $i$-row and
$i$-column and $\sigma _m (A \left| ij \right.)$ the symmetric
function with $A$ deleting the $i,j$-rows and $i,j$-columns.

\begin{defn}\label{def}
Let $1\le k\le n$ and let $h\in C^{2}(\mathcal C_{\theta})$. 
   We say that $h$ is \emph{$k$-convex} in $\mathcal C_{\theta}$ if
$$
\lambda\bigl(W[h]\bigr)\in\Gamma_k \quad \text{in } ~ ~\mathcal C_{\theta}.
$$
Moreover, $h$ is said to be \emph{convex} (respectively, \emph{strictly convex}) in $\mathcal C_{\theta}$ if $W[h]$ is positive semi-definite (respectively, positive definite) in $\mathcal C_{\theta}$.
\end{defn}

\begin{prop}\label{prop2.1}
Suppose that  $A=\{A_{ij}\}$ is diagonal matrix, and $1\leq k\leq n$. 
Then
\begin{eqnarray*}
\sigma_{k-1}^{ij}(A)= \begin{cases}
\sigma _{k- 1} (A\left| i \right.), &\text{ if } i = j, \\
0, &\text{ if } i \ne j,
\end{cases}
\end{eqnarray*}
where $\sigma_{k-1}^{ij}(A)\coloneqq \frac{{\partial \sigma _k (A)}} {{\partial A_{ij} }}$.
\end{prop}
\begin{prop}\label{pro-2.3}\ 
The following four properties hold:
    \begin{enumerate}
        \item  For $\lambda \in \Gamma_k$ and $k > l \geq 0$, $ r > s \geq 0$, $k \geq r$, $l \geq s$, 
\begin{eqnarray*}  
\left(\frac{H_{k}(\lambda)}{H_{l}(\lambda)}\right)^{\frac{1}{k-l}}\leq \left(\frac{H_{r}(\lambda)}{H_{s}(\lambda)}\right)^{\frac{1}{r-s}}.
\end{eqnarray*}
If $(k,l)\neq(r,s)$, equality holds if and only if
$\lambda_1=\cdots=\lambda_n>0$.
\item For $0\le l<k\le n$, $\left(\frac{H_k(\lambda)}{H_{l}(\lambda)}\right)^{\frac{1}{k-l}}$ is concave in $\Gamma_k$.

\item For $0\leq l<k\leq n$, denote  $F(\lambda)\coloneqq \left(\frac{\sigma_{k}(\lambda)}{\sigma_{l}(\lambda)}\right)^{\frac{1}{k-l}}$, $\lambda\in \Gamma_k$.
If $\lambda_{1}\geq \lambda_{2}\geq\cdots\geq \lambda_{n}$, then
    \begin{eqnarray*}
      \frac{\partial F}{\partial \lambda_{1}}\leq \frac{\partial F}{\partial\lambda_{2}}\leq \cdots\leq \frac{\partial F}{\partial\lambda_{n}}.  
    \end{eqnarray*}
\item For $0\leq l<k\leq n$. If $\l\in \Gamma_k$, then 
\begin{eqnarray*}
\sum\limits_{i=1}^{n} \frac{\p F(\l) }{\p \l_i} 
          \geq \left[\binom{n}{k}
            \Big/ \binom{n}{l}\right]^{\frac{1}{k-l}}\eqqcolon c_{n,k,l} .
		\end{eqnarray*} 
 \end{enumerate}
\end{prop}
For proofs of Propositions \ref{prop2.1}--\ref{pro-2.3}, see, e.g., \cite[Chapter XV, Section 4]{L96} and \cite[Lemma 2.10, Theorem 2.11, Theorem 2.16, Lemma 1.5]{Spruck}, respectively. 
 
\subsection{Geometric properties of capillary convex bodies}
In this subsection, we briefly review some properties of the capillary convex bodies, and for more details concerning the geometry of capillary convex bodies, we refer the reader to \cite[Section 2]{MWW-AIM} and \cite[Section 2]{MWWX}.

In this paper, we work in the smooth, strictly convex category.
We call $\widehat{\Sigma}$ a \emph{capillary convex body} if it is
a compact convex domain in $\ol{\RR^{n+1}_+}$ whose boundary consists
of a smooth, strictly convex capillary hypersurface
$\Sigma\subset\ol{\RR^{n+1}_+}$ and its flat boundary in
$\p\RR^{n+1}_+$, with
$\operatorname{int}\Sigma\subset\RR^{n+1}_+$. We denote
\begin{eqnarray*}
    \widehat{\p\Sigma} \coloneqq \p(\widehat{\Sigma}) \setminus \Sigma \subset \p\RR^{n+1}_+ .
\end{eqnarray*}
Let $\mathcal{K}_{\theta}$ be the collection of all capillary convex bodies in $\ol{\RR^{n+1}_+}$, and let $\mathcal{K}^{\circ}_{\theta}$ denote the subfamily consisting of those capillary convex bodies whose flat boundary contains the origin as a relative interior point.  
For $\widehat{\Sigma}\in \mathcal{K}_{\theta}$, let $\nu$ denote the unit outward normal vector field of $\Sigma$ in $\ol{\RR^{n+1}_+}$. The \emph{capillary Gauss map} of $\Sigma$ is defined by
\begin{eqnarray*}
    \wt \nu \coloneqq \nu + \cos\theta\, e .
\end{eqnarray*}Then $\wt \nu$ defines a diffeomorphism from $\Sigma$ onto the spherical cap $\C_{\theta}$, see \cite[Lemma 2.2]{MWWX}.  
The support function of $\Sigma$ is given by 
	\begin{eqnarray*}
		h(X)=\<X,\nu(X)\>.
	\end{eqnarray*}
	By the parametrization of $\wt\nu^{-1}$, we can view $h$ as a function on $\C_\theta$ by
	\begin{eqnarray}\label{support}
		h(\xi) =\<X(\xi),\nu(X(\xi ))\>=
		\< \wt \nu  ^{-1}(\xi)  , \xi -\cos\theta e\>,
	\end{eqnarray}and satisfies
 \begin{eqnarray} \label{robin}
     \n_\mu h=\cot\theta h, ~~~~~~~~~\text{ on } \p \C_\theta.
 \end{eqnarray} cf. \cite[Lemma 2.4]{MWWX}. We still call $h$ in \eqref{support} the support function of $\S$ (or of $\widehat{\S}$), which is only defined on $\C_\theta$.

By viewing $\wh\S\subset \ol{\RR^{n+1}_+}$ as a closed convex body in $\RR^{n+1}$, we recall several notions of radius for $\wh\S$ used in \cite[Section~2.2]{SW2024} and \cite[Section~2]{MWW-GCF}.   These notions will play a crucial role in controlling the geometry of $\S$, in particular in obtaining uniform diameter estimates.   The standard notion of the inner radius of  $\widehat{\S}$ is defined as
\begin{eqnarray*}    \rho_{-}(\widehat{\S})\coloneqq \sup \left\{\rho>0 \mid B_{\rho}^{+}(x_{0})\subset \widehat{\S}~\text{for~some~}x_{0}\in \partial  {\RR^{n+1}_{+}} \right\},\end{eqnarray*} and the outer radius of $\widehat\S$ is  

\begin{eqnarray*}    \rho_{+}(\widehat{\S})\coloneqq \inf \left\{\rho>0 \mid \widehat{\S}\subset B_{\rho}^{+}(x_{0})~\text{for~some~}x_{0}\in \partial  {\RR^{n+1}_{+}} \right\}, \end{eqnarray*} 
where $B^{+}_{\rho}(x_{0})\coloneqq B_{\rho}(x_{0})\cap \ol{\RR^{n+1}_+}$ and $B_{\rho}(x_0)$ is an open ball of radius $\rho$ centered at $x_0$ in $\RR^{n+1}$. The  notion of the capillary inner radius of $\wh\S\subset\ol{\RR^{n+1}_+}$ is defined as
\begin{eqnarray*}
\rho_{-}(\widehat\Sigma, \theta)\coloneqq\sup \left\{r>0 ~\mid~ \widehat{\C_{r,\theta}(x_0)}\subset \widehat\S\text{ for some } x_0\in \p   {\RR^{n+1}_+}\right\},\end{eqnarray*}and 
 the capillary outer radius of $\wh\Sigma$   as
\begin{eqnarray*}
\rho_+(\widehat\Sigma, \theta)\coloneqq \inf \left\{r>0 ~\mid~\widehat\S\subset \widehat{\C_{r,\theta}(x_0)} \text{ for some } x_0\in  \p  {\RR^{n+1}_+} \right\},
\end{eqnarray*}where 
\begin{eqnarray*}
    \C_{r,\theta}(x_0)\coloneqq \left\{x\in \ov{\RR^{n+1}_+} ~\mid~ |x-(x_0+r\cos\theta  e)|=r \right\}
\end{eqnarray*} 
is the spherical cap centered at $x_0+r\cos\theta e$ with radius $r>0$. 

The following proposition was shown in \cite[Proposition 2.3]{MWW-GCF} (or \cite[Proposition 2.4 and its proof]{SW2024}), which shows that for capillary convex bodies, the inner (outer) radius and the capillary inner (outer) radius can be mutually controlled when $\theta\in(0,\pi/2)$. 
\begin{prop}\label{prop-radius-control}
    Let $\widehat{\S}\in \K_{\theta}$ and $\theta \in (0, \pi/2)$. Then there hold
    \begin{eqnarray*}
         (1-\cos\theta )\rho_{-}(\widehat{\S}, \theta)\leq \rho_{-}(\widehat{\S})\leq \sin\theta \rho_{-}(\widehat{\S}, \theta),
    \end{eqnarray*}
    and 
    \begin{eqnarray*}
          (1-\cos\theta )\rho_{+}(\widehat{\S}, \theta)\leq \rho_{+}(\widehat{\S})\leq \sin\theta \rho_{+}(\widehat{\S}, \theta).
    \end{eqnarray*}
\end{prop}
Finally,  we introduce an important geometric lemma, which plays a crucial role in establishing the a priori estimates for the solution to \eqref{quo-equ} when $1\leq p<k-l+1$. This geometric lemma was originally obtained by Chou-Wang for smooth, closed, strictly convex hypersurfaces, as shown in \cite[Lemma 2.2]{CW00}, see also \cite[Lemma~15.11]{Ben-book}.  We obtain an analogous result for strictly convex capillary hypersurfaces in \cite[Lemma 2.3]{MWW-Weingarten}.

\begin{lem}[\cite{MWW-Weingarten}]\label{chou-wang-lemma}
    Let $\S$ be a strictly convex capillary hypersurface in $\ov{\RR^{n+1}_{+}}$ and $\theta\in (0, \pi/2)$. 
    Then there exists a positive constant $C=C(n)$ such that
    \begin{eqnarray*}
        \frac{\rho_{+}(\widehat{\S},\theta)^{2}}{\rho_{-}(\widehat{\S},\theta)}\leq C\sup\limits_{X\in  \S}\lambda_{X, \S} ,
   \end{eqnarray*}
 where $\lambda_{X,\Sigma}$ denotes the largest principal radius of curvature of $\Sigma$ at the point $X$.
\end{lem}

\subsection{The capillary \texorpdfstring{$k$-th $L_p$}{}-surface area measure}\label{sec-2.3}

In this subsection, we introduce the notion of the capillary $k$-th $L_p$ surface area measure for a capillary convex body
$\widehat{\Sigma}\in \mathcal{K}_{\theta}^{\circ}$. This notion extends the capillary $L_p$ surface area measure introduced in \cite[Definition~2.2]{MWW-Lp}. In \cite[Section~2]{MWW-Lp}, the capillary $L_p$ surface area measure was defined via the first variation of the volume of a capillary convex body under Firey’s $L_p$ perturbation. 
Motivated by Lutwak’s approach in \cite{Lut93}, we define the capillary $k$-th $L_p$ surface area measure by considering the variation of the $k$-th quermassintegral of capillary convex bodies under perturbations given by Firey’s $p$-sum.


For  $\varepsilon>0, p\geq 1$, and $K, L\in \mathcal{K}^{\circ}_{\theta}$, we denote  $K_\varepsilon\coloneqq K+_p\varepsilon L$ as the Firey $p$-sum of convex bodies $K$ and $L$.   From \cite[Definition 2.1, Proposition 2.1]{MWW-Lp}, we see that $K_{\varepsilon}$ is also a  capillary convex body in $\K_{\theta}^{\circ}$, and its support function $h_{K_{\varepsilon}}$ is given by 
\begin{eqnarray}\label{p-sum}
    h_{K_{\varepsilon}} (\cdot) \coloneqq \left[h_{K}^{p}(\cdot)+\varepsilon h_{L}^{p}(\cdot)\right]^{\frac{1}{p}},
\end{eqnarray}
where $h_{K}, h_{L}\in C^{2}(\C_{\theta})$ are the usual support functions of the convex bodies $K$ and  $L$, respectively. 
 From \cite[Eq. (1.8) and Proposition 2.6]{WWX-MA2024} and \cite[Lemma 2.14]{MWWX}, the $k$-th quermassintegral $\mathcal{V}_{k,\theta}(K)$ is given by 
\begin{eqnarray}\label{k-quermassintegral}
    \mathcal{V}_{k,\theta}(K)\coloneqq \frac{1}{n+1}\int_{\C_{\theta}}h_{K}H_{n-k}(W[h_{K}])d\s, ~~~~~0\leq k\leq n-1,
\end{eqnarray}
where $W[h_{K}]\coloneqq\n^{2}h_{K}+h_{K}\s$. 
\begin{defn}
    For $p\geq 1$ and $0\leq k\leq n-1$,  the capillary mixed $p$-quermassintegral of $K$ and $L$  is defined by 
\begin{eqnarray}\label{mixed-quemassintegral}
    \frac{n+1-k}{p}\mathcal{V}_{p, k}(K, L)\coloneqq \lim\limits_{\varepsilon \rightarrow 0^{+}}\frac{\mathcal{V}_{k, \theta}(K_{\varepsilon})-\mathcal{V}_{k, \theta}(K)}{\varepsilon}.
\end{eqnarray}
\end{defn}
Next, we establish the integral representation for the capillary mixed $p$-quermassintegral.
 
\begin{prop}
Let $K,L\in\mathcal{K}^{\circ}_{\theta}$ and let $0\le k\le n-1$. The capillary mixed $p$-quermassintegral defined in \eqref{mixed-quemassintegral} admits the integral representation:

\begin{eqnarray} \label{integral-pre}  \mathcal{V}_{p, k}(K,L)
=\frac{1}{(n+1)\binom nk}\int_{\C_{\theta}}h_{L}^{p}h_{K}^{1-p}dS_{K, n-k},
\end{eqnarray}
where $$dS_{K, n-k}\coloneqq  \sigma_{n-k}(W[h_{K}])d\s.$$
\end{prop}

\begin{proof} 
 Let $K, L\in \mathcal{K}^{\circ}_{\theta}$, we set $K_{\varepsilon}=K+_{p}\varepsilon L$ for $\varepsilon \geq 0$.  
From \eqref{p-sum}, 
\begin{eqnarray*}
   \frac{d}{d\varepsilon}\bigg|_{\varepsilon=0} h_{K_{\varepsilon}}=\frac{1}{p}h_{L}^{p}h_{K}^{1-p}.
\end{eqnarray*} 
Together with \eqref{mixed-quemassintegral} and \eqref{k-quermassintegral}, we derive
\begin{eqnarray}
    \frac{n+1-k}{p}\mathcal{V}_{p, k}(K,L)&=&\frac{d}{d\varepsilon}\bigg|_{\varepsilon=0}\mathcal{V}_{k,\theta}(K_{\varepsilon}) \notag \\
    &=&\frac{1}{p(n+1)}\left[\int_{\C_{\theta}}h_{L}^{p}h_{K}^{1-p}H_{n-k}(W[h_{K}])d\s \right.\notag \\
    &&\left.+\int_{\C_{\theta}}h_{K}H_{n-k-1}^{ij}(W[h_{K}])\left( (h_{L}^{p}h_{K}^{1-p})_{ij}+h_{L}^{p}h_{K}^{1-p}\sigma_{ij}\right)d\s \right].   \label{mixed-quer}
    \end{eqnarray}
In view of \eqref{robin}, it is readily verified that $$\n_{\mu}(h_{L}^{p}h_{K}^{1-p})=\cot\theta h_{L}^{p}h_{K}^{1-p}~~~  \text{ on } \partial\C_{\theta}.$$ From \cite[Lemma 18.30]{Ben-book}, we know that $$\sum\limits_{i=1}^{n}\n_i\left(H_{n-k-1}^{ij}(W[h_{K}])\right) =0.$$  Using \cite[Lemma 2.7]{MWWX} and integrating by parts twice, we have  
    \begin{align*}
&\int_{\C_{\theta}} h_K H_{n-k-1}^{ij}(W[h_K])
  \left((h_L^p h_K^{1-p})_{ij}
        +h_L^p h_K^{1-p}\sigma_{ij}\right)\,d\sigma\\
&\qquad=(n-k)\int_{\C_{\theta}}
  h_L^p h_K^{1-p}H_{n-k}(W[h_K])\,d\sigma.
\end{align*}
Together with \eqref{mixed-quer}, we obtain
\begin{eqnarray*}
    \mathcal{V}_{p,k}(K, L)=\frac{1}{n+1}\int_{\C_{\theta}}h_{L}^{p}h_{K}^{1-p}H_{n-k}(W[h_K])d\s.
\end{eqnarray*}
This completes the proof of \eqref{integral-pre}.
\end{proof}

In particular, when $\theta=\pi/2$, reflection across the supporting hyperplane
identifies twice the capillary mixed $L_p$ quermassintegral with
the corresponding mixed $p$-quermassintegral introduced by Lutwak \cite[Section~1]{Lut93} for the reflected convex bodies. This observation motivates the introduction of the following capillary $k$-th $L_p$ surface area measure for capillary convex bodies.
With the normalization used below, \eqref{integral-pre} becomes
\begin{eqnarray*}
\mathcal V_{p,k}(K,L)
=\frac{1}{(n+1)\binom nk}
\int_{\C_\theta}\left(\frac{h_L}{\ell}\right)^p
dS^c_{K,n-k,p}.
\end{eqnarray*}
Thus, for $p\geq1$, the measure is paired with the $p$-th power
of the capillary support function $h_L/\ell$ in the Firey variation.
For arbitrary $p\in\RR$, the following formula is an algebraic
definition for positive support functions.

\begin{defn}\label{defn:2.10}
Let $p\in \RR$, $1\le k\le n$, and let $\widehat{\Sigma}\in\mathcal K_{\theta}^{\circ}$. We define
\begin{eqnarray*}
dS_{\widehat{\Sigma},k,p}^{c} \coloneqq \ell^{p} h_{\widehat{\Sigma}}^{1-p} \, dS_{\widehat{\Sigma}, k}= \ell^{p} h_{\widehat{\Sigma}}^{1-p} \s_{k} \left(\n^2 h_{\wh\S}+h_{\wh\S} \s \right) d\s,
\end{eqnarray*}
and call $dS^{c}_{\widehat{\Sigma},k,p}$ the \emph{capillary $k$-th $L_p$ surface area measure} of the capillary convex body $\widehat{\Sigma}$.
\end{defn}

In particular, when $p=1$ and $1 \le k \le n$, $dS_{\widehat{\Sigma}, k, p}^{c}$ coincides with the \emph{capillary $k$-th area measure} introduced by the authors in \cite[Eq. (1.4)]{MWW-CM}. On the other hand, when $k=n$, $dS^c_{\wh\S,n,p}$ coincides with the \textit{capillary $L_p$ surface area measure} introduced in \cite[Definition 2.2 or Eq. (1.3)]{MWW-Lp}. Thus, the capillary \(k\)-th \(L_p\)-surface area measure provides a common generalization of both the capillary \(k\)-th area measure and the capillary \(L_p\) surface area measure.

To conclude this subsection, we formulate the \emph{capillary $L_p$ Christoffel--Minkowski problem}. The goal is to identify a capillary convex body $\widehat{\Sigma} \in \mathcal{K}^{\circ}_{\theta}$ whose capillary $k$-th $L_p$-surface area measure is prescribed on $\C_{\theta}$.

\

\noindent	\textbf{Capillary $L_p$ Christoffel--Minkowski problem:} \textit{Let $1\leq k\leq n$. Given a positive function $f\in C(\C_{\theta})$,  does  there exist a capillary convex body $\widehat{\Sigma}\in\mathcal{K}^{\circ}_{\theta}$ such that
\begin{eqnarray}\label{eqn-Lp-CM-problem}
    dS_{\widehat{\S}, k, p}^{c}= \ell^{p}f d\sigma~?
\end{eqnarray}}

From \cite[Proposition~2.6]{MWWX}, there exists a one-to-one correspondence between a strictly convex capillary hypersurface and a strictly convex function $h \in C^{2}(\C_{\theta})$ satisfying the Robin boundary condition \eqref{robin}. 
By applying arguments similar to those in \cite[Proposition~2.4]{MWW-AIM}, the capillary $L_p$ Christoffel--Minkowski problem can be reduced to the existence problem for strictly convex solutions of a fully nonlinear elliptic equation with a Robin boundary condition.  Namely, 

\begin{prop}\label{prop:cm-problem-lp} 
In the smooth category, the capillary $L_p$ Christoffel--Minkowski problem  \eqref{eqn-Lp-CM-problem} is equivalent to finding a positive function $h\in C^2(\C_\theta)$ with $W[h]>0$ satisfying the following Hessian equation with Robin boundary condition:
\begin{eqnarray*}
\left\{
\begin{array}{rcll}\vspace{1mm}
\s_k(\nabla^{2} h + h \, \sigma) &=& f \, h^{\,p-1}, & \quad \text{in } ~\C_{\theta},\\
\nabla_\mu h &=& \cot \theta \, h, & \quad \text{on }~ \partial \C_{\theta}.
\end{array}
\right.
\end{eqnarray*}
	\end{prop}

\section{A priori estimates}\label{sec-3}
In this section, we derive a priori estimates for solutions of Eq. \eqref{quo-equ} and state the main theorem as follows.
\begin{thm} \label{thm priori est}
Let $0\leq l<k\leq n$ and $\theta\in(0, \frac \pi 2)$.    For any positive smooth function $f$ on $\C_{\theta}$, the following results hold :
    \begin{enumerate}
        \item \label{thm-C2-est-p-large} If $p>k-l+1$ and $h$ is an admissible solution of Eq. \eqref{quo-equ}, then
    \begin{eqnarray}\label{thm-lower}
        \min\limits_{\C_{\theta}}h \geq c, 
    \end{eqnarray}
    and for any $\alpha \in (0, 1)$, 
    \begin{eqnarray}\label{thm-C2}
        \|h\|_{C^{4,\alpha}(\C_{\theta})}\leq C,
    \end{eqnarray}
   where the positive constants $c, C$ depend only on $n, p$ and $f$.

\item \label{thm-C2-est-p-middle} If $k-l+1< p\leq k-l+2$, then for any positive, admissible solution $h$ of Eq. \eqref{quo-equ}, the rescaled solution 
   \begin{eqnarray*}
       \widetilde{h}\coloneqq \frac{h}{\min\limits_{\C_{\theta}}h}
   \end{eqnarray*}
satisfies
\begin{eqnarray}\label{re-est}
    \|\widetilde{h}\|_{C^{4,\alpha}(\C_{\theta})}\leq C,
\end{eqnarray}
where $C>0$ depends only on $n$ and $f$, and is independent of $p$.
  \item \label{thm-C2-est-p-small}  

  If $1\leq p<k-l+1$ and $h$ is a positive, capillary even,
strictly convex solution of Eq. \eqref{quo-equ}, then
\eqref{thm-lower} and \eqref{thm-C2} hold.
For each fixed $p$, $\theta$ and $\alpha\in(0,1)$,
the constants in this part depend only on
$n,k,l,p,\theta,\alpha$, $\min\limits_{\C_\theta}f$,
and $\|f\|_{C^{4,\alpha}(\C_\theta)}$.
In particular, they are uniform along smooth families
of data with these bounds.
    \end{enumerate}
\end{thm}

To prove the main theorem, we first introduce some notation. We rewrite \eqref{quo-equ} as
\begin{eqnarray}\label{F-equ}
	 F(W)\coloneqq\left(\frac{\sigma_{k}(W)}{\sigma_{l}(W)}\right)^{\frac{1}{k-l}}=f^{\frac{1}{k-l}}h^{\frac{p-1}{k-l}}= \widehat f h^{a}~~~~~~\text{ in }~~ \C_\theta,
\end{eqnarray} 
where $\widehat{f}\coloneqq f^{\frac{1}{k-l}}$ and \begin{eqnarray}\label{eqn-defn-a}
    a\coloneqq \frac{p-1}{k-l}.
\end{eqnarray} Moreover, we denote
	\begin{eqnarray*}
		F^{ij}\coloneqq \frac{\partial F}{\partial W_{ij}},
		\quad  F^{ij,kl}\coloneqq \frac{\partial^{2}F}{\partial W_{ij}\partial W_{kl}}.
	\end{eqnarray*}
In the rest of this paper, for convenience, we adopt a local frame to express tensors and their covariant derivatives on the spherical cap $\C_\theta$. Throughout this paper, indices appearing as subscripts on tensors indicate covariant differentiation. For instance, given an orthonormal frame $\{e_i\}_{i=1}^n$ on $\C_\theta$, the notation $h_{ij}$ stands for $\nabla^2 h(e_i,e_j)$, and $W_{ijk} \coloneqq \nabla_{e_k}W_{ij}$, and so forth. We employ the Einstein summation convention: repeated indices are implicitly summed over, regardless of whether they appear as upper or lower indices. In cases where ambiguity may arise, summation will be indicated explicitly. The constant $C$ may vary from line to line, but its dependence will be specified when necessary.

\subsection{\texorpdfstring{$C^{0}$ and $C^{1}$}{} estimates}
In this subsection, when $p>k-l+1$, we derive the $C^{0}$ and $C^{1}$ estimates for admissible solutions to Eq. \eqref{quo-equ}.   Furthermore, we establish a logarithmic gradient estimate for positive, admissible solutions to Eq. \eqref{quo-equ}, this estimate plays a crucial role in dealing with the case $p=k-l+1$.  {For $1<p<k-l+1$, we establish a weighted gradient estimate
for every acute contact angle, using the argument
of \cite{hi-lpcurvature}.
The endpoint $p=1$ is treated separately in
Lemma~\ref{lem-p-one-C2} and in the proof of
Theorem~\ref{thm priori est}.}

\begin{lem}\label{lem kC0}
   Let $0\leq l< k\leq n$ and  $\theta \in (0, \pi)$. Suppose that $p \geq  k-l+1$ and that $h$ is a positive admissible solution to \eqref{quo-equ}. 
\begin{enumerate}
    \item If $\theta \in (0,\pi/2]$, then
   \begin{eqnarray}\label{c0 est}
\frac{\binom{n}{k} }{\binom{n}{l}}\frac{(1-\cos\theta)^{p+l-k-1}}{\max\limits_{\C_{\theta}}f \cdot (\sin\theta)^{2(p-1)}}\leq        h^{p+l-k-1}(\xi)
    \leq\frac{\binom{n}{k} }{\binom{n}{l}}\frac{(\sin\theta)^{2(p+l-k-1)}}{\min\limits_{\C_{\theta}}f \cdot (1-\cos\theta)^{p-1}}, ~~~\xi\in \C_\theta.
   \end{eqnarray}

\item If $\theta \in (\pi/2,\pi)$, then
    \begin{eqnarray}\label{c0 est-1}
\frac{\binom{n}{k} }{\binom{n}{l}}\frac{(\sin\theta)^{2(p+l-k-1)}}{\max\limits_{ \C_{\theta}}f \cdot (1-\cos\theta)^{p-1}} \leq h^{p+l-k-1}(\xi) \leq\frac{\binom{n}{k} }{\binom{n}{l}}\frac{(1-\cos\theta)^{p+l-k-1}}{\min\limits_{ \C_{\theta}}f \cdot (\sin\theta)^{2(p-1)}},~~~\xi\in \C_\theta.
    \end{eqnarray}
\end{enumerate}
Moreover, if $p>k-l+1$, then \eqref{c0 est} and \eqref{c0 est-1} provide uniform $C^0$ bounds for admissible solutions of Eq. \eqref{quo-equ}.  
\end{lem}

\begin{proof}
Consider the capillary support function (see, e.g.,  \cite[Eq. (1.6)]{MWW-AIM} and \cite[Eq. (2.6)]{WWX24})
    \begin{eqnarray}\label{u}
        u \coloneqq \ell^{-1} h,
    \end{eqnarray}where $\ell$ is given by \eqref{eqn-ell}. 
Then, by \eqref{quo-equ}, $u$ satisfies
\begin{eqnarray}\label{eq-monge-ampere-Lp-u}
\left\{\begin{array}{rcll}\vspace{2mm}\displaystyle
		\frac{\sigma_{k}( \ell \n^2 u+\cos \theta ( \n u \otimes e^T +e^T\otimes \n u)+  u \sigma)}{\sigma_{l}(\ell \n^2 u+\cos \theta ( \n u \otimes e^T +e^T\otimes \n u)+u\sigma)}&=&  f(u\ell)^{p-1},& \quad    \hbox{ in } ~ \C_{\theta},\\
			\n_\mu u &=& 0, &\quad  \hbox{ on } ~\p \C_\theta,
		\end{array}\right.
	\end{eqnarray}  
where $e^{T}$ denotes the tangential component of $e$ on $\mathcal{C}_{\theta}$, and we used that $\nabla_{\mu}\ell = \cot\theta\,\ell$ on $\partial \mathcal{C}_{\theta}$.

Suppose that $u$ attains its maximum at some point $\xi_{0}\in \mathcal{C}_{\theta}$. If $\xi_{0}\in \mathcal{C}_{\theta}\setminus \partial\mathcal{C}_{\theta}$, then
\begin{eqnarray}\label{deri12}
    \n u(\xi_{0})=0\quad \text{and}\quad \n^{2}u(\xi_{0})  \leq 0.
\end{eqnarray}
If $\xi_{0}\in \partial\C_{\theta}$, then the boundary condition in \eqref{eq-monge-ampere-Lp-u} implies that \eqref{deri12} still holds, see, e.g., \cite[Proof of Lemma 3.2 and Eq. (3.5)]{MWW-Weingarten}. In the following, we perform computations at $\xi_{0}$. Substituting \eqref{deri12} into the first equation in \eqref{eq-monge-ampere-Lp-u}, we obtain
\begin{eqnarray}\label{maximum-condition}
   f(u\ell)^{p-1}=\frac{\sigma_{k}( \ell \n^2 u+\cos \theta ( \n u \otimes e^T +e^T\otimes \n u)+  u \sigma)}{\sigma_{l}(\ell \n^2 u+\cos \theta ( \n u \otimes e^T +e^T\otimes \n u)+u\sigma)} \leq \frac{\binom{n}{k} }{\binom{n}{l}} u^{k-l}.
\end{eqnarray}
If $\theta \in (0,\pi/2]$, then by \eqref{maximum-condition} we obtain
\begin{eqnarray*}
    u^{p+l-1-k} (\xi_0)\leq\frac{\binom{n}{k} }{\binom{n}{l}} f^{-1}\ell^{1-p}\leq  \frac{\binom{n}{k} }{\binom{n}{l}}\frac{1}{\min\limits_{\C_{\theta}}f \cdot (1-\cos\theta)^{p-1}}.
\end{eqnarray*}
Hence, for any $\xi\in \mathcal C_{\theta}$,
\begin{eqnarray*}
    h^{p+l-k-1}(\xi)&=&\left(u\ell\right)^{p+l-k-1}(\xi)
    \leq  u^{p+l-k-1}(\xi_{0})\cdot \left(\max \limits_{\xi\in \C_{\theta}}\ell(\xi)\right)^{p+l-k-1}\notag \\
    &\leq&\frac{\binom{n}{k} }{\binom{n}{l}}\frac{(\sin\theta)^{2(p+l-k-1)}}{\min\limits_{\C_{\theta}}f \cdot (1-\cos\theta)^{p-1}}.
\end{eqnarray*}
Similarly, there holds
\begin{eqnarray*}
    h^{p+l-k-1}(\xi)\geq\frac{\binom{n}{k} }{\binom{n}{l}}\frac{(1-\cos\theta)^{p+l-k-1}}{\max\limits_{ \C_{\theta}}f \cdot (\sin\theta)^{2(p-1)}}.
\end{eqnarray*}
This proves \eqref{c0 est}. The proof of \eqref{c0 est-1} is analogous.
The final assertion follows immediately from \eqref{c0 est} and \eqref{c0 est-1}.
This completes the proof.
\end{proof}

Next, we establish the $C^{1}$ estimates for the solution to Eq. \eqref{F-equ} when $p>k-l+1$. The construction of the auxiliary function is inspired by the ideas of Ma--Xu \cite{MX} and Ma--Qiu--Xu \cite{MXQ} for studying the Neumann boundary problems. To begin with, we define the function
 \begin{eqnarray}\label{d-func}
 d(\xi)\coloneqq {\rm{dist}}(\xi, \partial\C_{\theta}).
 \end{eqnarray}
 The function $d(\xi)$ is smooth near $\partial\C_{\theta}$, we 
  extend it to a smooth function on   $\C_{\theta}$ and still denote it as $d$. The function $d$ satisfies  
  \begin{eqnarray}\label{eqn:bdry-distance}
    d\big |_{\partial\C_{\theta}}=0, \quad {\rm{and}} \quad \n d\big |_{\partial\C_{\theta}}=-\mu,  
  \end{eqnarray} and 
\begin{eqnarray}\label{2-big}
   \frac{1}{2} \leq 1+\cot\theta d\leq 2.
    \end{eqnarray}

\begin{lem}\label{lem c1}
  Let $0\leq l<k\leq n$  and $\theta\in (0, {\pi})$. Suppose that $p>k-l+1$ and  $h$ is a positive admissible solution to Eq. \eqref{quo-equ}, then 
    \begin{eqnarray}\label{eqn:C1-est}
       |\n h|\leq C, 
    \end{eqnarray}
 where $C>0$ depends on $n, k, l, p, \min\limits_{\C_{\theta}}f$ and $\|f\|_{C^{1}(\C_{\theta})}$.
 \end{lem}

\begin{proof}
Consider the function
\begin{eqnarray*}
\Phi(x)\coloneqq \frac{1}{2}\log\left(1+|\n \widehat{h}|^{2}\right)-\frac{1}{2}\log \left(S-\widehat{h}-d\right)-a_{0}d,
\end{eqnarray*}
where  
$$\widehat{h}\coloneqq (1+\cot\theta d)h,$$ and  we choose the constants \begin{eqnarray}\label{eqn:choice-S}
    S\coloneqq 2+2\|h\|_{C^{0}(\C_{\theta})}+\|d\|_{C^{0}(\C_{\theta})},~~~~~a_{0}\coloneqq \min\{\cot\theta, 0\}.
\end{eqnarray} 
It follows from \eqref{eqn:choice-S} that 
\begin{eqnarray}\label{eqn:log-C1-est-two-sided-bound}
1\leq S-\widehat{h}-d\leq 2S\quad {\text{in}}\quad \C_{\theta}.
\end{eqnarray}
Moreover, by \eqref{robin} and \eqref{eqn:bdry-distance}, we have
\begin{eqnarray}\label{eqn-h-normal-derivative}
    \n_{\mu}\widehat{h}=\cot\theta h \n_{\mu}d +\n_{\mu}h=0\quad \text{on}\quad \p\C_\theta.
\end{eqnarray}
Assume $\Phi$ attains its maximum value at some point, say $\xi_{0}\in \C_{\theta}$. In view of the position of $\xi_{0}$, we divide the proof into the following two cases. For the remainder of the argument, all calculations are performed at $\xi_{0}$. We may further assume that \(|\nabla \wh h|(\xi_{0})\ge 1\), since otherwise the desired estimate is immediate by combining with Lemma \ref{lem kC0}.

\

\noindent\textbf{Case 1} $\xi_{0}\in \partial \C_{\theta}$.   
Choose a local orthonormal frame $\{e_{i}\}_{i=1}^{n}$ around $\xi_{0}$ such that $e_{n}=\mu$ and $\n \wh h= |\n \wh h| e_1$ at $\xi_{0}$.
By \eqref{eqn-h-normal-derivative} and the Gauss--Weingarten formula for $\p \C_\theta \subset \C_\theta$, 
\begin{eqnarray*}
\widehat{h}_{n1}=\n_{e_{1}}(\n_{e_{n}}\widehat h)-\s(\n_{e_{1}}e_{n}, \n \widehat h)=-\cot\theta \widehat{h}_{1},
\end{eqnarray*}
together with \eqref{eqn:choice-S} and \eqref{eqn-h-normal-derivative}, we derive
\begin{eqnarray*}
0\leq \n_{n}\Phi 
&=&\sum_{k=1}^n\frac{\widehat{h}_{k}\widehat{h}_{kn}}{1+|\n \widehat{h}|^{2}}+\frac{\widehat{h}_{n}+d_{n}}{2(S-\widehat{h}-d)}+a_{0}
=\frac{\widehat{h}_{1}\widehat{h}_{1n}}{1+|\n \widehat{h}|^{2}}-\frac{1}{2(S-\widehat{h}-d)}+a_{0}\notag \\
&=&-\left(\frac{\cot\theta \widehat h_{1}^{2}}{1+|\n \widehat{h}|^{2}}+\frac{1}{2(S-\widehat h-d)}
\right)+a_{0}<0,
\end{eqnarray*}
which is a contradiction.

\

\noindent\textbf{Case 2}  $\xi_{0}\in \C_{\theta}\setminus \partial \C_{\theta}$. Choose a local orthonormal frame $\{e_{i}\}_{i=1}^{n}$ around $\xi_0$ such that 
$$W=\n^2 h +h\s$$ is diagonal at $\xi_{0}$. Consequently, by Proposition~\ref{prop2.1}, $$ F^{ij}= \frac{\p F}{\p W_{ij}} $$ is also diagonal at \(\xi_{0}\). Since \(\Phi\) attains its maximum at \(\xi_{0}\), we have
\begin{eqnarray}\label{one deri cond}
0=\n_{e_{i}}\Phi=\frac{\widehat{h}_k\widehat{h}_{ki}}{1+|\n \widehat{h}|^{2}}+\frac{\widehat{h}_{i}+d_{i}}{2(S-\widehat{h}-d)}-a_{0}d_{i}.
\end{eqnarray}
Differentiating once more and using \eqref{one deri cond}, we obtain
\begin{eqnarray}\label{c1-contract cond}
0\geq F^{ij}\n_{ij}\Phi
&=& \frac{F^{ij}\widehat{h}_{ki}\widehat{h}_{kj}+F^{ij}\widehat{h}_{k}\widehat{h}_{kij}}{1+|\nabla \widehat{h}|^{2}}+\frac{F^{ij}(\widehat{h}_{ij}+d_{ij})}{2(S-\widehat{h}-d)}+\frac{2a_{0}F^{ij}(\widehat{h}_{i}+d_{i})d_{j}}{S-\widehat{h}-d} \notag  \\
&&-2a_{0}^{2}F^{ij}d_{i}d_{j}-a_{0}F^{ij}d_{ij}.
\end{eqnarray}
Since $F(A)$ is homogeneous of degree one in \eqref{F-equ}, we have 
\begin{eqnarray}
    F^{ij}\wh h_{ij} &= & (1+\cot \theta d) \wh f h^a-(1+\cot \theta d) h \sum_{i=1}^n F^{ii}+2\cot\theta F^{ij} d_i h_j+h \cot\theta F^{ij}d_{ij}  \notag 
    \\& \geq & - C \left(\|h\|_{C^0(\C_\theta)} +\|\n h\|_{C^0(\C_\theta)} \right) \sum_{i=1}^n F^{ii}, \label{eqn:Fij-h-ij-C1-log-est}
\end{eqnarray}where $C>0$ depends only on $n$.
Taking the derivative of Eq. \eqref{F-equ} in the $e_k$ direction, we obtain
\begin{eqnarray}\label{one deri condition-1}
F^{ij}h_{ijk}= (\widehat{f}h^{a})_k-h_k\sum\limits_{i=1}^{n}F^{ii}.
\end{eqnarray}
On the spherical cap $\C_\theta$, the commutation formula gives:
\begin{eqnarray}\label{third comm}
h_{kij}=h_{ijk}+h_k\d_{ij}-h_j\d_{ki}.
\end{eqnarray}
Combining \eqref{one deri condition-1} and \eqref{third comm}, and applying the Cauchy-Schwarz inequality,  we derive
\begin{eqnarray}\label{combi-2}
F^{ij}\widehat{h}_k\widehat{h}_{kij}&=&F^{ii}\widehat{h}_k \left(\widehat{h}_{iik}+\widehat{h}_k-\widehat{h}_{i}\delta_{ki} \right)\geq F^{ii}\widehat{h}_k\widehat{h}_{iik}\notag \\  
   &\geq &2\cot\theta F^{ii}d_{i}h_{ik}\widehat{h}_k+\cot\theta F^{ii}d_{k}h_{ii}\widehat{h}_k-C \left(1+|\n \widehat h|^{2}\right)\sum\limits_{i=1}^{n}F^{ii}\notag 
   \\
   &\geq& -\frac{1}{4}(1+\cot\theta d)^{2}F^{ii}W_{ii}^{2}-C\left(1+|\n \widehat h|^{2} \right)\sum\limits_{i=1}^{n}F^{ii},
\end{eqnarray}
and 
\begin{eqnarray}\label{combi-1}
    F^{ij}\widehat{h}_{ki}\widehat{h}_{kj}
    &\geq& (1+\cot\theta d)^{2}F^{ii}W_{ii}^{2}-C|\n \widehat{h}|F^{ii}|W_{ii}|-C \left(1+|\n \widehat{h}|^{2} \right)\sum\limits_{i=1}^{n}F^{ii}\notag \\ 
    &\geq &\frac{1}{2}(1+\cot\theta d)^{2}F^{ii}W_{ii}^{2}-C \left(1+|\nabla \widehat{h}|^{2} \right ) \sum\limits_{i=1}^{n}F^{ii},
\end{eqnarray}
where $C>0$ depends only on $n, k, l , \min\limits_{\C_{\theta}}f, \|f\|_{C^{1}(\C_{\theta})}$ and $\|h\|_{C^{0}(\C_{\theta})}$.

Substituting \eqref{eqn:Fij-h-ij-C1-log-est}, \eqref{combi-2}  and \eqref{combi-1}  into \eqref{c1-contract cond}, using \eqref{eqn:log-C1-est-two-sided-bound}, we obtain
\begin{eqnarray}
\begin{aligned}\label{key ineq}
    0\geq& \frac{1}{1+|\n \widehat{h}|^{2}}\left[\frac{1}{4}\left(1+\cot\theta d\right)^{2}F^{ii}W_{ii}^{2}
    -C \left(1+|\n \widehat{h}|^{2} \right)\sum\limits_{i=1}^{n}F^{ii}\right]\\
    &-C \left(|\n\widehat{h}|+1 \right )\left(\sum\limits_{i=1}^{n}F^{ii}+1\right).
    \end{aligned}
\end{eqnarray}
 Define $$\mathcal T\coloneqq \left\{~i~~\Big|~ 1\leq i\leq n, ~|\widehat{h}_{i}|\geq \frac{1}{\sqrt{n}}|\n \widehat{h}| \right\}.$$ 
Then $\mathcal T\neq\emptyset$.

Fix $i_0\in\mathcal T$.
By \eqref{one deri cond},
\begin{align}
\widehat{h}_{i_{0}i_{0}}&=-\frac{1+|\n \widehat{h}|^{2}}{\widehat{h}_{i_{0}}}\left(\frac{\widehat{h}_{i_{0}}+d_{i_{0}}}{2(S-\widehat{h}-d)} -a_{0}d_{i_{0}} \right)-\frac{1}{\widehat{h}_{i_{0}}}\sum\limits_{i\neq i_{0}} \widehat{h}_{i}\widehat{h}_{ii_{0}} \notag\\
&=-(1+|\n\widehat{ h}|^{2})\left[ \frac{1}{2(S-\widehat{h}-d)} \left(1+\frac{d_{i_{0}}}{\widehat{h}_{i_{0}}} \right)-\frac{a_{0}d_{i_{0}}}{\widehat{h}_{i_{0}}} \right] \notag \\& \quad -\frac{\cot\theta}{\widehat{h}_{i_{0}}}\sum\limits_{i\neq i_{0}}\widehat{h}_{i}(h d_{ii_{0}}+ d_{i}h_{i_{0}}+ d_{i_{0}}h_{i}).\quad \quad \quad \label{hii}
\end{align}
Observe that 
\begin{eqnarray}
\left|    \frac{\cot\theta}{\widehat{h}_{i_{0}}} \sum\limits_{i\neq i_{0}}\widehat{h}_{i} \left(h d_{ii_{0}}+ d_{i}h_{i_{0}}+ d_{i_{0}}h_{i} \right) \right| &\leq & \left|    \frac{\cot\theta \sqrt{n}}{|\n \wh h|}\sum\limits_{i\neq i_{0}}\widehat{h}_{i} \left(h d_{ii_{0}}+ d_{i}h_{i_{0}}+ d_{i_{0}}h_{i} \right) \right|  \notag
\\& \leq & C_1 \left(1+|\n \wh h| \right), \label{eqn:log-C1-term2}
\end{eqnarray} where $C_1>0$ depends only on $n$ and $\|h\|_{C^0(\C_\theta)}$.
We may assume that \begin{eqnarray}\label{eqn:log-C1-large-nable-h}
    |\n \wh h|\geq \max\left\{6\sqrt{n}, ~ 12 S|a_0|\sqrt{n}, ~  24S \left(1+\|h\|_{C^0(\C_\theta)} \right)+C_1\right\},
\end{eqnarray} since otherwise the desired gradient estimate \eqref{eqn:C1-est} is immediate. Together with \eqref{eqn:log-C1-est-two-sided-bound}, it yields
\begin{eqnarray}
    \frac{1}{2(S-\widehat{h}-d)} \left(1+\frac{d_{i_{0}}}{\widehat{h}_{i_{0}}} \right)-\frac{a_{0}d_{i_{0}}}{\widehat{h}_{i_{0}}} &\geq& \frac 1 {4S} \left(1-\frac {\sqrt{n}|\n d|}{|\n \wh h|} \right)- \frac{|a_0|\sqrt{n}|\n d|}{|\n \wh h|} \notag
    \\&\geq & \frac 1 {4S} \left(1-\frac 1 6\right)- \frac{1}{12S}=\frac 1 {8S}.\label{eqn:log-C1-term1}
\end{eqnarray}
Substituting \eqref{eqn:log-C1-term2} and \eqref{eqn:log-C1-term1} into \eqref{hii}, and using \eqref{eqn:log-C1-large-nable-h}, we infer that
    \begin{eqnarray}\label{Wii}
        W_{i_{0}i_{0}}\coloneqq h_{i_{0}i_{0}}+h\leq -\frac{|\n\widehat{h}|^{2}+1}{16S}.
    \end{eqnarray} 
Without loss of generality, assume that $$W_{11}\leq W_{22}\leq \cdots\leq W_{nn}.$$ From \eqref{Wii}, 
\begin{eqnarray}\label{w11}
    W_{11}\leq W_{i_{0}i_{0}}\leq  -\frac{|\n \widehat h|^{2} +1 }{16S}.\end{eqnarray}
Moreover, Proposition \ref{pro-2.3} (3), (4) yields
\begin{eqnarray}\label{key ine-2}
    F^{11}\geq\frac{1}{n}\sum\limits_{i=1}^{n}F^{ii}\geq \frac{c_{n,k,l}}{n}.
\end{eqnarray}
Substituting \eqref{w11} and \eqref{key ine-2} into \eqref{key ineq},  using \eqref{2-big},  we obtain 
\begin{eqnarray*}
    0&\geq &\frac{1}{4}\frac{(1+\cot\theta d)^{2}}{1+|\n \widehat h|^{2}}F^{11}W_{11}^{2}-C \left(|\n \widehat h|+1 \right)\left(\sum\limits_{i=1}^{n}F^{ii}+1\right)\notag \\
    &\geq &\frac{(1+\cot\theta d)^{2}}{4096S^{2}}  \left(1+|\n \widehat h|^{2}\right) F^{11}-C \left(|\n \widehat h|+1\right)\left(\sum\limits_{i=1}^{n}F^{ii}+1\right)\\&\geq & 
 \frac{1}{4096nS^{2}} \left(1+|\n \widehat h|^{2}\right)  \sum_{i=1}^n F^{ii}  -C \left(|\n \widehat h|+1 \right)\left(\sum\limits_{i=1}^{n}F^{ii}+1\right).
    \end{eqnarray*}
Hence, by Proposition \ref{pro-2.3} (4), it follows that
$$\max\limits_{\C_{\theta}}|\n \wh h|\leq C,$$ 
and therefore
$$\max\limits_{\C_{\theta}}|\n h|\leq C,$$ where, in view of Lemma \ref{lem kC0}, the constant \(C>0\) depends only on $n, k, l, p, \min\limits_{\C_{\theta}}f$ and $\|f\|_{C^{1}(\C_{\theta})}$.
This completes the proof.
 \end{proof}

Next, we establish a logarithmic gradient estimate for positive, admissible solutions to Eq. \eqref{quo-equ}. This estimate plays an important role in treating the case $p=k-l+1$. Since  Eq. \eqref{quo-equ} is dilation invariant when $p=k-l+1$,  we cannot obtain the $C^{0}$ estimates as in the case $p>k-l+1$. We adopt the strategy presented in Guan--Lin  \cite[Section 3]{GL2000} and Hu--Ma--Shen \cite[Section 3]{HMS2004}, which uses an approximating argument. This logarithmic gradient estimate will be used to derive a Harnack-type estimate for the approximating solutions when $p\in(k-l+1,k-l+2]$. It extends our previous result \cite[Lemma~3.5]{MWW-Lp} for the capillary Minkowski problem. The choice of the auxiliary function is motivated by the ideas developed in \cite{MX, MQ} for Neumann boundary problems. 

\begin{lem}\label{lem kC1}
  Let $0\leq l<k\leq n$ and $\theta\in (0,  {\pi}/{ 2})$.  Suppose that $p>k-l+1$ and $h$ is a positive admissible solution of Eq. \eqref{quo-equ}, then there exists a  positive constant $C$ depending on $n, k, l,p,  \min\limits_{\C_{\theta}}f$ and $\|f\|_{C^{1}(\C_{\theta})}$, such that 
  \begin{eqnarray}\label{gradient logh}
\max\limits_{\C_{\theta}}|\n\log h|\leq C. 
  \end{eqnarray}
  For each fixed $\theta\in(0,\pi/2)$, the constant $C$ can be chosen
independently of $p$ when $k-l+1<p\le k-l+2$.
  \end{lem}

\begin{proof}
Define $$v\coloneqq \log h.$$ From Eq. \eqref{F-equ}, we deduce that $v$ satisfies 
\begin{eqnarray}\label{G-equ}
	 G(B)\coloneqq\left(\frac{\sigma_{k}(B)}{\sigma_{l}(B)}\right)^{\frac{1}{k-l}}=\wt{f}~~~~~~\text{ in }~~ \C_\theta,
\end{eqnarray} 
where $B \coloneqq \n^2 v+\n v\otimes \n v+\sigma$ and $\wt{f}\coloneqq f^{\frac{1}{k-l}}e^{p_{0}v}$ with $p_{0}\coloneqq \frac{p-(k-l+1)}{k-l}$. Moreover, by \eqref{robin}, the boundary condition becomes 
\begin{eqnarray} \label{v-bry}
    \n_{\mu} v=\cot\theta ~~~~~~~~\text{ on }\p \C_\theta .
\end{eqnarray}

 Define  
\begin{eqnarray*}
\Psi\coloneqq  |\n \omega|^{2}.
\end{eqnarray*}
where  \begin{eqnarray}\label{t-l}
\omega\coloneqq v-\cot\theta \tilde{\ell},~~~\text{ and }~~~     \tilde{\ell}\coloneqq  \frac{\ell}{\sin\theta \cos\theta}. 
    \end{eqnarray}
A direct computation shows that
\begin{eqnarray}\label{ell-bry}
    \tilde{\ell}|_{\partial\C_{\theta}}=\tan\theta,\quad\quad \n\tilde{\ell}|_{\partial\C_{\theta}} =\mu,
\end{eqnarray}
and 
\begin{eqnarray}\label{ell-positive}
    \n^{2}\tilde{\ell}=\frac{\n^{2}\ell}{\sin\theta\cos\theta}=\frac{(1-\ell)}{\sin\theta\cos\theta}\sigma. 
\end{eqnarray}

Suppose that $\Psi$ attains its maximum at some point $\xi_0 \in \C_\theta$. We divide the proof into two cases: either  $\xi_0\in \p \C_\theta$ or $\xi_0\in \C_\theta\setminus \p \C_\theta$.

 \

\noindent\textbf{Case 1}~ $\xi_{0}\in \partial \C_{\theta}$. Let $\{e_{i}\}_{i=1}^{n}$ be an orthonormal frame around $\xi_0$ such that $e_{n}=\mu$. By \eqref{v-bry} and \eqref{ell-bry}, we have \begin{eqnarray}\label{eqn:omega-n=0}
 \o_n=  \n_\mu \omega=0~~~~~~~\text{ on } ~~\p \C_\theta. 
\end{eqnarray}
By the Gauss-Weingarten equation for $\partial\C_{\theta}\subset \C_{\theta}$,  
\begin{eqnarray}\label{w-2}
    \omega_{\alpha n}=\n_{\alpha}(\omega_{n})-\sigma(\n \omega, \n_{e_{\alpha}}e_{n})=-\cot\theta \omega_{\alpha}.
\end{eqnarray}
 Since \(\Psi\) attains its maximum at \(\xi_{0}\), it follows from \eqref{eqn:omega-n=0} and \eqref{w-2} that 
\begin{eqnarray*} 
    0\leq \n_{n}\Psi &=& 2\sum_{i=1}^n \omega_i\omega_{in}=2\sum\limits_{\alpha=1}^{n-1}\omega_{\alpha}\omega_{\alpha n}
    \\&=& -2\cot\theta \sum\limits_{\alpha=1}^{n-1}\omega_{\alpha}^{2}.
\end{eqnarray*}
Since $\theta\in(0, \frac \pi 2)$, we obtain
\begin{eqnarray*}
  \sum\limits_{\alpha=1}^{n-1}\omega_{\alpha}^{2}\leq 0.
\end{eqnarray*}
Therefore, $$|\nabla \omega|=0,$$ at \(\xi_{0}\), where we have used \eqref{eqn:omega-n=0} again. Consequently,
\begin{eqnarray*}
    |\n \log h|\leq \frac 1{\sin^2\theta} \sup_{\C_\theta} |\n   \ell| 
\end{eqnarray*}
and hence \eqref{gradient logh} holds.

 \

\noindent \textbf{Case 2}~ $\xi_{0}\in \C_{\theta}\setminus \partial\C_{\theta}$. Choose an orthonormal frame  $\{e_{i}\}_{i=1}^{n}$ around $\xi_{0}$, such that  $$B_{ij}\coloneqq v_{ij}+v_iv_j+\s_{ij} $$ is diagonal at $\xi_0$. Consequently, $$G^{ij}\coloneqq \frac{\partial G(B)}{\partial B_{ij}}$$ is also diagonal at $\xi_0$. 

We may assume that \(|\nabla v|\) is sufficiently large, which implies that \(|\nabla\omega|\) is sufficiently large as well; otherwise, \eqref{gradient logh} follows immediately. Since \(\Psi\) attains its maximum at the interior point \(\xi_{0}\), we have
\begin{eqnarray}\label{Psi-1}
 0= ( \log\Psi)_{i} =\frac{2\omega_{s}\omega_{si}}{ |\n \omega|^{2}},  
\end{eqnarray}
and using \eqref{Psi-1}, \begin{eqnarray} \label{Psi-two}
0\ge G^{ij} (\log\Psi)_{ij}&=&\frac{2}{|\n\omega|^{2}} G^{ij}\left(\omega_{s}\omega_{sij}+\omega_{si}\omega_{sj}\right)\notag \\
&=&\frac{2\sum\limits_{s=1}^{n}G^{ii}\omega_{si}^{2}}{|\nabla \omega|^{2}}+\frac{2G^{ij}\omega_{s}(v_{sij}-\cot\theta \tilde{\ell}_{sij})}{|\nabla \omega|^{2}}.
\end{eqnarray}
For the standard metric on the spherical cap $\C_\theta$, the commutation formula gives
		\begin{eqnarray}\label{3-com}
		v_{sij}=v_{ijs}+v_s\d_{ij}-v_j\d_{si}. 
		\end{eqnarray}
Differentiating Eq. \eqref{G-equ}, using the fact that $p_{0}\geq 0$ and \eqref{t-l}, we have 
  \begin{eqnarray}
      G^{ij} \omega_{s}\left(v_{ijs}+2v_{si}v_{j}\right)&=&\wt{f}_{s}\omega_{s}=\wt{f}\left(\frac{1}{k-l}\frac{f_{s}}{f}+p_{0}v_{s}\right)\omega_{s} \notag \\
      &\geq& \wt{f}\left(\frac{1}{k-l}\frac{f_{s}}{f}+p_{0}\cot\theta \tilde{\ell}_{s}\right)\omega_{s}.   \label{tf1-times}
  \end{eqnarray}
By Lemma \ref{lem kC0} (1), we see 
\begin{eqnarray}\label{tf-1}
    \wt{f}=\left( fh^{p+l-k-1}\right)^{\frac{1}{k-l}}\leq\left(\frac{\binom{n}{k} \max\limits_{\C_{\theta}} f }{\binom{n}{l} \min\limits_{\C_{\theta}} f}\frac{  (\sin\theta)^{2(p+l-k-1)}}{ (1-\cos\theta)^{p-1}}\right)^{\frac{1}{k-l}}.
\end{eqnarray}
Inserting  \eqref{tf-1} into \eqref{tf1-times}, we obtain
\begin{eqnarray}
      G^{ij} \omega_{s}\left(v_{ijs}+2v_{si}v_{j}\right)\geq -C|\n \omega|,   \label{1-times}
  \end{eqnarray}
where $C>0$ depends only on $n,k,l, p, \min\limits_{\C_{\theta}}f$ and $\|f\|_{C^1(\C_{\theta})}$. 

Moreover, if $$k-l+1< p\leq k-l+2,$$ then $0<p_{0}\leq\frac{1}{k-l}$ and the estimate \eqref{tf-1} improves to
\begin{eqnarray}\label{tf-2}
    \wt{f} \leq\left( \frac{\binom{n}{k} \max\limits_{\C_{\theta}} f }{(1-\cos\theta)^{k-l+1} \binom{n}{l} \min\limits_{\C_{\theta}} f} \right)^{\frac{1}{k-l}}.
\end{eqnarray}
Consequently, the constant \(C\) in \eqref{1-times} can be chosen independently of \(p\).

  Inserting \eqref{3-com} and \eqref{1-times} into \eqref{Psi-two}, we obtain 
  \begin{eqnarray}\label{Psi-3}
      0&\geq&\frac{2}{|\n \omega|^{2}}\left(\sum\limits_{s=1}^{n}G^{ii}\omega_{si}^{2} -2G^{ii}\omega_{s}v_{si}v_{i} \right)+\frac{2}{|\n \omega|^{2}} \left(
   |\n v |^{2}\sum\limits_{i=1}^{n}G^{ii}
  -G^{ii}v_{i}^{2} \right) \notag \\
  &&-\frac{C}{|\n \omega|}\left(1+\sum\limits_{i=1}^{n}G^{ii}\right).
  \end{eqnarray}
  
Set $$\mathcal{J}\coloneqq 2\sum\limits_{s=1}^{n}G^{ii}\omega_{si}^{2}-4G^{ii}\omega_s v_{si}v_{i}.$$  A direct calculation yields
\begin{eqnarray}\label{J}
    \mathcal{J}&=&2\sum\limits_{s=1}^{n}G^{ii}(B_{si}-v_{s}v_{i}-\delta_{si}-\cot\theta \t{\ell}_{si})^{2}-4G^{ii}v_{i}\omega_{s}(B_{si}-v_{i}v_{s}-\delta_{si}) \notag \\
    &\geq &-4G^{ii}B_{ii}v_{i}\omega_{i}+4G^{ii}v_{i}^{2}v_{s}\omega_{s}+4G^{ii}v_{i}\omega_{i}+2G^{ii}B_{ii}^{2}-4G^{ii}B_{ii}v_{i}^{2} \notag  \\
    &&+2|\n v|^{2}G^{ii}v_{i}^{2}+4G^{ii}v_{i}^{2} -4\cot\theta G^{ii}B_{ii}\t{\ell}_{ii}+4\cot\theta G^{ii}\t{\ell}_{si}v_{s}v_{i}\notag  \\
   && -C \left(1+\sum\limits_{i=1}^{n}G^{ii}\right).
\end{eqnarray}
 It follows from \eqref{Psi-1} that
\begin{eqnarray*}
    2\omega_{s}\left(B_{si}-v_{s}v_{i}-\delta_{si}-\cot\theta \t{\ell}_{si}\right)=2\omega_{s}\omega_{si}=0, \quad \forall~ 1\leq i\leq n.
\end{eqnarray*}
Therefore,
\begin{eqnarray*}
    \omega_{i}B_{ii}=\cot\theta \t{\ell}_{si}\omega_{s}+\omega_{i}+\omega_{s}v_{s}v_{i}, \quad \forall~ 1\leq i\leq n. 
\end{eqnarray*}
Consequently,
\begin{eqnarray}
    -G^{ii}v_{i}B_{ii}\omega_i= -\cot\theta G^{ii}\t{\ell}_{si}\omega_{s}v_{i}-G^{ii}v_{i}\omega_{i}-G^{ii}\omega_{s}v_{s}v_{i}^2. \label{logg-keyine}
\end{eqnarray}
Substituting \eqref{logg-keyine} into \eqref{J}, we obtain 
\begin{eqnarray}\label{J-2}
   \mathcal{J}&\geq & 2G^{ii}B_{ii}^{2}-4G^{ii}B_{ii}v_{i}^{2}+(2|\n v|^{2}+4)G^{ii}v_{i}^{2}-
4\cot\theta G^{ii}B_{ii}\t{\ell}_{ii}\notag \\
&&+4\cot\theta G^{ii}\t{\ell}_{si}v_{s}v_{i}-C(1+|\n v|) \left(1+\sum\limits_{i=1}^{n} G^{ii} \right).
   \end{eqnarray}
By \eqref{ell-positive}, we have
\begin{eqnarray}\label{J-4}
    4\cot\theta G^{ii}\t{\ell}_{si}v_{s}v_{i}=\frac{4(1-\ell)}{\sin^{2}\theta}G^{ii}v_{i}^{2}\geq 0.
\end{eqnarray}
Moreover, using \eqref{tf-1} and \eqref{tf-2} once again, we obtain
\begin{eqnarray}\label{J-5}
    4\cot\theta G^{ii}B_{ii}\t{\ell}_{ii}=\frac{4(1-\ell)}{\sin^{2}\theta}G^{ii}B_{ii}=\frac{4(1-\ell)\t{f}}{\sin^{2}\theta}\leq C,
\end{eqnarray}
where $C>0$ depends only on $n, k,l, p, \min\limits_{\C_{\theta}}f$ and $\max\limits_{\C_{\theta}}f$. Moreover, when
\(k-l+1<p\leq k-l+2\), the constant \(C\) can be chosen independently of \(p\).

Substituting \eqref{J-4}, \eqref{J-5} into \eqref{J-2},  we derive that
\begin{eqnarray}\label{J-3}
    \mathcal{J}&\geq&   2G^{ii}B_{ii}^{2}-4G^{ii}B_{ii}v_{i}^{2}+(2|\n v|^{2}+4)G^{ii}v_{i}^{2}-C(1+|\n v|) \left(1+\sum\limits_{i=1}^{n} G^{ii} \right)\notag \\
    &= &2G^{ii} (B_{ii}-v_i^2)^2+2 G^{ii}v_i^2(|\nabla v|^2-v_i^2+2) -C(1+|\n v|) \left(1+\sum\limits_{i=1}^{n} G^{ii} \right)
    \notag \\
    &\geq & 4G^{ii}v_{i}^{2}-C(1+|\n v|) \left(1+\sum\limits_{i=1}^{n} G^{ii} \right).
    \end{eqnarray}
Inserting \eqref{J-3} into \eqref{Psi-3}, we obtain 
\begin{eqnarray*}
    0 \geq 2 \sum\limits_{i=1}^{n}G^{ii}-\frac{C}{|\n \omega|}\left(1+\sum\limits_{i=1}^{n}G^{ii}\right),
\end{eqnarray*}
together with Proposition \ref{pro-2.3} (4), this implies that $$|\n \omega|\leq C.$$ 
Consequently, we obtain the desired estimate \eqref{gradient logh}. Furthermore, if $k-l+1< p\leq k-l+2$, the constant $C$ in \eqref{gradient logh} is independent of $p$.
This completes the proof of Lemma~\ref{lem kC1}.
\end{proof}

 In the following, we establish a weighted gradient estimate for solutions to Eq. \eqref{F-equ} when $1<p<k-l+1$. Such estimates play an important role in deriving a priori bounds, particularly the \(C^{0}\) estimate for solutions to the \(L_{p}\) Christoffel--Minkowski problem. A weighted gradient estimate of the form $\frac{|\n h|^2}{h^\gamma}$ was previously employed by Guan--Lin \cite{GL2000} for $L_p$ Minkowski problem, and was subsequently further developed in \cite{GX, Huang-Lu2013} and, in particular, in Guan \cite[Section 2]{Guan2023} for $L_p$ Christoffel--Minkowski problem; see also the later works \cite{HI, MWW-Quotient}. In the capillary setting considered here, we employ the auxiliary function $P\coloneqq \ell^\a \frac{|\n u|^2}{h^\gamma}$, 
 which was introduced in \cite{hi-lpcurvature} in the study of Eq. \eqref{quo-equ} for the case $k=n$ and $\theta\in(0,\pi/2)$. For the reader's convenience, we provide a proof here.   

\begin{lem}\label{lem weighted gra}
Let $0\leq l<k\leq n$, $1< p<k-l+1$ and $\theta\in (0, \frac \pi 2)$. Suppose that $h$ is a positive and strictly convex solution to Eq. \eqref{quo-equ}. Then there exists a  positive constant $C$ depending on $n, k, l, \min\limits_{\C_{\theta}}f$ and $\|f\|_{C^{1} (\C_{\theta})}$ and a  constant $\gamma \in \left(0,   \frac{2(p-1)}{k-l}    \right)$, such that 
  \begin{eqnarray}\label{weighed gredient est}
\frac{|\n h|^{2}}{h^{\gamma}}\leq C  \max\limits_{\C_{\theta}}h^{2-\gamma}.
  \end{eqnarray} 
  \end{lem}
 \begin{proof}
 By \eqref{u}, it suffices to prove that
      \begin{eqnarray}\label{u-weight-est}
       P\coloneqq   \frac{ \ell ^{\a} |\n u|^{2}}{h^{\gamma}}\leq  N \max_{\C_\theta} \ell^{\a-\g} \max\limits_{\C_{\theta}}u^{2-\gamma}
      \end{eqnarray}
for some positive constant $N$ and $\a\in \RR$. We proceed with the proof of \eqref{u-weight-est} by contradiction. Suppose that \eqref{u-weight-est} does not hold for some sufficiently large constant $N$. Assume further that $P$ attains its maximum value at some point $\xi_{0}\in \C_{\theta}$, 
 \begin{eqnarray}\label{assum}
     \frac{ \ell^{\a } |\n u|^{2}}{h^{\gamma}}(\xi_{0})=\max\limits_{\xi\in \C_{\theta}}P >  N \max_{\C_\theta} \ell^{\a-\g}   
   \max\limits_{\C_{\theta}}u^{2-\gamma},
 \end{eqnarray}
 it follows from \eqref{assum} that 
\begin{eqnarray}\label{N-big}
    \frac{|\n u|^{2}}{u^{2}}(\xi_{0})> N.
\end{eqnarray}
   Next, we divide the proof into two cases according to whether $\xi_0$ is an interior point or not. In the following, all the calculations are performed at the point $\xi_{0}$.

  \

 \noindent \textbf{Case~1}  $\xi_{0}\in \partial\C_{\theta}$. Let $\{e_{i}\}_{i=1}^{n}$ be an orthogonal frame around $\xi_{0}$ such that $e_{n}=\mu$. Since $\n_{\mu}u=0,$ then for all $1\leq \b\leq n-1$,  from the Gauss-Weingarten equation of $\partial\C_{\theta}\subset \C_{\theta}$,  
  \begin{eqnarray}\label{weight-u-mixed}
      \n^{2}u(e_{\b}, e_{n})=\n_{e_{\b}}(u_{n})-\n u(\n_{e_{\b}}e_{n})=-\cot\theta u_{\b}.
  \end{eqnarray}
By the maximal value condition, \eqref{robin}, \eqref{eq-monge-ampere-Lp-u} and \eqref{weight-u-mixed},  
\begin{eqnarray*}
   0&\leq &   \n_{\mu} \log P= \sum_{k=1}^n \frac{2u_ku_{kn}}{|\n u|^2} +\frac{\a \n_\mu \ell}{\ell}-\frac{\g \n_\mu h}{h}  \notag \\
  & \leq &  (   \a -2-\g) \cot\theta  
  \\& <& 0 , \notag 
\end{eqnarray*}where the last line follows by our choice of $\alpha= 2$ and $\theta<\frac\pi 2$. 
This yields a contradiction. Hence $\xi_0\in \p\C_\theta$ is impossible.

\noindent\textbf{Case~2} $\xi_{0}\in \C_{\theta}\setminus\partial\C_{\theta}$. We reformulate Eq. \eqref{F-equ} as
\begin{eqnarray}\label{weight-u-F-equ}
    \widetilde F(A)\coloneqq \left(\frac{\sigma_{k}(\widehat A)}{\sigma_{l}(\widehat A)}\right)^{\frac{1}{k-l}}=f^{\frac{1}{k-l}}(\ell u)^{\frac{p-1}{k-l}}\coloneqq u^{a}g
\end{eqnarray}
where $a$ is defined in \eqref{eqn-defn-a}, $g\coloneqq f^{\frac{1}{k-l}}\ell^{\frac{p-1}{k-l}} $ and  $\widehat A\coloneqq (\widehat A_{ij})=(\ell u_{ij}+u_{i}\ell_{j}+u_{j}\ell_{i}+u\delta_{ij})$.  Choose an orthogonal frame $\{e_{i}\}_{i=1}^{n}$ around $\xi_{0}$, such that $\widehat A$ is diagonal at $\xi_{0}$, then $\left(\widetilde F^{ij} \right)\coloneqq \left(\frac{\partial \widetilde F(\widehat A)}{\partial \widehat A_{ij}} \right)$ is also diagonal at $\xi_{0}$. 
The maximal condition implies
\begin{eqnarray}\label{weight-one-der}
    0=\n_{i} \log P= \frac{2u_k u_{ki} }{|\n u|^{2}}-\frac{\gamma h_{i}}{h}+\a \frac{\ell_i}{\ell} ,
\end{eqnarray}
and 
\begin{eqnarray}
    0&\geq &\widetilde F^{ij} \n_{ij} \log P =\widetilde{F}^{ij}\left(\frac{|\n u|^{2}_{ij}}{|\n u|^{2}}-\frac{|\n u|^{2}_{i}|\n u|^{2}_{j}}{|\n u|^{4}}-\frac{\gamma h_{ij}}{h}+\frac{\gamma h_{i}h_{j}}{h^{2}}+\alpha  \frac{\ell_{ij} }{\ell}-\a \frac{\ell_i\ell_j}{\ell^2} \right)\notag \\
    &=&\frac{2\widetilde{F}^{ij}u_ku_{kij}+2\widetilde F^{ij}u_{ki}u_{kj}}{|\n u|^{2}}-\wt F^{ij} \left( \g \frac{h_i} h-\a \frac{\ell_i}{\ell} \right) \left( \g \frac{h_j} h-\a \frac{\ell_j}{\ell} \right)  \notag  
    \\&& -\gamma \frac{\wt F^{ij}h_{ij}}  h +\g  \frac{\wt F^{ij}h_ih_j} {h^2} +\a \frac{\wt F^{ij} \ell_{ij}}{\ell} -\a \frac{\wt F^{ij} \ell_i\ell_j}{\ell^2} \notag 
    \\& =& \frac{2\widetilde{F}^{ij}u_ku_{kij}+2\widetilde F^{ij}u_{ki}u_{kj}}{|\n u|^{2}} -\gamma \frac{\wt F^{ij}h_{ij}}  h +(\g-\g^2)  \frac{\wt F^{ij}h_ih_j} {h^2} \notag 
    \\&& -(\a^2+\a) \frac{\wt F^{ij} \ell_i\ell_j }{\ell^2} +2\g\a \frac{ \wt F^{ij}  h_i\ell_j}{h\ell } + {\a}   \frac{ 1-\ell}{\ell} \sum_i \wt F^{ii}.  \label{weight-important-formula}
\end{eqnarray}

 Since $\widetilde{F}(\widehat A)$ is homogeneous of degree one, we have
  \begin{eqnarray}\label{weight-degree-one-1}
     \wt F^{ij} h_{ij}= u^ag- h \sum_{i=1}^n \wt F^{ii},
 \end{eqnarray} and 
 \begin{eqnarray} \label{weight-degree-one}
     \widetilde{F}^{ij}u_{ij}&=&\frac{1}{\ell}\widetilde{F}^{ij}\left(\widehat A_{ij}-u_{i}\ell_{j}-u_{j}\ell_{i}-u\delta_{ij}
     \right)\\ \notag
     &=&\frac{u^{a}g}{\ell}-\frac{2\widetilde{F}^{ij}u_{i}\ell_{j}}{\ell}-\frac{u}{\ell}\sum\limits_{i=1}^{n}\widetilde{F}^{ii}.
 \end{eqnarray}
 Taking the first derivatives of Eq. \eqref{weight-u-F-equ} in the $e_k$ direction, we have 
  \begin{eqnarray*}
\widetilde{F}^{ij}(\ell u_{ijk}+\ell_ku_{ij}+2u_{ik}\ell_{j}+2u_{i}\ell_{jk}+u_k\delta_{ij})=au^{a-1}u_k g+u^{a}g_k,
  \end{eqnarray*}
Together with \eqref{third comm}, \eqref{weight-degree-one} and \eqref{weight-one-der},  we obtain 
\begin{eqnarray}\label{weight-third-deri-gradient}
u_k\widetilde{F}^{ij}u_{kij}&=&
 \frac{u_k}{\ell }\left[au^{a-1}u_kg+u^{a}g_k-\widetilde{F}^{ij}(\ell_ku_{ij}+2u_{ik}\ell_{j}+2u_{i}\ell_{jk}+u_k\delta_{ij})
    \right]\notag \\
    &&+|\n u|^{2}\sum\limits_{i=1}^{n}\widetilde{F}^{ii}-\widetilde{F}^{ij}u_{i}u_{j} \notag \\
    &=&\frac{au^{a-1}g}{\ell}|\n u|^{2}+\frac{u^{a}}{\ell}\<\n u,\n g\>-\frac{u^{a}g}{\ell^2} \<\n u,\n \ell\>+\frac{2}
    {\ell^2}\<\n u,\n\ell\> \widetilde{F}^{ij}u_{i}\ell_{j}  \notag \\&& +\frac{u }{\ell^2} \<\n u,\n \ell\> \sum\limits_{i=1}^{n}\widetilde{F}^{ii}-\gamma |\n u|^{2}  \frac{\widetilde{F}^{ij} \ell_{i}h_{j}}{h\ell} + {\a |\n u|^{2}}  \frac{\widetilde{F}^{ij}\ell_{i}\ell_{j}}{\ell^2} \notag \\
    && +\left(1-\frac{2}{\ell}\right)\widetilde{F}^{ij}u_iu_j+\left(1-\frac{1}{\ell} \right)|\n u|^{2}\sum\limits_{i=1}^{n}\widetilde{F}^{ii}.~~~~~~
\end{eqnarray}
From the Cauchy--Schwarz inequality and  \eqref{weight-one-der}, 
\begin{eqnarray}\label{N-2}
   \frac{2}{|\n u|^{2}} \sum\limits_{k=1}^{n}\widetilde{F}^{ij}u_{ki}u_{kj} \geq \frac{\widetilde{F}^{ii}}{2}\left(\frac{2u_ku_{ki}}{|\n u|^{2}}\right)^{2}=\frac{\gamma^{2}}{2} \widetilde{F}^{ii} \frac{h_i^2}{h^2}+\frac{\a^{2}}{2}   \widetilde{F}^{ii} \frac{\ell_i^2}{\ell^2} -\a\gamma  \widetilde{F}^{ii} \frac{h_i\ell_i}{h\ell} . 
\end{eqnarray}
In view of
\begin{eqnarray*}
    \wt F^{ij}u_iu_j= \wt F^{ij} \frac{h_ih_j}{\ell^2}- \frac{2\wt F^{ij} h h_i\ell_j}{\ell^3} + \frac{h^2}{\ell^4} \wt F^{ij}\ell_i\ell_j ,
\end{eqnarray*}and 
\begin{eqnarray*}
     \wt F^{ij} u_i\ell_j = \wt F^{ij} h_i (\log \ell)_j - h \wt F^{ij} (\log \ell)_i (\log \ell)_j,
\end{eqnarray*}
and by inserting \eqref{N-2}, \eqref{weight-third-deri-gradient},  \eqref{weight-degree-one}, and \eqref{weight-degree-one-1}  into \eqref{weight-important-formula}, we get
\begin{eqnarray}
    0&\geq&  \frac{u^{a-1}g}{\ell}\left[ (2a-\g) +\frac{2u}{|\n u|} \left\< \frac{\n u}{|\n u|}, \n \log \left( \frac g \ell \right)  \right\>  \right]  \notag 
    \\
    &&+\left[ \frac{\a(1-\ell)}{\ell } +  \g   +2\left(1-\frac 1 \ell\right) + \frac {2u}{\ell |\n u|} \left\< \frac {\n u}{|\n u|} , \n\log \ell\right\> \right]\sum\limits_{i=1}^{n}\widetilde{F}^{ii} \notag
    \\&&   + \left[ \gamma-\frac{\g^2}2 + 2\left(1-\frac 2 \ell \right) \frac{u^2}{|\n u|^2} \right] \wt F^{ii} (\log h)_i ^2  %
      \notag \\&& +   \left[ \a-\frac{\a^2}2  - \frac{4u}{|\n u|} \left\< \frac{\n u}{|\n u|}, \n \log \ell\right\> +2 \left(1-\frac 2 {\ell}\right) \frac{ u^2}{|\n u|^2} \right ] \wt F^{ii} (\log \ell)_i^2 
    \notag \\&&  + \left[ (\a-2)\g  -4 \left(1-\frac 2 \ell\right) \frac{u^2}{|\n u|^2} + \frac{4u}{ |\n u|} \left\< \frac{\n u}{|\n u|}, \n \log \ell\right\> \right]  \wt F^{ii} (\log h)_i (\log \ell) _i.   \quad \label{less zero-condi}
\end{eqnarray}
At $\xi_0$, it follows from \eqref{N-big} that
\begin{eqnarray*} 
 \frac{2u}{|\n u|} \left\< \frac{\n u}{|\n u|}, \n \log  \left(\frac g \ell \right) \right\> \geq -\frac{C}{\sqrt{N}},~~~   \frac {2u}{\ell |\n u|} \left\< \frac {\n u}{|\n u|} , \n\log \ell\right\>  \geq - \frac{C }{  \sqrt{ N}}, 
\end{eqnarray*}
and
\begin{eqnarray*} 
 \frac{4u}{  |\n u|} \left\< \frac{\n u}{|\n u|}, \n \log \ell\right\>\geq -\frac C{\sqrt{N}},~~    4\left(1-\frac 2 \ell \right) \frac{u^2}{|\n u|^2} \geq -\frac C{N},
\end{eqnarray*}for some positive constant $C>0$, which only depends on $\|\ell\|_{C^1}$ and $\|g\|_{C^1}$. Moreover,
\begin{eqnarray*}
    -\frac C {\sqrt{N}} \wt F^{ii} (\log h)_i (\log \ell)_i \geq  -\frac C  {\sqrt{N}}  \wt F^{ii} |(\log h)_i| \geq -\frac C  {\sqrt{N}} \left( \sum_i \wt F^{ii}+\sum_i \wt F^{ii}(\log h)_i^2\right) ,
\end{eqnarray*}and
\begin{eqnarray*}
    -\frac C N \wt F^{ii} (\log \ell)_i^2 \geq-\frac C N \sum_i \wt F^{ii}.
\end{eqnarray*}
Hence \eqref{less zero-condi} yields
\begin{eqnarray}
  0 &\geq &  \frac{u^{a-1} g}{\ell}\left(2a-\gamma -\frac{C}{\sqrt N}\right)+\left[\frac{\a(1-\ell)}{\ell} +  \g   +2\left(1-\frac 1 \ell\right)- \frac{C}{\sqrt N} 
    \right]\sum\limits_{i=1}^{n}\widetilde{F}^{ii}
   \notag  \\&& +   (\a-2)\g \widetilde{F}^{ii}(\log h)_{i} (\log \ell)_{i}  + \left(\g-\frac{\g^2}2 -\frac{C}{ N }\right) \widetilde{F}^{ii}(\log h)_{i}^{2}  \notag\\&& +  \left( \a-\frac{\a^2} 2 -\frac{C}{N} \right) \wt F^{ii} (\log \ell)_i^2 . ~~~~~~~~
    \label{eqn-weighted-C1-final}
\end{eqnarray} 
We now choose $N>0$ sufficiently large such that
\begin{eqnarray}\label{eqn-large-N-control-varepsilon}
\frac C{\sqrt{N}}\leq \frac {2a-\g}2,~~~~ \frac{C}{ N }+\frac C {\sqrt{ N}} < \frac{  \gamma}{2},  ~ ~~~ \frac C  N \leq \frac 1 2  \left(\g-\frac {\g^2} 2\right) .
\end{eqnarray}Note that $\a=2$, it follows from \eqref{eqn-weighted-C1-final} and \eqref{eqn-large-N-control-varepsilon} that 
\begin{eqnarray*}
  0&\geq &\frac{u^{a-1} g}{\ell} \left( 2a-\gamma-\frac C {\sqrt{ N}} \right)+\left( \g -\frac C {\sqrt{N}} -\frac{C}{ N }\right)  \sum\limits_{i=1}^{n}\widetilde{F}^{ii}
  \\&&  + \left(\g-\frac{\g^2}2   -\frac C N \right) \widetilde{F}^{ii}(\log h)_{i}^{2}   
   \\&\geq & \frac{u^{a-1} g}{2\ell}(2a-\gamma)+\frac{\g}2 \sum\limits_{i=1}^{n}\widetilde{F}^{ii} +\frac 1 2 \left(\g-\frac{\g^2}{2}\right) \wt F^{ii}(\log h)_i^2 
  \\&> & 0, 
\end{eqnarray*}
 a contradiction. Therefore, \eqref{assum} cannot hold. This proves \eqref{u-weight-est} and hence completes the proof of Lemma \ref{lem weighted gra}.
\end{proof}

\subsection{\texorpdfstring{$C^{2}$}{} estimate}\label{sec3.2}
In this subsection, we establish the \(C^{2}\) estimate for solutions to Eq. \eqref{quo-equ} with \(p\ge 1\). When \(p> k-l+1\), the argument follows the strategy developed by Lions--Trudinger--Urbas \cite{LTU}, Ma--Qiu \cite{MQ}, and Chen--Zhang \cite{CZ2021} for the Neumann boundary problems for fully nonlinear elliptic equations. The range \(1\le p<k-l+1\) is more delicate due to the lack of a priori \(C^{0}\) estimate. To overcome this difficulty, we derive a quantitative dependence of the largest eigenvalue of \(W=\n^{2}h+h\s\) (equivalently, the maximal principal curvature radius of \(\S\)) on \(\|h\|_{C^{0}(\C_{\theta})}\). This relationship will subsequently yield the desired \(C^{0}\) estimate.

We begin by reducing the global \(C^2\) estimate of solution to Eq. \eqref{quo-equ} to the double normal derivative estimate, which yields a stronger and more refined result than that in \cite[Lemma 3.6]{MWW-Lp}. Denote  $$M\coloneqq \max\limits_{\partial \C_{\theta}}|\n^{2}h(\mu, \mu)|.$$
\begin{lem}\label{lemma-global-C2-double-normal}
    Let $0\leq l<k\leq n$ and $\theta\in (0, \pi/2)$.  
    \begin{enumerate}
        \item If $p> k-l+1$, and  $h$ is a positive admissible  solution to  Eq. \eqref{quo-equ}, then there holds
   \begin{eqnarray}\label{uniform C2}
\max\limits_{\C_{\theta}}|\n^{2}h|\leq M+C,
    \end{eqnarray}
    where $C>0$ depends only on $n, k, l, p, \min\limits_{\C_{\theta}}f, \|h\|_{C^{1}(\C_{\theta})}$, and $\|f\|_{C^{2}(\C_{\theta})}$. Furthermore, if $k-l+1< p\leq k-l+2$, the constant $C$ is independent of $p$.

    \item If $1<p<k-l+1$ and $h$ is a positive, capillary even,
strictly convex solution of Eq. \eqref{quo-equ}, then
\begin{eqnarray}\label{uniform C2-1}
\max\limits_{\C_\theta}|\n^2h|
\leq C\left[M+(\max\limits_{\C_\theta}h)^{2-\gamma}
(\min\limits_{\C_\theta}h)^{\gamma-1}
+\|h\|_{C^0(\C_\theta)}+1\right],
\end{eqnarray}
where $\gamma$ is chosen as in Lemma~\ref{lem weighted gra},
and $C>0$ depends only on $n,k,l,p,\theta,\gamma$,
$\min\limits_{\C_\theta}f$, and $\|f\|_{C^2(\C_\theta)}$.

\item If $p=1$ and $h$ is a positive, capillary even,
strictly convex solution of Eq. \eqref{quo-equ}, then
\begin{eqnarray}\label{uniform C2-p-one}
\max\limits_{\C_\theta}|\n^2h|
\leq C\left[M+\|h\|_{C^0(\C_\theta)}+1\right],
\end{eqnarray}
where $C>0$ depends only on $n,k,l,\theta$,
$\min\limits_{\C_\theta}f$, and $\|f\|_{C^2(\C_\theta)}$.
    \end{enumerate}
\end{lem}

\begin{proof}
    Consider the function
    \begin{eqnarray*}
        \mathcal{P}(\xi, \Xi)\coloneqq \n^{2}{h}(\Xi,\Xi)+{h}(\xi),
    \end{eqnarray*}
    for $\xi\in\C_{\theta}$ and a unit vector $\Xi\in T_{\xi}\C_{\theta}$. Suppose that $\mathcal{P}$ attains its maximum value at some point $\xi_{0}\in \C_{\theta}$ and some unit vector $\Xi_{0}\in T_{\xi_{0}}\C_{\theta}$. We divide the proof into two cases according to whether $\xi_{0}$ is an interior point or not.
    
    \
        
  \noindent\textbf{Case~1}~ $\xi_{0}\in \C_{\theta}\setminus \partial \C_{\theta}$. In this case, we choose an orthonormal frame $\{e_{i}\}_{i=1}^{n}$ around $\xi_{0}$ such that ${W}$ is diagonal and $\Xi_{0}=e_{1}$ at $\xi_{0}$, then $F^{ij}$ is also diagonal at $\xi_{0}$. The maximality of $\mathcal{P}$ implies 
    \begin{eqnarray*}
        0=\mathcal{P}_{i}={h}_{11i}+{h}_{i},
    \end{eqnarray*}
    and 
    \begin{eqnarray}\label{2-deri}
        0\geq F^{ij}\mathcal{P}_{ij}=F^{ii}{h}_{11ij}+F^{ii}{h}_{ii}.
    \end{eqnarray}
    On $\C_{\theta}$,  the commutation formula gives
    \begin{eqnarray*}
h_{rsij}=h_{ijrs}+2h_{rs}\d_{ij}-2h_{ij}\d_{rs}+h_{si}\d_{rj}-h_{rj}\d_{is}.
\end{eqnarray*}
Inserting this into \eqref{2-deri}, using the one-homogeneity of \(F\) and \eqref{F-equ}, we obtain
    \begin{eqnarray}
0&\geq&F^{ij}\mathcal{P}_{ij}
=F^{ii}({h}_{ii11}+2{h}_{11}-2{h}_{ii})+F^{ii}h_{ii}\notag \\
&=&F^{ii}W_{ii11}+(h_{11}+h)\sum\limits_{i=1}^{n}F^{ii}-\widehat f h^{a}.  \label{2-key-ine}
    \end{eqnarray}
Differentiating  Eq. \eqref{F-equ}  once and twice in the $e_{1}$ direction, respectively,  we get 
    \begin{eqnarray}\label{1-time-deri}
        F^{ii}W_{ii1}= \widehat{f}_{1}h^{a}+ah^{a-1}h_{1}\widehat{f},
    \end{eqnarray}
    and 
    \begin{eqnarray} \label{2-times-deri}
        F^{ii}W_{ii11}&=&h^{a-1}\left(h\widehat{f}_{11}+2a \widehat{f}_{1}{h}_{1}+a(a-1)\frac{ \widehat{f} h_{1}^{2}}{h}+ah_{11} \widehat{f}
        \right)-F^{ij,pq}W_{ij1}W_{pq1} \notag\\
        &\geq & h^{a-1}\left( h\widehat{f}_{11}+2a \widehat{f}_{1}{h}_{1}+a(a-1)\frac{\widehat{f}  h_{1}^{2}}{h}+ah_{11} \widehat{f}
        \right),
    \end{eqnarray}
 where we used the concavity of $F$ in the last inequality. Substituting \eqref{2-times-deri} into \eqref{2-key-ine}, we get
\begin{eqnarray}
    0&\geq& (h_{11}+h)\sum\limits_{i=1}^{n}F^{ii}+ h^{a-1}\bigg[ h\widehat{f}_{11}+2a \widehat{f}_{1}{h}_{1}+a(a-1)\frac{\widehat{f}  h_{1}^{2}}{h} \notag \\
    &&+a(h_{11}+h) \widehat{f}-(1+a)h\widehat{f}
        \bigg].\label{key-ine-C2thm}  ~~~~~~~~  
\end{eqnarray}
Next, we divide the argument into three subcases according to the range of \(p\).

\

\noindent\textbf{Subcase~1.1}~ $p> k-l+1$.  In this case, \(a> 1\). It follows from \eqref{key-ine-C2thm} that
\begin{eqnarray*}
   0\geq (h_{11}+h)\sum\limits_{i=1}^{n}F^{ii}+h^{a-1}\left(h \widehat{f}_{11}+2a\widehat{f}_{1}h_{1}-(1+a)h\widehat{f}
   \right).
\end{eqnarray*}
Together with Proposition \ref{pro-2.3} (4) and Lemma~\ref{lem kC0}, we deduce that
\begin{eqnarray}\label{1}
    h_{11}+h\leq C,
\end{eqnarray}
where $C>0$ depends only on $n, k, l, p,  \min\limits_{\C_{\theta}}f, \|h\|_{C^{1}(\C_{\theta})}$ and $\|f\|_{C^{2}(\C_{\theta})}$. 

Moreover, if $k-l+1< p\leq k-l+2$, then \(1< a\le 2\), and hence the constant \(C\) in \eqref{1} can be chosen independently of \(p\).

\

\noindent\textbf{Subcase~1.2}  $1<p<k-l+1$.  In this case, we have $0<a< 1$. Applying \eqref{key-ine-C2thm}  again yields 
\begin{eqnarray}\label{sub-1.2}
h_{11}+h \leq (1-a)\frac{h_{1}^{2}}{h}-  \frac{2\wh f_1}{\wh f}  h_{1}-\frac{\widehat{f}_{11} }{ a \wh f} h+\left(1+\frac{1}{a}\right)h .
\end{eqnarray}
By Lemma  \ref{lem weighted gra},   for some $\g\in(0, 2)$,  we have  
\begin{eqnarray*}
    \frac{h_{1}^{2}}{h}\leq  \frac{|\n h|^{2}}{h^{\gamma}}h^{\gamma-1}\leq C (\max\limits_{\C_{\theta}} h)^{2-\gamma} h^{\gamma-1}.
\end{eqnarray*}
If $0<\gamma\leq 1$, then  
\begin{eqnarray}\label{weighted-gre-control-1}
    \frac{h_1^2} h \leq C (\max_{\C_\theta} h)^{2-\g} (\min_{\C_\theta} h)^{\g-1}.
\end{eqnarray}
On the other hand, if $1<\gamma<2$, then
\begin{eqnarray}\label{weighted-gre-control-2}
    \frac{h_{1}^{2}}{h}\leq  \frac{|\n h|^{2}}{h^{\gamma}}h^{\gamma-1}\leq C (\max\limits_{\C_{\theta}} h)^{2-\gamma} h^{\gamma-1}\leq C\max\limits_{\C_{\theta}}h.
\end{eqnarray}
In view of \eqref{weighted-gre-control-1} and \eqref{weighted-gre-control-2}, we obtain
\begin{eqnarray}\label{weighted-gre-control-final}
    \frac{h_{1}^{2}}{h}\leq C (\max_{\C_\theta} h)^{2-\g} (\min_{\C_\theta} h)^{\g-1}+ C\max\limits_{\C_{\theta}}h.
\end{eqnarray}
Moreover, by \cite[Lemma 3.3 and its proof]{MWW-Weingarten}, we derive that
\begin{eqnarray}\label{hc1}
    \max\limits_{\C_{\theta}}|\n h|\leq (1+\cot^{2}\theta)^{\frac{1}{2}}\|h\|_{C^{0}(\C_{\theta})}.
\end{eqnarray}
Substituting \eqref{hc1} and \eqref{weighted-gre-control-final} into \eqref{sub-1.2}, we conclude that
\begin{eqnarray}\label{2}
    h_{11}+h\leq C \left[ (\max_{\C_\theta} h)^{2-\g} (\min_{\C_\theta} h)^{\g-1}+\|h\|_{C^{0}(\C_{\theta})}+1\right],
\end{eqnarray}
 where $C$ depends only on $n, k,l, p,  \min\limits_{\C_{\theta}}f$ and $\|f\|_{C^{2}(\C_{\theta})}$.

 \

\noindent\textbf{Subcase 1.3} $p=1$.  In this case, $a=0$. It follows from \eqref{key-ine-C2thm} that 
 \begin{eqnarray*}
     0\geq (h_{11}+h)\sum\limits_{i=1}^{n}F^{ii}+\widehat{f}_{11}-\widehat{f}.
 \end{eqnarray*}
Invoking Proposition \ref{pro-2.3} (4) once again, we obtain
 \begin{eqnarray}\label{3}
     h_{11}+h\leq C,
 \end{eqnarray}
where $C$ depends only on $n$, $k$, $l$, $\min\limits_{\C_\theta}f$, and $\|f\|_{C^2(\C_\theta)}$.

\
		
\noindent\textbf{Case 2}  $ \xi_{0} \in \partial \C_{\theta}$. In this case, the argument proceeds as in \cite[Proof of Lemma 3.3, Case 2]{MWW-AIM}, yielding
		\begin{eqnarray}\label{1n}
			\n^2	h(\Xi_0,\Xi_0)  \leq  |\n^{2}h(\mu,\mu)|(\xi_{0})+2 \|h\|_{C^0(\C_{\theta})}.
		\end{eqnarray}

        \
        
Finally, combining \eqref{1} and \eqref{1n} gives the estimate
in the supercritical case.
For $1<p<k-l+1$, \eqref{2} and \eqref{1n} bound the largest
eigenvalue of $W$, and $W>0$ bounds every eigenvalue of
$\n^2h$ from below by $-h$. This proves \eqref{uniform C2-1}.
For $p=1$, the same argument, using \eqref{3} instead of
\eqref{2}, proves \eqref{uniform C2-p-one}.

\end{proof}

Next, we derive an estimate for the double normal derivative of $h$ on the boundary, namely Lemma \ref{c2 est-admissible soluion-big case} and Lemma \ref{c2 est-convex soluion-small case} in the cases $p>k-l+1$ and $1\leq p<k-l+1$, respectively. To this end, in each case we construct two barrier functions (more precisely, \eqref{eqn:overline-Q}, \eqref{eqn:underline-Q}, and \eqref{eqn:overline-fancy-Q}, \eqref{eqn:underline-fancy-Q} for the cases $p>k-l+1$ and  $1\leq p<k-l+1$, respectively) for $h_\mu$ in a small neighborhood of the boundary. The approach and construction are inspired by Lions--Trudinger--Urbas \cite{LTU}, Ma--Qiu \cite[Section 4]{MQ}, and, in particular, Chen--Zhang \cite[Section 4]{CZ2021}. 
 We  first introduce the function 
	\begin{eqnarray*}
		\zeta(\xi)\coloneqq e^{-  d(\xi)}-1,
	\end{eqnarray*}where the function ${d}$ is defined as in \eqref{d-func}. 
Note that $\zeta(\xi)$ is well-defined on $\xi\in \C_\theta\setminus \{(1-\cos\theta) E_{n+1}\}$. This function was previously used by  Guan \cite[Lemma~3.1]{Gb99} and also in our recent work \cite[Section 3.3]{MWW-AIM}. A direct computation gives \begin{eqnarray}\label{eqn-bdry-zeta}
    \zeta|_{\partial \C_{\theta}}=0\quad  {\rm{and}}\quad  \n \zeta|_{\partial \C_{\theta}}=\mu.
\end{eqnarray} 
In the domain considered below, $d$ is the geodesic distance to
$\partial\C_\theta$, before the smooth extension used in
\eqref{d-func}. Set
\begin{eqnarray*}
\Omega_{\delta_0}
\coloneqq\{\xi\in\C_\theta:
{\rm dist}(\xi,\partial\C_\theta)\leq\delta_0\},
\qquad
\epsilon_0\coloneqq
\frac{c_{n,k,l}}{2(c_{n,k,l}+1)}\min\{\cot\theta,1\}.
\end{eqnarray*}
For $0<\delta_0\leq\min\{\theta/2,1/2\}$, direct computation gives
\begin{eqnarray*}
\n^2_{ij}\zeta
=e^{-d}\left[\cot(\theta-d)(\s_{ij}-d_i d_j)+d_i d_j\right].
\end{eqnarray*}
Consequently,
\begin{eqnarray}\label{hessian of zeta}
(\n^2_{ij}\zeta)
\geq\frac{c_{n,k,l}+1}{c_{n,k,l}}\epsilon_0\s
\quad\text{in }\Omega_{\delta_0}.
\end{eqnarray}
By Proposition~\ref{pro-2.3} (4), this implies
\begin{eqnarray}\label{eqn:Fij-zeta-ij-sum}
F^{ij}\zeta_{ij}
\geq\epsilon_0\left(\sum\limits_{i=1}^{n}F^{ii}+1\right).
\end{eqnarray}

\begin{lem}\label{c2 est-admissible soluion-big case}
	Let $0\leq l<k\leq n$ and $\theta\in (0, \pi/2)$.  Suppose that  $p> k-l+1$ and $h$ is an admissible solution to Eq. \eqref{quo-equ}. Then
		\begin{eqnarray*}
			\max\limits_{\p \C_{\theta}}|\n^{2}h(\mu,\mu)|\leq C,
		\end{eqnarray*}
		where $C>0$ depends only on $n, k, l, p,  \min\limits_{\C_{\theta}}f, \min\limits_{\C_{\theta}}h, \|f\|_{C^{1}(\C_{\theta})}$ and  $\|h\|_{C^1(\C_{\theta})}$. Furthermore, if $k-l+1< p\leq k-l+2$,  the constant $C$ is independent of $p$.
	\end{lem}
\begin{proof}
Define the function

    \begin{eqnarray}\label{eqn:overline-Q}
       \overline Q(\xi)\coloneqq -(1+\beta  d)\left(\<\n {h}, \n  d\>+\cot\theta {h}\right)-\left(B_{1}+\frac{M}{2} \right)\zeta, \quad \xi\in \Omega_{\delta_{0}},
    \end{eqnarray}
where $\beta>1$, $B_{1}$, and $\delta_0$ are positive constants to be determined later. Clearly, by \eqref{eqn:bdry-distance} and \eqref{eqn-bdry-zeta}, \begin{eqnarray}\label{eqn:Q-bdry}
 \overline   Q=0, ~~~\partial \C_\theta.
\end{eqnarray}

We next prove that
\begin{eqnarray}\label{eqn-Q-nonnegative}
	 \overline	Q(\xi)\geq 0,~~\quad{\rm{in} } \quad \Omega_{\delta_{0}}.
\end{eqnarray} 
To this end, let $\xi^*\in\Omega_{\delta_{0}}\setminus\partial\Omega_{\delta_{0}}$ be a point at which $ \overline Q$ attains its minimum.
Choose a local orthonormal frame $\{e_i\}_{i=1}^{n}$ around $\xi^*$ such that $W=\n^2 h+h\s$ is diagonal at $\xi^*$. Then $F^{ij}$ is also diagonal at $\xi^*$. A direct calculation yields
\begin{eqnarray}\label{Q-one-deri}
 0=  \overline Q_{i}=-\beta d_{i}\left(\<\n h, \n d\>+\cot\theta h\right)-(1+\beta  d)\left(h_{ki} d_{k}+h_k d_{ki}+\cot\theta h_{i}\right)-\left(B_{1}+\frac{M}{2}\right)\zeta_{i},
\end{eqnarray}
and
\begin{eqnarray}\label{twice deri}
     \overline Q_{ij}&=&-\beta  d_{ij}\left(\<\n h, \n d\>+\cot\theta {h}\right)-\beta  d_{i}({h}_{kj} d_k+{h}_k d_{kj}+\cot\theta {h}_{j})-\beta  d_{j}({h}_{ki} d_k+{h}_k d_{ki} \notag \\
    &&+\cot\theta {h}_{i})-(1+\beta  d)({h}_{kij} d_k+{h}_{ki} d_{kj}+{h}_{kj} d_{ki}+{h}_k d_{kij}+\cot\theta {h}_{ij}) \notag \\
    &&-\left(B_{1}+\frac{M}{2}\right)\zeta_{ij}.
\end{eqnarray}
Recall that the commutation formula on $\C_\theta$,
\begin{eqnarray}\label{eqn:C2-double-normal-third commutator}
h_{kij}=h_{ijk}+h_k\d_{ij}-h_j\d_{ki}=W_{ijk}-h_j\delta_{ki}.
\end{eqnarray}
 Contracting \eqref{twice deri} with $F^{ij}$, and using Lemmas \ref{lem kC0}--\ref{lem c1} and \eqref{eqn:C2-double-normal-third commutator}, at $\xi=\xi^*$,  we obtain
\begin{eqnarray}
   F^{ij}  \overline Q_{ij} 
   &=&-\beta(\<\n h,\n d\>+\cot\theta{h})F^{ij} d_{ij}-2\beta F^{ij} d_{i}({h}_{kj} d_{k}+{h}_{k} d_{kj}+\cot\theta {h}_{j}) \notag \\
   &&- (1+\beta  d)F^{ij}({h}_{kij} d_k+2{h}_{ki} d_{kj}+{h}_{k} d_{kij}+\cot\theta {h}_{ij}) \notag \\
   && -\left(B_{1}+\frac{M}{2} \right)F^{ij}\zeta_{ij} \notag
   \\&
     \leq & -2\beta  F^{ij}W_{kj}  d_id_k -(1+\beta  d)F^{ij}W_{ijs} d_{s}-2(1+\beta  d)F^{ij}W_{ik} d_{kj}\notag\\
 &&- \left(B_{1}+\frac{M}{2} \right)F^{ij}\zeta_{ij}+C(\beta+1)\sum\limits_{i=1}^{n}F^{ii}. \label{eqn:FijQij}
\end{eqnarray}
From \eqref{1-time-deri}, we see
\begin{eqnarray}\label{s}
    |F^{ij}W_{ijs}d_{s}|=|\< \n \widehat{f}, \n d\> h^{a}+a\widehat{f}h^{a-1}\< \n h, \n d\>|\leq C.
\end{eqnarray}

For any $\beta>1$, choose
$0<\delta_0\leq\min \left\{\theta/2,1/(9\beta+10) \right\}$.
Then
\begin{eqnarray}\label{eqn:choice-delta-0}
   1\leq 1+\beta d(\xi)\leq \frac{10}{9},\quad  \text{and}\quad \frac{9}{10}\leq e^{-d(\xi)}\leq 1,~~~\text{ on }~~~ \xi\in \O_{\delta_0}.
\end{eqnarray}
Inserting \eqref{s} into \eqref{eqn:FijQij}, and using \eqref{eqn:Fij-zeta-ij-sum} together with Proposition \ref{pro-2.3} (4), we obtain
\begin{align}
F^{ij}\overline Q_{ij}
&\leq -2\beta F^{ii}W_{ii}d_i^2
      -2(1+\beta d)F^{ii}W_{ii}d_{ii}\notag\\
&\quad-\left[\left(B_1+\frac M2\right)\epsilon_0
             -C(\beta+1)\right]
       \left(\sum_{i=1}^n F^{ii}+1\right).
\label{eqn-control-term}
\end{align}
Here the constant $C>0$ depends only on $n, k, l, a, \min\limits_{\C_{\theta}}f, \min\limits_{\C_{\theta}}h, \|f\|_{C^{2}(\C_{\theta})}$, and  $\|h\|_{C^1(\C_{\theta})}$.

\ 

We now \textbf{claim} that, by choosing the constants $B_{1}$, $\beta$, and $\delta_{0}$ appropriately (cf. \eqref{eqn-beta-delta0-fix} and \eqref{eqn-B-1-fix}), one can guarantee that
\begin{eqnarray}\label{F>0}
    F^{ij}  \overline Q_{ij}<0, ~~~\text{ at } ~~\xi=\xi^*,
\end{eqnarray}
and 
\begin{eqnarray}\label{bry-condi}
   \overline  Q> 0, ~~~~\text{ on } ~\partial\Omega_{\delta_{0}}\setminus \partial \C_{\theta}. 
   \end{eqnarray}

\ 

  Assuming this \textbf{claim}, \eqref{F>0} contradicts the minimality of $ \overline Q$ at the interior point
$\xi^*\in \Omega_{\delta_{0}}\setminus \partial \Omega_{\delta_{0}}$.
Hence, $ \overline Q$ attains its minimum on $\partial\Omega_{\delta_0}$.
Together with \eqref{eqn:Q-bdry} and \eqref{bry-condi}, this yields \eqref{eqn-Q-nonnegative}.   
Let $\eta_0 \in \partial \C_\theta$ be a point at which $h_{\mu\mu}$ attains its maximum on $\partial \C_\theta$, i.e.,
\begin{eqnarray*}
    h_{\mu\mu}(\eta_0)=\sup_{\partial \C_\theta} h_{\mu\mu}>0.
\end{eqnarray*}
By \eqref{eqn:Q-bdry}, \eqref{eqn-Q-nonnegative}, and \eqref{eqn:bdry-distance}, we obtain
\begin{eqnarray*}
0&\geq &  \overline Q_{\mu}(\eta_{0})\\
&\geq &-(h_{s\mu}d_{s}+h_{s}  d_{s\mu})-\left(B_{1}+\frac{M}{2}\right)\zeta_{\mu}-\cot\theta h_{\mu}\\
&\geq& h_{\mu\mu}(\eta_{0})-\left(B_{1}+\frac{M}{2}\right)-C,
\end{eqnarray*}
which  implies  
\begin{eqnarray}\label{sup estimate-1}
\sup\limits_{\partial \C_{\theta}}h_{\mu\mu}\leq C+\frac{M}{2}.
\end{eqnarray}
Here $C>0$ depends only on $n, k, l, a=\frac{p-1}{k-l},  \min\limits_{\C_{\theta}}f, \min\limits_{\C_{\theta}}h, \|f\|_{C^{2}(\C_{\theta})}$, and  $\|h\|_{C^1(\C_{\theta})}$.

\ 

We now turn to the proof of the \textbf{claim} that \eqref{F>0} and \eqref{bry-condi} hold.
To this end, we partition the index set $\{1,2,\cdots,n\}$ into
\begin{eqnarray*}
    \mathcal{B}\coloneqq \left\{i: \beta d_{i}^{2} <\frac{1}{n}  \right\},~~~\text{ and } ~~~
    \mathcal{G}\coloneqq \left\{i: \beta d_{i}^{2}\geq \frac{1}{n}  \right\}.
\end{eqnarray*}
Since $\sum\limits_{i=1}^{n}d_{i}^{2}=|Dd|^{2}=1$, there exists an index $  i_{0}\in \{1,2,\cdots, n\}$ such that $d_{i_{0}}^{2}\geq \frac{1}{n}$. Hence $i_{0}\in\mathcal{G}$. Without loss of generality, we assume that $i_{0}=1$. By the concavity of $d$ in $\O_{\delta_0}$, it then follows from \eqref{eqn-control-term} that
\begin{eqnarray}
     F^{ij} \overline Q_{ij}&\leq& -2\beta \sum\limits_{i\in \mathcal{G}}F^{ii}d_{i}^{2}W_{ii}-2\beta\sum\limits_{i\in \mathcal{B}}F^{ii}d_{i}^{2}W_{ii}-2(1+\beta d)\sum\limits_{W_{ii}>0}F^{ii}W_{ii}d_{ii}\notag \\
    &&-2(1+\beta d)\sum\limits_{W_{ii}<0}F^{ii}W_{ii}d_{ii}-\left[ \left(B_{1}+\frac{M}{2} \right) \epsilon_0 -C(\beta+1)\right]\left(\sum\limits_{i=1}^{n}F^{ii}+1\right).\notag \\
    &\leq &-2\beta \sum\limits_{i\in \mathcal{G}}F^{ii}d_{i}^{2}W_{ii}-2\beta\sum\limits_{i\in \mathcal{B}}F^{ii}d_{i}^{2}W_{ii}+C\sum\limits_{W_{ii}>0}F^{ii}W_{ii}\notag \\
    &&-\left[ \left(B_{1}+\frac{M}{2} \right) \epsilon_0 -C(\beta+1)\right]\left(\sum\limits_{i=1}^{n}F^{ii}+1\right).\label{key-F-ineq} 
\end{eqnarray}

For each $i\in \mathcal{G}$, it follows from \eqref{Q-one-deri} that
\begin{eqnarray}\label{Wii-1}
    W_{ii}=\frac{e^{-d}}{1+\beta d}\left(B_{1}+\frac{M}{2}\right)-\frac{\beta}{1+\beta d}\left(\<\n h,\n d\> +\cot\theta h\right)+\frac{1}{d_{i}}\left(
    hd_{i}-h_{s}d_{si}-\cot\theta h_{i}\right).
\end{eqnarray}
By \eqref{eqn:choice-delta-0}, for $i\in \mathcal{G}$, we have
\begin{eqnarray}\label{control}
    \left|  -\frac{\beta}{1+\beta d}\left(\<\n h,\n d\>+\cot\theta h\right)+\frac{1}{d_{i}}\left(
    hd_{i}-h_{s}d_{si}-\cot\theta h_{i}\right)  \right| \leq C_1 \beta,
\end{eqnarray}where $C_1>0$ depends only on $n$ and $\|h\|_{C^1(\C_\theta)}$. 
Choosing \begin{eqnarray}\label{B11}
    B_{1}\geq B_{11}\coloneqq 5 {C_1\beta}, 
\end{eqnarray}and combining \eqref{control}, \eqref{Wii-1}, and \eqref{eqn:choice-delta-0}, we obtain
\begin{eqnarray}\label{good-index}
0<   \frac{3B_{1}}{5}+\frac{2M}{5}\leq W_{ii}\leq \frac{6B_{1}}{5}+\frac{M}{2},\quad \forall~ i\in \mathcal{G},
\end{eqnarray}
consequently, 
\begin{eqnarray}\label{good}
    -2\beta\sum\limits_{i\in \mathcal{G}} F^{ii}W_{ii}d_{i}^{2}\leq -2\beta F^{i_{0}i_{0}}W_{i_{0}i_{0}}d_{i_{0}}^{2}\leq -\frac{2\beta}{n}F^{11}W_{11}.
\end{eqnarray}
Since $F(A)$ is homogeneous of degree one, we have 
\begin{eqnarray*}
    -2\beta \sum\limits_{i\in \mathcal{B}}F^{ii}d_{i}^{2}W_{ii}&\leq& -2\beta \sum\limits_{i\in \mathcal{B}, W_{ii}<0} F^{ii}d_{i}^{2}W_{ii}\leq -\frac{2}{n}\sum\limits_{W_{ii}<0}F^{ii}W_{ii} \notag \\
    &=& -\frac{2}{n}\left(\widehat{ f}h^{a}-\sum\limits_{W_{ii}>0}F^{ii}W_{ii}\right)\leq \frac{2}{n}\sum\limits_{W_{ii}>0}F^{ii}W_{ii}.
    \end{eqnarray*}
Together with \eqref{good}, it follows from \eqref{key-F-ineq} that 
    \begin{eqnarray}
     F^{ij}  \overline Q_{ij} &\leq & -\frac{2\beta}{n}F^{11}W_{11}+\left(\frac{2}{n}+C_{0}\right)\sum\limits_{W_{ii}>0}F^{ii}W_{ii} \notag \\& & -\left[ \left(B_{1}+\frac{M}{2} \right)\epsilon_0 -C_0(\beta+1)\right]\left(\sum\limits_{i=1}^{n}F^{ii}+1\right). \label{key-F-ine-1}
    \end{eqnarray}
Here $C_0>0$ depends only on $n, k, l, a,  \min\limits_{\C_{\theta}}f, \min\limits_{\C_{\theta}}h, \|f\|_{C^{1}(\C_{\theta})}$, and  $\|h\|_{C^1(\C_{\theta})}$. 
    
We now distinguish three cases. Without loss of generality, we assume that $$W_{22}\geq W_{33}\geq \cdots\geq W_{nn}.$$

\noindent\textbf{Case~1}~ $W_{nn}\geq 0$. In this case,  
\begin{eqnarray*}
    \sum\limits_{W_{ii}>0}F^{ii}W_{ii}=F^{ii}W_{ii}=\widehat{f}h^{a}.
\end{eqnarray*}

Substituting this into \eqref{key-F-ine-1}, we obtain\begin{eqnarray*}
F^{ij}\overline Q_{ij}
&\leq&\left(\frac{2}{n}+C_0\right)\widehat f h^a
-\left[\left(B_1+\frac{M}{2}\right)\epsilon_0
-C_0(\beta+1)\right]
\left(\sum\limits_{i=1}^{n}F^{ii}+1\right)\\
&\leq&C_2
-\left[\left(B_1+\frac{M}{2}\right)\epsilon_0
-C_0(\beta+1)\right]
\left(\sum\limits_{i=1}^{n}F^{ii}+1\right).
\end{eqnarray*}
Here $C_2>0$ depends only on $n, k, l, a,  \min\limits_{\C_{\theta}}f, \min\limits_{\C_{\theta}}h, \|f\|_{C^{1}(\C_{\theta})}$, and  $\|h\|_{C^1(\C_{\theta})}$. 
Choosing 
\begin{eqnarray}\label{B12}
    B_{1}>B_{12}\coloneqq \frac{C_{0}(\beta+1)}{ \epsilon_0 }+\frac{C_2}{\epsilon_0(c_{n,k,l}+1) } , 
\end{eqnarray} together with Proposition \ref{pro-2.3}  (4), we obtain \eqref{F>0}.

\

\noindent\textbf{Case 2} $W_{nn}<0$ and $-W_{nn}\leq \frac{\epsilon_0}{10 \left(\frac{2}{n}+C_{0}\right)}W_{11}$. A direct calculation yields
\begin{eqnarray*}
    \left(\frac{2}{n}+C_{0}\right) \sum\limits_{W_{ii}>0}F^{ii}W_{ii}&=&\left(\frac{2}{n}+C_{0}\right)\left(\widehat{f}h^{a}-\sum\limits_{W_{ii}<0}F^{ii}W_{ii}\right)\\
    &\leq &\left(\frac{2}{n}+C_{0}\right)\left(\widehat{f}h^{a}-W_{nn}\sum\limits_{i=1}^{n}F^{ii}\right)\\
    &\leq &C_3+\frac{\epsilon_0}{10}W_{11}\sum\limits_{i=1}^{n}F^{ii}.
\end{eqnarray*}Here $C_3>0$ depends only on $n, \|h\|_{C^0(\C_\theta)}$ and $C_0$. Together  with \eqref{good-index}, this gives
\begin{eqnarray*}
    \left(\frac{2}{n}+C_{0}\right) \sum\limits_{W_{ii}>0}F^{ii}W_{ii}\leq 2C_3+\frac{\epsilon_0}{10}\left(\frac{6B_{1}}{5}+\frac{M}{2}\right)\sum\limits_{i=1}^{n}F^{ii}.
\end{eqnarray*}
Inserting the above inequality into \eqref{key-F-ine-1}, we get
\begin{eqnarray*}
    F^{ij}  \overline Q_{ij}\leq 2C_3-\left[ \left(B_{1}+\frac{M}{2} \right)\epsilon_0-C_{0}(\beta+1)-\frac{\epsilon_0}{10}\left(\frac{6B_{1}}{5}+\frac{M}{2}\right)\right]\left(\sum\limits_{i=1}^{n}F^{ii}+1\right).
\end{eqnarray*}
Choosing\begin{eqnarray}\label{B13}
    B_{1}>B_{13}\coloneqq\frac{25}{22\epsilon_0}\left[C_{0}(\beta+1) +\frac{2C_3}{c_{n,k,l}+1}  \right]+1,
\end{eqnarray} together with Proposition \ref{pro-2.3}  (4), we again obtain \eqref{F>0}.

\

\noindent\textbf{Case 3} $W_{nn}<0$ and $-W_{nn}\geq \frac{\epsilon_0}{10 \left(\frac{2}{n}+C_{0}\right)}W_{11}$.
By \eqref{uniform C2}, there exists a constant $C_4>1$, depending only on $n, k, l, p, \min\limits_{\C_{\theta}}f, \|h\|_{C^{1}(\C_{\theta})}$, and $\|f\|_{C^{2}(\C_{\theta})}$, such that
\begin{eqnarray*}
    W_{22}\leq M+C\leq C_4\left(1+M\right).
    \end{eqnarray*}   
From \eqref{B13}, we have $B_{1}>1$. Combining this with \eqref{good-index}, it follows that
\begin{eqnarray}\label{eqn:case3-W11-comparable-22}
    W_{11}\geq \frac 2 5 (1+M) \geq  \frac{2}{5C_4} W_{22}.
\end{eqnarray}
When $k\geq 2$, by \cite[Lemma~2.7, Eq. (20)]{CZ2021},  using \eqref{eqn:case3-W11-comparable-22}, we have  
\begin{eqnarray*}
    F^{11}\geq \epsilon_{1}\sum\limits_{i=1}^{n}F^{ii},
\end{eqnarray*}
where $\epsilon_{1}>0$ depends only on $n, \epsilon_0, C_{0}$, and $C_4$. 
\\ When $k=1$, then $l=0$ and $F=\sigma_{1}(W)$, thus
\begin{eqnarray*}
    F^{11}= \frac 1n \sum\limits_{i=1}^{n}F^{ii}.
\end{eqnarray*}
It follows from \eqref{key-F-ine-1} that
\begin{eqnarray*}
    F^{ij} \overline Q_{ij}\leq -\frac{2\beta \epsilon_{1}}{n}\left(\frac{3B_{1}}{5}+\frac{2M}{5}\right)\sum\limits_{i=1}^{n}F^{ii}+ \left(\frac{2}{n}+C_{0}\right)C_4(1+M)\sum\limits_{i=1}^{n}F^{ii}.
\end{eqnarray*}
Choosing \begin{eqnarray}\label{eqn:beta-choice}
\beta\geq \frac{5nC_4  }{2\epsilon_{1} } \left(\frac{2}{n}+C_{0}\right),
\end{eqnarray} we deduce that \eqref{F>0} holds.

On $\partial\Omega_{\delta_{0}}\setminus\partial\C_{\theta}$, we have 
\begin{eqnarray*}
   \overline  Q\geq -2(|\n h|+\cot\theta |h|)+ \left(B_{1}+\frac{M}{2}\right)\left(1-e^{-\delta_{0}}\right) \geq -C_5+B_{1}\left(1-e^{-\delta_{0}}\right).
\end{eqnarray*}
Thus, we choose \begin{eqnarray}\label{B14}
    B_{1}>B_{14}\coloneqq \frac{C_5}{1-e^{-\delta_{0}}},
\end{eqnarray}
then \eqref{bry-condi} follows.

In summary, we fix $\beta, \delta_0$ and $B_1$ in \eqref{eqn:overline-Q} as 
\begin{eqnarray}\label{eqn-beta-delta0-fix}
\beta=\frac{5nC_4}{2\epsilon_1}
\left(\frac{2}{n}+C_0\right)+1,
\qquad
\delta_0=\min\left\{\frac{\theta}{2},
\frac{1}{9\beta+10}\right\},
\end{eqnarray}
and 
\begin{eqnarray}\label{eqn-B-1-fix}
    B_{1}=\max \left\{B_{11}, B_{12}, B_{13}, B_{14}\right\}+1. 
\end{eqnarray}These constants depend only on $n, k, l, a=\frac{p-1}{k-l},  \min\limits_{\C_{\theta}}f, \min\limits_{\C_{\theta}}h, \|f\|_{C^{1}(\C_{\theta})}$, and  $\|h\|_{C^1(\C_{\theta})}$.  
With these choices, \eqref{eqn:choice-delta-0}, \eqref{B11}, \eqref{B12}, \eqref{B13}, \eqref{eqn:beta-choice}, and \eqref{B14} are all satisfied.
Combining the conclusions in the above three cases, we see that both \eqref{F>0} and \eqref{bry-condi} hold. This completes the proof of the \textbf{claim}.
  
\ 

Having established the upper bound in \eqref{sup estimate-1} for the double normal derivative, we briefly outline the derivation of the lower bound similarly. The constants $C_0,C_1,C_4,C_5$
may be chosen large enough in the preceding argument to satisfy
the estimates below as well. For the domain $\Omega_{\delta_0}$ with $\delta_0$ fixed in \eqref{eqn-beta-delta0-fix}, we define the function
	\begin{eqnarray}\label{eqn:underline-Q}
		\underline{Q}(\xi)\coloneqq -(1+\underline \beta d)\left(\<\n{h}, \n  d \>+\cot\theta{h}\right)+\left({B}_{2}+\frac{M}{2} \right)\zeta, \quad \xi\in \Omega_{\delta_{0}},
	\end{eqnarray}
where 	$\underline\beta, B_2>0$ are positive constants. If $\underline Q$ attains its maximum at some interior point $\xi_*\in\Omega_{\delta_{0}}\setminus\partial\Omega_{\delta_{0}}$. For suitable choices of $\underline \beta$, $\delta_0$ and $B_2$, a similar argument as above yields
\begin{eqnarray}\label{eqn:lower-barrier-h-mu}
    F^{ij}  \underline Q_{ij}>0, ~~~\text{ at } ~~\xi=\xi_*,
\end{eqnarray}
and 
\begin{eqnarray}\label{eqn:lower-barrier-h-mu-bdry}
   \underline  Q < 0,  ~~~~\text{ on } ~\partial\Omega_{\delta_{0}}\setminus \partial \C_{\theta}.
   \end{eqnarray}
However, \eqref{eqn:lower-barrier-h-mu} contradicts the fact that $\underline Q$ attains its maximum at an interior point $\xi_*$. Hence, by $\underline Q=0$ on $\partial\C_\theta$ and \eqref{eqn:lower-barrier-h-mu-bdry}, we deduce that
\begin{eqnarray*}
    \underline{Q}(\xi)\leq 0 ~~~\text{ in }   ~ \Omega_{\delta_{0}}.
\end{eqnarray*}
Since $\underline Q=0$ on $\partial\C_\theta$, arguing as before yields
	\begin{eqnarray}\label{eqn:double-normal-h-lower-bound}
		\inf\limits_{\partial \C_{\theta}}{h}_{\mu\mu}\geq -C-\frac{M}{2}.
	\end{eqnarray}

Finally, in view of \eqref{sup estimate-1} and \eqref{eqn:double-normal-h-lower-bound},  we conclude that  
	\begin{eqnarray}\label{eqn:C2-double-normal-p-large}
		\sup\limits_{\partial \C_{\theta}}|{h}_{\mu\mu}|\leq C,
	\end{eqnarray}
    where the constant $C$ depends on $n, k, l, a,  \min\limits_{\C_{\theta}}f, \min\limits_{\C_{\theta}}h, \|f\|_{C^{1}(\C_{\theta})}$ and  $\|h\|_{C^1(\C_{\theta})}$.
Furthermore, if
 $k-l+1< p\leq k-l+2$, we have $1< a=\frac{p-1}{k-l}\leq 2$, then the constant $C$ in
\eqref{eqn:C2-double-normal-p-large} is independent of $p$. 
This completes the proof.
\end{proof}

As a direct consequence of Lemma \ref{c2 est-admissible soluion-big case} and Lemma \ref{lemma-global-C2-double-normal} (1), we obtain the global $C^2$ estimate of solutions to Eq. \eqref{quo-equ} for $p>k-l+1$. 
\begin{cor}\label{Cor-c2 est-admissible soluion-p-large}
	Let $0\leq l<k\leq n$ and $\theta\in (0, \pi/2)$.  Suppose that  $p>k-l+1$ and $h$ is an admissible solution to Eq. \eqref{quo-equ}. Then
		\begin{eqnarray*}
			\max\limits_{\C_{\theta}}|\n^{2}h|\leq C,
		\end{eqnarray*}
		where $C>0$ depends only on $n, k, l, p,  \min\limits_{\C_{\theta}}f, \min\limits_{\C_{\theta}}h, \|f\|_{C^{2}(\C_{\theta})}$ and  $\|h\|_{C^1(\C_{\theta})}$. Furthermore, if $k-l+1< p\leq k-l+2$,  the constant $C$ is independent of $p$.
	\end{cor}

Next, we consider the case $1< p<k-l+1$ and establish a quantitative estimate for the double normal derivative of the positive, strictly convex solution $h$ to \eqref{quo-equ}. This estimate is a key ingredient in our treatment of the range $1\le p<k-l+1$. Moreover, the estimate obtained here extends our previous work \cite[Lemma~3.5]{MWW-Weingarten} for the case $k=n$ and $p=1$. The endpoint $p=1$ is
handled separately in Lemma~\ref{lem-p-one-C2}. 

\begin{lem}\label{c2 est-convex soluion-small case}
		Let $0\leq l<k\leq n$  and $\theta\in (0, \pi/2)$.  Suppose that  $1< p<k-l+1$ and $h$ is a capillary even, strictly convex solution to Eq. \eqref{quo-equ}
        , then		 \begin{eqnarray*}
        \max\limits_{ \p\C_{\theta}}|\n^{2}h(\mu,\mu)|\leq C\left[1+\|h\|_{C^{0}(\C_{\theta})}\left(1+(\min\limits_{\C_{\theta}} h)^{\frac{p-1}{k-l}-1 }\right)\right]. 
     \end{eqnarray*}  
	where $C>0$ depends only on $n, k, l, p, \min\limits_{\C_{\theta}}f$, and $\|f\|_{C^{1}(\C_{\theta})}$.
	\end{lem}
\begin{proof}

Consider the function
		\begin{eqnarray}\label{eqn:overline-fancy-Q}
			 \overline{\mathcal{Q}}(\xi)\coloneqq\<\n h, \n \zeta\>-\left(\mathcal{B}_1+\frac{M}{2}\right)\zeta(\xi)-\cot\theta h(\xi), \quad  \xi\in \Omega_{\delta_{0}},
		\end{eqnarray}
		where	$\mathcal{B}_1$ is a positive constant to be determined later.
		
		Assume that $ \overline{\mathcal{Q}} $ attains its minimum  at some point $\xi_{0}\in  \Omega_{\delta_{0}}\setminus \partial \Omega_{\delta_{0}}$, and choose an orthonormal frame  $\{e_{i}\}_{i=1}^{n}$ around $\xi_{0}$  such that  $(W_{ij})$ is diagonal at $\xi_{0}$. Using \eqref{third comm}, \eqref{1-time-deri}, and \eqref{hessian of zeta}, at $\xi_0$, we obtain
		\begin{eqnarray}
			0&\leq & 	 F^{ij}  \overline{\mathcal{Q}}_{ij} \notag 
			\\&=&F^{ij}h_{sij}\zeta_{s}+F^{ij}h_{s}\zeta_{sij}+2F^{ij}h_{si}\zeta_{sj}-\left(\mathcal{B}_1+\frac{M}{2}\right)F^{ii}\zeta_{ii}-\cot\theta F^{ii}h_{ii}
			\notag \\&= &F^{ii}\left(W_{iis}-h_{i}\delta_{si}
			\right)\zeta_{s}+2F^{ii}(W_{ii}-h)\zeta_{ii}-\left(\mathcal{B}_1+\frac{M}{2}\right)F^{ii}\zeta_{ii}\notag    \\
        &&+F^{ii}h_{s}\zeta_{sii}-\cot\theta(F(W)-h\sum\limits_{i=1}^{n}F^{ii})
        \notag \\
         &\leq &\zeta_{s} \left(\widehat{f}_{s}h^{a}+ah^{a-1}\widehat{f}h_{s}\right)-F^{ii}h_{i}\zeta_{i}+2F^{ii}W_{ii}\zeta_{ii}+F^{ii}h_{s}\zeta_{sii}\notag \\
         &&+h\cot\theta \sum\limits_{i=1}^{n}F^{ii}-\left(\mathcal{B}_1+\frac{M}{2}\right)F^{ii}\zeta_{ii}.\label{c2 key-1}
\end{eqnarray}
The convexity of $h$ implies
\begin{eqnarray}\label{fw}
    F^{ii}W_{ii}\zeta_{ii}\leq CF(W).
\end{eqnarray}
Combining \eqref{hc1}, \eqref{c2 key-1}, \eqref{fw} and \eqref{hessian of zeta},  we deduce that 
\begin{eqnarray*}
    0\leq F^{ij}\ov{\mathcal{Q}}_{ij}\leq C\left[1+\|h\|_{C^{0}(\C_{\theta})}\left(1+(\min\limits_{\C_{\theta}} h)^{a-1}\right)\right]\left(1+\sum\limits_{i=1}^{n}F^{ii} \right)-\epsilon_0 \left(\mathcal{B}_1+\frac{M}{2} \right)\sum\limits_{i=1}^{n}F^{ii},
\end{eqnarray*}
where $C>0$ depends only on $n, k, l, p, \min\limits_{\C_{\theta}}f$, and $\|f\|_{C^{1}(\C_{\theta})}$.

Choose  
		\begin{eqnarray}\label{chosen of B}
		\mathcal{B}_1&=&
\frac{4C \left( \frac 1 {c_{n,k,l} }
+1 \right)\left[1+\|h\|_{C^{0}(\C_{\theta})}\left(1+(\min\limits_{\C_{\theta}} h)^{a-1}\right)\right]}{\epsilon_0}\notag \\
&&+\frac{1}{1-e^{-\delta_{0}}}\max\limits_{   \C_{\theta}}(|\n h|+\cot\theta  h), ~~~~~~~~
		\end{eqnarray}
this gives $F^{ij}\ov{\mathcal Q}_{ij}<0$ at $\xi_{0}$, which contradicts  $F^{ij}  \overline{\mathcal{Q}}_{ij}\geq 0$ at $\xi_0$. Consequently, we derive that $\xi_0\in \p \O_{\d_0}$.
		
		If  $\xi_0\in \partial \C_{\theta}\cap \p \O_{\d_0}$, it is easy to see that $  \overline{\mathcal{Q}}(\xi_0)=0$. 
		
		If $\xi_0\in \partial \Omega_{\delta_{0}}\setminus  \partial \C_{\theta} $, we have ${d}(\xi_0)=\delta_0$, and from \eqref{chosen of B},
		\begin{eqnarray*}
			 \overline{\mathcal{Q}}(\xi)\geq -|\n h|+\mathcal{B}_1(1-e^{-\delta_{0}})-\cot\theta h\geq 0.
		\end{eqnarray*}
		
		In conclusion, we deduce that 
		\begin{eqnarray*}
			 \overline{\mathcal{Q}}(\xi)\geq 0,\quad{\rm{in}}\quad \Omega_{\delta_{0}}.
		\end{eqnarray*}  
	
		Assume  $h_{\mu\mu}(\eta_{0})\coloneqq\sup\limits_{\partial \C_{\theta}}h_{\mu\mu}>0$ for some $\eta_0\in \p \C_\theta$. In view of  \eqref{hc1} and $ \ov{\mathcal  Q}\equiv 0$ on $\p \C_\theta$, 
		\begin{eqnarray*}
			0&\geq &  \overline{\mathcal{Q}}_{\mu}(\eta_{0})\\
			&\geq &(h_{s\mu}\zeta_{s}+h_{s}\zeta_{s\mu})-\left(\mathcal{B}_1+\frac{M}{2}\right)\zeta_{\mu}-\cot\theta h_{\mu}\\
			&=& h_{\mu\mu}(\eta_{0})-\left(\mathcal{B}_1+\frac{M}{2}\right)+h_{s}\zeta_{s\mu}-\cot^{2}\theta h,
		\end{eqnarray*}
		which yields  
		\begin{eqnarray}\label{sup estimate}
			\max\limits_{\partial \C_{\theta}}h_{\mu\mu}
            \leq C \left[1+\|h\|_{C^{0}(\C_{\theta})}\left(1+(\min\limits_{\C_{\theta}} h)^{a-1}\right)\right]            +\frac{M}{2},
		\end{eqnarray}
  where $C$ is a  constant depending on $n, k, l, p,   \min\limits_{\C_{\theta}}f$ and    $\|f\|_{C^1(\C_{\theta})}$.

  Since $W$ is positive definite, we have 
  \begin{eqnarray}
      h_{\mu\mu}>-h\geq -\|h\|_{C^{0}(\C_{\theta})}.\label{inf estimate}
  \end{eqnarray}
Therefore, \eqref{sup estimate} and \eqref{inf estimate} together yield
		 \begin{eqnarray*}\max\limits_{ \p\C_{\theta}}|\n^{2}h(\mu,\mu)|\leq C\left[1+\|h\|_{C^{0}(\C_{\theta})}\left(1+(\min\limits_{\C_{\theta}} h)^{a-1}\right)\right].
     \end{eqnarray*}This completes the proof.
     \end{proof}

  As a direct consequence of Lemma \ref{c2 est-convex soluion-small case} and Lemma \ref{lemma-global-C2-double-normal} (2), we obtain the global $C^2$ estimate of solutions to Eq. \eqref{quo-equ} for $1\leq p<k-l+1$. 
\begin{cor}\label{cor-global-c2-p-small case}
Let $0\leq l<k\leq n$ and $\theta\in(0,\pi/2)$.
Suppose that $1<p<k-l+1$ and $h$ is a positive, capillary even,
strictly convex solution of Eq. \eqref{quo-equ}. Then
\begin{eqnarray*}
\max\limits_{\C_\theta}|\n^2h|
\leq C\left[1+\|h\|_{C^0(\C_\theta)}
\left(1+(\min\limits_{\C_\theta}h)^{\frac{p-1}{k-l}-1}\right)
+(\max\limits_{\C_\theta}h)^{2-\gamma}
(\min\limits_{\C_\theta}h)^{\gamma-1}\right],
\end{eqnarray*}
where $\gamma$ is chosen as in Lemma~\ref{lem weighted gra},
and $C>0$ depends only on $n,k,l,p,\theta,\gamma$,
$\min\limits_{\C_\theta}f$, and $\|f\|_{C^2(\C_\theta)}$.
\end{cor}

We next treat the endpoint $p=1$ separately.
\begin{lem}\label{lem-p-one-C2}
Let $0\leq l<k\leq n$ and $\theta\in(0,\pi/2)$.
Suppose that $h$ is a positive, capillary even, strictly convex
solution of Eq. \eqref{quo-equ} with $p=1$. Then
\begin{eqnarray}\label{eqn-p-one-C2}
\max\limits_{\C_\theta}|\n^2h|
\leq C\left(1+\|h\|_{C^0(\C_\theta)}\right),
\end{eqnarray}
where $C>0$ depends only on $n,k,l,\theta$,
$\min\limits_{\C_\theta}f$, and $\|f\|_{C^2(\C_\theta)}$.
\end{lem}
\begin{proof}
The estimate \eqref{hc1} gives
$\max\limits_{\C_\theta}|\n h|\leq C_\theta \|h\|_{C^0(\C_\theta)}$.
Choose a fixed region $\O_{\delta_0}$ on which
$\zeta=e^{-d}-1$ is smooth and
$\n^2\zeta\geq\epsilon_0\s$, as in
\eqref{hessian of zeta}. The constants $\delta_0,\epsilon_0>0$
depend only on the geometry of $\C_\theta$.

Consider
\begin{eqnarray}\label{eqn:underline-fancy-Q}
\overline{\mathcal Q}
\coloneqq\<\n h,\n\zeta\>
-\left(\mathcal B_1+\frac M2\right)\zeta-\cot\theta h.
\end{eqnarray}
When $p=1$, we have $a=0$ and $F(W)=\widehat f$.
In particular, the term $ah^{a-1}\widehat f h_s$ in
\eqref{c2 key-1} vanishes. By \eqref{fw},
\begin{eqnarray*}
\sum\limits_{i=1}^nF^{ii}W_{ii}\zeta_{ii}
\leq CF(W)=C\widehat f.
\end{eqnarray*}
Thus the computation in \eqref{c2 key-1} yields
\begin{eqnarray}\label{eqn-p-one-barrier}
F^{ij}\overline{\mathcal Q}_{ij}
\leq C_0(1+\|h\|_{C^0(\C_\theta)})\left(1+\sum\limits_iF^{ii}\right)
-\epsilon_0\left(\mathcal B_1+\frac M2\right)
\sum\limits_iF^{ii}.
\end{eqnarray}
By Proposition~\ref{pro-2.3} (4),
$\sum\limits_iF^{ii}\geq c_{n,k,l}>0$.
Choose $\mathcal B_1\coloneqq A(1+\|h\|_{C^0(\C_\theta)})$, with
\begin{eqnarray*}
A>\frac{2C_0(1+c_{n,k,l}^{-1})}{\epsilon_0},
\qquad
A>\frac{C_\theta+\cot\theta}{1-e^{-\delta_0}}+1.
\end{eqnarray*}
Then the right-hand side of \eqref{eqn-p-one-barrier} is negative.
Moreover, $\overline{\mathcal Q}=0$ on $\p\C_\theta$, while
on the inner boundary $\{d=\delta_0\}$ we have
\begin{eqnarray*}
\overline{\mathcal Q}
\geq-|\n h|-\cot\theta h
+\mathcal B_1(1-e^{-\delta_0})\geq0.
\end{eqnarray*}
The minimum principle gives
$\overline{\mathcal Q}\geq0$ in $\O_{\delta_0}$.

Differentiating at $\p\C_\theta$ in the outward conormal
direction, and using $\n\zeta=\mu$ there, we obtain
\begin{eqnarray*}
0\geq\n_\mu\overline{\mathcal Q}
&=&h_{\mu\mu}-\mathcal B_1-\frac M2
+h_s\zeta_{s\mu}-\cot^2\theta h.
\end{eqnarray*}
Consequently,
\begin{eqnarray*}
\max\limits_{\p\C_\theta}h_{\mu\mu}
\leq C(1+\|h\|_{C^0(\C_\theta)})+\frac M2.
\end{eqnarray*}
Strict convexity also gives $h_{\mu\mu}>-h\geq-\|h\|_{C^0(\C_\theta)}$.
These inequalities imply $M\leq C(1+\|h\|_{C^0(\C_\theta)})$.
Substituting this bound into \eqref{uniform C2-p-one}
proves \eqref{eqn-p-one-C2}.
\end{proof}

\subsection{Proof of Theorem \ref{thm priori est}}

We now complete the proof of Theorem \ref{thm priori est}.
\begin{proof}[\textbf{Proof of Theorem \ref{thm priori est}}]
We first prove part \eqref{thm-C2-est-p-large}.   Suppose that $p>k-l+1$. Combining Lemma \ref{lem kC0} and Lemma \ref{lem c1} with Corollary \ref{Cor-c2 est-admissible soluion-p-large}, we obtain 
   \begin{eqnarray*}
       \min\limits_{\C_{\theta}} h \geq c,\quad {\text{and}}\quad \|h\|_{C^{2}(\C_{\theta})}\leq C.
   \end{eqnarray*}
The higher order estimate in \eqref{thm-C2}  then follows from the theory of fully nonlinear uniformly elliptic equations with oblique boundary conditions,  together with the Schauder estimates (see, e.g., \cite[Theorem 1.1]{LT}, or \cite[Section 3]{LTU}).

\ 

Next, we prove part~\eqref{thm-C2-est-p-middle}.
The rescaled function $\widetilde h$ defined there satisfies

\begin{eqnarray*}
\left\{\begin{array}{rcll}\vspace{2mm}\displaystyle
	  \frac{\sigma_{k}(\n^{2} \wt h+ \wt h\sigma)}{\sigma_{l}(\n^{2}\wt h+\wt h\s)}&=&f \left (\min\limits_{\C_{\theta}}h \right)^{p+l-k-1} \wt h^{p-1},& \quad    \hbox{ in } \C_{\theta},\\			\n_\mu \wt h &=& \cot\theta \wt h, &\quad  \hbox{ on } \p \C_\theta.
		\end{array}\right.
	\end{eqnarray*}  
Denote $$\wt{g}\coloneqq f \left(\min\limits_{\C_{\theta}}h\right)^{p+l-k-1}.$$ By Lemma~\ref{lem kC0} (1) and $p\in (k-l+1, k-l+2]$, we have \begin{eqnarray}    \frac{\binom{n}{k}}{\binom{n}{l}}\frac{1-\cos\theta}{\max\limits_{\C_{\theta}}f \cdot (\sin\theta)^{2(k-l)}}\leq  \left(\min_{\C_\theta} h \right)^{p+l-k-1}\leq \frac{\binom{n}{k}}{\binom{n}{l}}\frac{1}{ \min\limits_{\C_{\theta}}f\cdot (1-\cos\theta)^{k-l+1}}.   \label{min-con}\end{eqnarray} then \eqref{min-con} implies
\begin{eqnarray*}
    \min\limits_{\C_{\theta}}\wt g\geq C, \quad \text{and}\quad \|\wt g\|_{C^{2}(\C_{\theta})}\leq C^{-1},
\end{eqnarray*}
where $C>0$ depends only on $n, k,l, \min\limits_{\C_{\theta}}f$ and $\|f\|_{C^{2}(\C_{\theta})}$, and is independent of $p$.

Moreover, Lemma \ref{lem kC1} implies that
\begin{eqnarray*}
    1\leq \wt{h}\leq C, \quad {\rm{and}}\quad |\n \wt h|\leq C.
\end{eqnarray*}
Applying Corollary~\ref{Cor-c2 est-admissible soluion-p-large} to $\widetilde h$ yields 
\begin{eqnarray*}
    \|\wt h\|_{C^2(\C_\theta)}\leq C.
\end{eqnarray*} The higher order estimate \eqref{re-est} then follows as in part~\eqref{thm-C2-est-p-large}.

\ 
Finally, we prove part~\eqref{thm-C2-est-p-small}. 
Since $0\leq p-1<k-l$, the extremum comparisons in the proof of
Lemma~\ref{lem kC0} yield
\begin{eqnarray*}
(\min\limits_{\C_\theta}u)^{k-l-(p-1)}
\leq\frac{\max\limits_{\C_\theta}(f\ell^{p-1})}{c^{k-l}_{n,k,l}},
\qquad
(\max\limits_{\C_\theta}u)^{k-l-(p-1)}
\geq\frac{\min\limits_{\C_\theta}(f\ell^{p-1})}{c^{k-l}_{n,k,l}}.
\end{eqnarray*}
These comparisons remain valid at boundary extrema because
$\n_\mu u=0$: the gradient vanishes there and the full Hessian
has the appropriate sign. In particular,
\begin{eqnarray*}
\min\limits_{\C_\theta}h\leq C_0,
\qquad \max\limits_{\C_\theta}h\geq c_0>0.
\end{eqnarray*}
 
Capillary evenness allows the comparison spherical caps to be centred
on the symmetry axis. Indeed, reflecting a comparison cap
in this axis gives a comparison spherical cap with the opposite horizontal centre, and convexity allows these centres to be replaced by their midpoint.
Together with the radius comparisons, this yields
\begin{eqnarray*}
c_\theta \rho_-(\widehat\S,\theta) \leq\min\limits_{\C_\theta}h\leq C_\theta \rho_-(\widehat\S,\theta),
\qquad
c_\theta \rho_+(\widehat\S,\theta) \leq\max\limits_{\C_\theta}h\leq C_\theta \rho_+(\widehat\S,\theta).
\end{eqnarray*}
Consequently,
\begin{eqnarray}\label{small key-1}
\rho_-(\widehat\S,\theta) \leq C_1,
\end{eqnarray}
and
\begin{eqnarray}\label{small key-2}
\rho_+(\widehat\S,\theta)\geq c_1>0.
\end{eqnarray}

\noindent\textbf{Case 1} $p=1$.
Lemma~\ref{lem-p-one-C2} gives
\begin{eqnarray*}
\max\limits_{\C_\theta}\lambda_{\max}(W)\leq C(1+\rho_+(\widehat\S,\theta)).
\end{eqnarray*}
By Lemma~\ref{chou-wang-lemma},
\begin{eqnarray}\label{eqn-p-one-radius}
\frac{\rho_+(\widehat\S,\theta)^2}{\rho_-(\widehat\S,\theta)}\leq C(1+\rho_+(\widehat\S,\theta)).
\end{eqnarray}
Using \eqref{small key-1}, we obtain $\rho_+(\widehat\S,\theta)\leq C_2$.
Then \eqref{small key-2} and \eqref{eqn-p-one-radius} give
\begin{eqnarray*}
\rho_-(\widehat\S,\theta) \geq\frac{\rho_+(\widehat\S,\theta)^2}{C(1+\rho_+(\widehat\S,\theta))}\geq c_2>0.
\end{eqnarray*}
The radius comparisons imply $0<c\leq h\leq C$.
Applying Lemma~\ref{lem-p-one-C2} once more gives
$\|h\|_{C^2(\C_\theta)}\leq C$.

\noindent\textbf{Case 2} $1<p<k-l+1$.
Recall that $a$ is \eqref{eqn-defn-a}, thus $a\in (0, 1)$.
Corollary~\ref{cor-global-c2-p-small case} and
Lemma~\ref{chou-wang-lemma} give
\begin{eqnarray*}
\frac{\rho_+(\widehat\S,\theta)^2}{\rho_-(\widehat\S,\theta)}
\leq C\left[1+\rho_+(\widehat\S,\theta)(1+\rho_-(\widehat\S,\theta)^{a-1})+\rho_+(\widehat\S,\theta)^{2-\gamma}\rho_-(\widehat\S,\theta)^{\gamma-1}\right].
\end{eqnarray*}
Multiplying by $\rho_-(\widehat\S,\theta)$ and using \eqref{small key-1}, we obtain
\begin{eqnarray*}
\rho_+(\widehat\S,\theta)^2\leq C\left(1+\rho_+(\widehat\S,\theta)+\rho_+(\widehat\S,\theta)^{2-\gamma}\right).
\end{eqnarray*}
Since $0<\gamma<2$, this bounds $\rho_+(\widehat\S,\theta)$ from above.
Using \eqref{small key-2} in the same inequality before
discarding the powers of $r$, we find
\begin{eqnarray*}
1\leq C\left(\rho_-(\widehat\S,\theta)+\rho_-(\widehat\S,\theta)^a+\rho_-(\widehat\S,\theta)^\gamma\right).
\end{eqnarray*}
As $a,\gamma>0$, this bounds $\rho_-(\widehat\S,\theta)$ away from zero.
The radius comparisons and
Corollary~\ref{cor-global-c2-p-small case} again give
$0<c\leq h\leq C$ and $\|h\|_{C^2(\C_\theta)}\leq C$.

In both cases, the equation and the Newton--Maclaurin
inequalities imply
\begin{eqnarray*}
0<c\leq fh^{p-1}
=\frac{\sigma_k(W)}{\sigma_l(W)}
\leq C\sigma_k(W)^{\frac{k-l}{k}},
\end{eqnarray*}
so $\sigma_k(W)\geq c'>0$.
Together with $W\geq0$ and $|W|\leq C$, this confines the
eigenvalues of $W$ to a compact subset of $\Gamma_k$.
The equation is therefore uniformly elliptic.
The global oblique boundary estimates and Schauder estimates
used in part~\eqref{thm-C2-est-p-large} yield the required
higher order bounds. For fixed $p$, $\theta$ and H\"older
exponent, the constants depend only on $n,k,l$ and the
corresponding positive lower bound and smooth norms of $f$.
This completes the proof of Theorem~\ref{thm priori est}.
\end{proof}

\section{Convexity}\label{sec-4}

For $1\leq k\leq n-1$, admissibility alone does not imply
strict convexity.
We show that a convex solution is strictly convex under
the stated structural conditions on $f$.
The interior argument follows \cite[Lemma~5]{HMS2004}
and \cite[Theorem~3.5]{MWW-Quotient},
while the treatment of boundary degeneracy uses the
boundary-point argument of 
\cite[Theorem~3.1]{HIS2025}.
For $k=n$, strict convexity already follows from admissibility.
 
\begin{thm}\label{thm-convex}
    Let $0\leq l< k\leq n-1$ and $\theta \in (0, \pi/2)$. Suppose that $p\geq 1$ and $f$ is a positive smooth function satisfying 
    \begin{eqnarray}\label{f-convex}
        \n^{2}f^{-\frac{1}{p+k-l-1}}+f^{-\frac{1}{p+k-l-1}}\sigma \geq 0,  \quad {\text{in}}\quad \C_{\theta},
    \end{eqnarray}
    and 
    \begin{eqnarray}\label{f-bry}
        \n_{\mu} f+\cot\theta (p+k-l-1)f\geq 0,\quad {\text{on}}\quad \partial\C_{\theta}.
    \end{eqnarray}
Let $h\in C^{4}(\mathcal{C}_{\theta})$ be a positive admissible solution to \eqref{quo-equ}. If the associated matrix $$W=\n^2 h+h\s $$ is positive semi-definite, then \(W\) is positive definite. Equivalently, the solution \(h\) is strictly convex.
\end{thm}

\begin{proof}
 We argue by contradiction. Suppose that \(W\) is attains its minimum rank $r$ at point $\xi_{0}\in\C_{\theta}$ with $k\leq r\leq n-1$, then there exist a small neighborhood $\mathcal{O}\subset \C_{\theta}$ of $\xi_{0}$ and a small constant $c_{0}$, such that 
    \begin{eqnarray}\label{eqn:sigma-r-W}
        \sigma_{r}(W)(\xi)>c_{0},~\text{in}~\mathcal{O}, \end{eqnarray}
        and \begin{eqnarray} \label{zero}
        \sigma_{r+1}(W)(\xi_{0})=0.
    \end{eqnarray}
    We define the  function 
\begin{eqnarray*}
    \eta(\xi)\coloneqq  \left\{
\begin{array}{llll} 
\sigma_{r+1}(W)+\frac{\sigma_{r+2}(W)}{\sigma_{r+1}(W)}, & \text{ if } x\in\mathcal{N}, \\
0 , & \text{ if } x\in \mathcal{O}\setminus\mathcal{N},
\end{array} \right.
\end{eqnarray*}
where $\mathcal{N}\coloneqq \left\{x\in\mathcal{O}: \s_{m+1}(W)(x)>0\right\}$. 
The local regularity result of Bian-Guan \cite[Corollary~2.2]{BG} gives $\eta\in C^{1.1}_{\text{loc}} \left(\mathcal{O}\cap \text{int}(\mathcal C_\theta) \right)$. 
Moreover, Newton--Maclaurin inequality of Proposition \ref{pro-2.3} (1) and \eqref{eqn:sigma-r-W} imply
\begin{eqnarray*}
   0\le\eta\le(1+c_0^{-1})\sigma_{r+1}(W), 
\end{eqnarray*}
so the indicated zero extension is continuous in \(\mathcal O\).

\

 \noindent\textbf{Case 1} $\xi_{0}\in \C_{\theta}\setminus \partial\C_{\theta}$. 
 By \cite[Theorem 3.5, (3.41)]{MWW-Quotient} and the structural condition \eqref{f-convex}, there exists a positive constant $C$ depending on $n, k,l, \|h\|_{C^{3}(\C_{\theta})}, \min\limits_{\C_{\theta}}f$ and $\|f\|_{C^{2}(\C_{\theta})}$, such that  
  \begin{eqnarray*}
   F^{ij}\n^{2}_{ij}\eta-C(|\n \eta|+\eta)\leq 0,~~~~ {\rm{ in }}~~~ \mathcal O.
 \end{eqnarray*}
Since $\eta\ge0$ and $\eta(\xi_0)=0$, the strong maximum
principle yields
\begin{eqnarray*}
\eta\equiv0\qquad\text{in }\mathcal O.
\end{eqnarray*}
Consequently, $W$ has constant rank $r$ in $\mathcal O$.

Let
\begin{eqnarray*}
\mathcal Z\coloneqq  \left\{\xi\in\operatorname{int}\mathcal C_\theta:
                  \operatorname{rank}W(\xi)=r \right\}.
\end{eqnarray*}
The preceding interior argument shows that \(\mathcal Z\) is open in
\(\operatorname{int}\mathcal C_\theta\). Continuity of \(W\)
and the minimality of \(r\) show that it is also closed there.
Since \(\operatorname{int}\mathcal C_\theta\) is connected and
\(\xi_0\in \mathcal Z\), we conclude that
\begin{eqnarray*}
\mathcal Z=\operatorname{int}\mathcal C_\theta.
\end{eqnarray*}
By continuity, \(\sigma_{r+1}(W)=0\) also on the boundary;
the global minimality of \(r\) implies
\(\sigma_r(W)>0\) there. 
Thus
\begin{eqnarray}\label{eqn:sigma-r+1-W=0}
\sigma_{r+1}(W)\equiv0
\quad\text{and}\quad
\sigma_r(W)>0
\quad\text{in }\mathcal C_\theta.
\end{eqnarray}

To justify the capillary Minkowski formula when $W$ is only positive
semi-definite, apply it first to $h_\varepsilon=h+\varepsilon\ell$.
The Robin boundary condition is preserved and
\begin{eqnarray*}
\n^2h_\varepsilon+h_\varepsilon\s=W+\varepsilon \s>0.
\end{eqnarray*}
Passing to the limit as $\varepsilon\downarrow0$ gives the identity
for $h$.

On the other hand, the capillary Minkowski formula (see, e.g., \cite[Corollary~2.10]{MWWX} or \cite[Proposition 2.6]{WWX-MA2024}) asserts that
    \begin{eqnarray*}
      (n-r)\int_{\C_{\theta}}h(\xi)\sigma_{r}(W)d\s=(r+1)\int_{\C_{\theta}}\ell(\xi)\sigma_{r+1}(W)d\s.
    \end{eqnarray*}
By \eqref{eqn:sigma-r+1-W=0}, the right-hand side vanishes. Therefore,
\begin{eqnarray*}    \int_{\C_{\theta}}h(\xi)\sigma_r(W)\,d\sigma=0.
\end{eqnarray*}
This is impossible, since \(h>0\) and \(\sigma_r(W)>0\) in \(\C_{\theta}\), by \eqref{eqn:sigma-r-W}. 

\

\noindent\textbf{Case 2} $\xi_{0}\in \partial \C_{\theta}$.  Let $\{e_{i}\}_{i=1}^{n}$ be an orthonormal frame around $\xi_0$ such that $e_{n}=\mu$. For any $\xi\in \partial\C_{\theta}$,  the boundary condition implies (cf. \cite[Proposition~2.8]{MWW-AIM})
    \begin{eqnarray}\label{h1n}
        \n^{2}h(e_{\alpha}, e_{n})=\n_{e_{\a}}(h_{n})-\<\n h, \n_{e_{\alpha}}e_{n}\>=0,~~~~1\leq \a \leq n-1. 
    \end{eqnarray}
    Hence, $W_{\alpha n}=0$ on $\partial\C_{\theta}$. By further diagonalizing the tangential block of $W$, we may assume that $W$ is diagonal at $\xi_{0}$. 
    
    We first \textbf{claim} that
\begin{eqnarray}\label{claim}
    W_{nn}(\xi_{0})=0.
\end{eqnarray}

Indeed, suppose, by contradiction, that $W_{nn}(\xi_{0})>0$. Together with \eqref{zero}, this implies there exists $\alpha$, $(1\leq \a\leq n-1)$, such that   $$W_{\alpha\alpha}(\xi_{0})=0.$$ 
Using the Gauss-Weingarten formula of $\p \C_\theta \subset \C_\theta$ and combining with \eqref{h1n}, we have 
\begin{eqnarray*}
			h_{n\alpha \alpha}&=& \n_{e_\a}(\n^2h(e_\alpha,e_n))-\n^2h(\n_{e_\a}e_\a,e_n)-\n^2h(e_\a,\n_{e_\a}e_n)
			\\&=&\cot\theta h_{nn}-\cot\theta h_{\a\a}.
		\end{eqnarray*} 
By \eqref{third comm}, we obtain  
\begin{eqnarray} \label{Waan}  
   W_{\alpha\alpha n}=\cot\theta (W_{nn}-W_{\a\a}).
\end{eqnarray}
Since $W_{\alpha\alpha}\geq0$ in $\C_\theta$ and
$W_{\alpha\alpha}(\xi_0)=0$, its one-sided outward conormal
derivative satisfies
\begin{eqnarray*}
 0\geq W_{\alpha\alpha n}(\xi_{0})=\cot\theta\left(W_{nn}(\xi_{0})-W_{\a\a}(\xi_{0})\right),
\end{eqnarray*}
which implies $$W_{nn}(\xi_{0})\leq0,$$ contradicting the assumption
$W_{nn}(\xi_{0})>0$. This proves \eqref{claim}.

Differentiating Eq. \eqref{F-equ} in the normal direction yields
\begin{eqnarray*}
    \sum\limits_{\alpha=1}^{n-1}F^{\alpha \alpha}W_{\a\a n}+F^{nn}W_{nnn}=\n_{\mu}(\widehat{f}h^{a}),
\end{eqnarray*}
thus
\begin{eqnarray*}
    W_{nnn}=\frac{1}{F^{nn}}\left(\n_{\mu}(\widehat{f}h^{a})-\sum\limits_{\a=1}^{n-1}F^{\a\a}W_{\a\a n}\right).
\end{eqnarray*}
Let $$Q(W)\coloneqq\sigma_{r+1}(W)~~~\text{ and }~~~  Q^{ij}\coloneqq\frac{\partial Q}{\partial W_{ij}}.$$  
Using \eqref{zero}, \eqref{claim}, and \eqref{Waan}, we deduce
\begin{eqnarray}\label{qmu}
\n_{\mu}Q&=&Q^{nn}W_{nnn}+\sum\limits_{\a=1}^{n-1}Q^{\alpha\alpha}W_{\alpha\alpha n}   \notag \\
    &=&\frac{Q^{nn}}{F^{nn}}\n_{\mu}(\widehat{f}h^{a})-\frac{\cot\theta }{F^{nn}}\sum\limits_{\alpha=1}^{n-1}(F^{\alpha\alpha}Q^{nn}-Q^{\alpha\alpha}F^{nn})(W_{nn}-W_{\a\a})\notag \\
    &=&\frac{Q^{nn}}{F^{nn}} \left(\n_{\mu}(\widehat{f}h ^{a})+\cot\theta \widehat f h^{a}\right).\notag 
\end{eqnarray}
By the boundary condition \eqref{f-bry},
\begin{eqnarray}\label{f-c}
    \n_{\mu}(\widehat{f }h^{a})+\cot\theta \widehat{f}h^{a}= \frac{f^{\frac{1}{k-l}-1}h^{a}}{k-l}\left( \n_{\mu}f+\cot\theta (p+k-l-1)f
    \right)\geq 0,
\end{eqnarray} this implies
\begin{eqnarray}\label{Q}
    \n_\mu Q(\xi_0)\geq 0.
\end{eqnarray}

Since \(Q\ge0\) in \(\mathcal C_\theta\) and
\(Q(\xi_0)=0\), its one-sided outward conormal derivative
satisfies \(\nabla_\mu Q(\xi_0)\le0\). Together with \eqref{Q},
this gives
\begin{eqnarray*}
\nabla_\mu Q(\xi_0)=0.
\end{eqnarray*}
Newton--Maclaurin inequality of Proposition \ref{pro-2.3} (1) and \eqref{eqn:sigma-r-W} imply, wherever \(Q>0\),
\begin{eqnarray*}
0\le\frac{\sigma_{r+2}(W)}{\sigma_{r+1}(W)}
\le\frac{\sigma_{r+1}(W)}{\sigma_r(W)}
\le c_0^{-1}Q.
\end{eqnarray*}
Here we use the convention that \(\sigma_{n+1}=0\) when \(r=n-1\). Therefore, with the
quotient extended by zero where \(Q=0\),
\begin{eqnarray*}
0\le\eta\le(1+c_0^{-1})Q.
\end{eqnarray*}
Let \(\gamma(t)\coloneqq \exp_{\xi_0}(-t\mu)\) for sufficiently small
\(t\ge0\). Since \(Q\in C^{1,1}\), \(Q(\xi_0)=0\), and
\(\nabla_\mu Q(\xi_0)=0\), we obtain
\begin{eqnarray}
0\le\eta(\gamma(t))
\le(1+c_0^{-1})Q(\gamma(t))=O(t^2).
\label{phi-bry}
\end{eqnarray}

The argument in Case 1 excludes any interior point of rank
\(r\); hence \(Q>0\) in the interior. The same local
constant-rank differential inequality can be written there as
\begin{eqnarray*}
F^{ij}\nabla_{ij}^2\eta+b^i\nabla_i\eta-C\eta\le0,
\qquad
b\coloneqq 
\begin{cases}
-C\,\dfrac{\nabla\eta}{|\nabla\eta|},
   & |\nabla\eta|>0,\\   
0, & |\nabla\eta|=0.
\end{cases}
\end{eqnarray*}
For the fixed smooth admissible solution, \((F^{ij})\) is
uniformly elliptic and \(b\) is bounded. Applying the Hopf
boundary point argument in an interior ball tangent to
\(\partial\mathcal C_\theta\) at \(\xi_0\) yields
\begin{eqnarray*}
\liminf_{t\downarrow0}\frac{\eta(\gamma(t))}{t}>0,
\end{eqnarray*}
contradicting \eqref{phi-bry}. This completes the proof of Theorem~\ref{thm-convex}.
\end{proof}

\section{Proof of main Theorems}\label{sec-5}

\subsection{Proof of Theorem \ref{thm-admissible-existence}} \
In this subsection,  we adapt the argument of \cite[Section~3.2]{MWW-Quotient} to prove Theorem \ref{thm-admissible-existence}.  While most of the arguments are analogous to those in \cite[Section~3.2]{MWW-Quotient}, the presence of the Robin boundary condition requires several additional modifications. We begin with the uniqueness part of Theorem \ref{thm-admissible-existence} in the case $p>k-l+1$.
\begin{thm}\label{thm-unique}
Let $0\leq l<k\leq n$ and  $\theta\in (0, \pi)$. Suppose that $p>k-l+1$. Then the positive admissible solution to Eq. \eqref{quo-equ}  is unique.
\end{thm}
\begin{proof}
We argue by contradiction.  Suppose that there exist two positive, admissible solutions $h_{1}$ and $h_{2}$ to Eq. \eqref{quo-equ}. Set $$\mathcal{R}\coloneqq \log \frac{h_{1}}{h_{2}}.$$ It then follows from \eqref{quo-equ} that
\begin{eqnarray}\label{eqn-R-bdry-cond}
    \n_{\mu}\mathcal{R}=0,~~~~\text{ on } \p\C_\theta.
\end{eqnarray}  Assume that $\mathcal{R}$ attains its minimum value at some point $\xi_0\in \C_{\theta}$. If $\xi_0\in \C_\theta\setminus \p\C_\theta$, then at $\xi_{0}$, we have
\begin{eqnarray*}
    0=\n_i\mathcal{R}=\frac{\n_{e_{i}}h_{1}}{h_{1}}-\frac{\n_{e_{i}}h_{2}}{h_{2}}, ~~~1\leq i\leq n,
\end{eqnarray*}
and
\begin{eqnarray*}
    0\leq \left(\n^2\mathcal{R}(e_i,e_j) \right)= \left( \frac{\n^{2}h_{1}(e_{i}, e_{j})}{h_{1}}-\frac{
    \n^{2}h_{2}(e_{i}, e_{j})}{h_{2}} \right), ~~1\leq i,j\leq n.
\end{eqnarray*}
Consequently,
\begin{eqnarray}\label{eqn-inequality-W-1-W-2}
    \frac{W_{1}}{h_{1}}\geq \frac{W_{2}}{h_{2}},
\end{eqnarray}
where $$W_s\coloneqq \n^{2}h_s+h_{s}\s, ~~~~~ s=1,2.$$

If instead $\xi_{0} \in \partial \C_{\theta}$, then \eqref{eqn-R-bdry-cond} yields the same inequality \eqref{eqn-inequality-W-1-W-2} at $\xi_{0}$. Therefore, in either case, we obtain
\begin{eqnarray*}
    h_{1}^{p-1}f(\xi_{0})=\frac{\sigma_{k}(W_1)}{\sigma_{l}(W_1)}\geq \left(\frac{h_{1}}{h_{2}} \right)^{k-l}\frac{\sigma_{k}(W_2)}{\sigma_{l}(W_2)}=h_{1}^{k-l}h_{2}^{p+l-k-1}f(\xi_{0}). 
\end{eqnarray*}
Since $p>k-l+1$, we deduce that $$h_{1}\geq h_{2}.$$ By reversing the roles of $h_{1}$ and $h_{2}$, we similarly obtain $$h_{1} \leq h_{2}.$$ Hence, $h_{1}\equiv h_{2}$ in $\C_{\theta}$.   
\end{proof}

In particular, if $$f=\frac{\binom{n}{k}}{\binom{n}{l}} \ell^{1-p},$$ in Eq. \eqref{quo-equ}, we derive the following rigidity result.
\begin{cor}\label{cor:5.2-uniqueness}
    Let $0\leq l<k\leq n$ and $\theta\in (0, \pi)$. Suppose that $p>k-l+1$, and let $h$ be an admissible solution of Eq.
\begin{eqnarray}\label{Hkl-eq}
\left\{\begin{array}{rcll}\vspace{2mm}\displaystyle
\frac{\sigma_{k}(W)}{\sigma_{l}(W)}&=&\frac{\binom{n}{k}  }{\binom{n}{l}}\left(\frac{h}{\ell}\right)^{p-1},& \quad    \text{ in }~ \C_{\theta},\\			\n_\mu h &=& \cot\theta h, &\quad  \text{ on } ~ \p \C_\theta.
		\end{array}\right.
	\end{eqnarray}
    Then $h=\ell$.
\end{cor}
If $h$ is strictly convex, then Eq. \eqref{Hkl-eq} is equivalent to the following prescribed $L^{p}$ curvature problem for a strictly convex capillary hypersurface 
\begin{eqnarray}\label{Lp curvature pro} 
    \frac{\sigma_{n-l}(\kappa)}{\s_{n-k}(\kappa)}=\frac{\binom{n}{l}}{\binom{n}{k}}\left(\frac{\ell}{h}\right)^{p-1}=\frac{\binom{n}{l}}{\binom{n}{k}}u^{1-p}. 
\end{eqnarray}
When $p\geq 1$, Gao--Li \cite[Theorem 1.6]{GLi} proved that if a strictly convex capillary hypersurface $\S\subset\ol{\RR^{n+1}_+}$ satisfies Eq. \eqref{Lp curvature pro}, then $\S$ must be a spherical cap.
In particular, in the case $p=1$ and $k=n$, Jia--Wang--Xia--Zhang \cite[Corollary 1.2]{JWXZ}  showed that the same conclusion remains valid for $(n-l)$-admissible embedded capillary hypersurfaces, by employing a capillary version of the Heintze--Karcher inequality. In what follows, we establish an analogous rigidity result for admissible solutions of \eqref{Lp curvature pro} for general $0\leq l<k\leq n$ and $p=1$, which may be of independent interest.  Our argument is inspired by the approach used by Guan--Ma--Zhou in \cite[Proposition 3.1 and its proof]{gmz2006} in the context of the closed hypersurface setting.

\begin{thm}\label{thm-constant-p-one}
Let $0\leq l<k\leq n$ and $\theta\in(0,\pi/2)$.
Suppose that $h\in C^4(\C_\theta)$ is an admissible solution of
\begin{eqnarray*}
\left\{\begin{array}{rcll}\vspace{2mm}\displaystyle
\frac{\sigma_k(\n^2h+h\s)}{\sigma_l(\n^2h+h\s)}
&=&\frac{\binom{n}{k}}{\binom{n}{l}},
&\quad\hbox{in }\C_\theta,\\
\n_\mu h&=&\cot\theta h,
&\quad\hbox{on }\p\C_\theta.
\end{array}\right.
\end{eqnarray*}
Then
\begin{eqnarray}\label{h-uni}
h=\ell+\sum\limits_{i=1}^n a_i\<\xi,E_i\>,
\end{eqnarray}
where $a_i\in\RR$ and $E_1,\ldots,E_n$ are the horizontal
coordinate unit vectors. In particular, if $h$ is capillary
even, then $h=\ell$.
\end{thm}
\begin{proof}
Let $F$ be as in \eqref{F-equ}, and set
\begin{eqnarray*}
\bar H\coloneqq\Delta h+nh=\sum\limits_{i=1}^nW_{ii}.
\end{eqnarray*}
Since $F(W)=F(I)$, concavity at $I$ gives
\begin{eqnarray*}
F(I)=F(W)
&\leq&F(I)+F^{ij}(I)(W_{ij}-\delta_{ij})\\
&=&F(I)+\frac{F(I)}n(\bar H-n).
\end{eqnarray*}
Consequently,
\begin{eqnarray}\label{H-2}
\bar H\geq n.
\end{eqnarray}
Concavity at $W$ and homogeneity also give
\begin{eqnarray*}
\sum\limits_{i=1}^nF^{ii}(W)\geq F(I).
\end{eqnarray*}
Differentiating $F(W)=F(I)$ twice and using concavity, we obtain
\begin{eqnarray*}
F^{ij}W_{ijss}
=-F^{ij,pq}W_{ijs}W_{pqs}\geq0,
\end{eqnarray*}
where repeated indices are summed. The commutation identities
therefore yield, throughout the interior,
\begin{eqnarray}\label{uni-1}
F^{ij}\bar H_{ij}
&=&F^{ij}W_{ijss}-nF^{ij}W_{ij}
+\bar H\sum\limits_iF^{ii}\notag\\
&\geq&F(I)(\bar H-n)\geq0.
\end{eqnarray}

At a boundary point, choose an orthonormal frame with
$e_n=\mu$ and diagonalize the tangential block of $W$.
The Robin condition gives $W_{\alpha n}=0$.
Differentiating the equation in the conormal direction and
using \eqref{Waan}, we find
\begin{eqnarray*}
\n_\mu\bar H
=\frac{\cot\theta}{F^{nn}}
\sum\limits_{\alpha=1}^{n-1}
(F^{nn}-F^{\alpha\alpha})(W_{nn}-W_{\alpha\alpha})
\leq0.
\end{eqnarray*}
Here we used Proposition~\ref{pro-2.3} (3) and $\cot\theta>0$.
For the fixed smooth admissible solution, the operator in
\eqref{uni-1} is uniformly elliptic. If $\bar H$ were
nonconstant, the strong maximum principle and the Hopf
boundary point lemma would give $\n_\mu\bar H>0$ at a
boundary maximum, a contradiction. Thus $\bar H$ is constant.
Equation \eqref{uni-1} now implies
\begin{eqnarray*}
\bar H\leq n.
\end{eqnarray*}
Together with \eqref{H-2}, this gives $\bar H\equiv n$.

Recall that
\begin{eqnarray*}
\Delta\ell+n\ell=n\quad\hbox{in }\C_\theta,
\qquad
\n_\mu\ell=\cot\theta\ell\quad\hbox{on }\p\C_\theta.
\end{eqnarray*}
Hence $w\coloneqq h-\ell$ satisfies
\begin{eqnarray*}
\left\{\begin{array}{rcll}
\Delta w+nw&=&0,&\quad\hbox{in }\C_\theta,\\
\n_\mu w&=&\cot\theta w,&\quad\hbox{on }\p\C_\theta.
\end{array}\right.
\end{eqnarray*}
Multiplying the first equation by $\ell$ and integrating
by parts gives
\begin{eqnarray*}
0=\int_{\C_\theta}\ell(\Delta w+nw)d\s
=n\int_{\C_\theta}w\,d\s.
\end{eqnarray*}
Equality therefore holds in the capillary Alexandrov--Fenchel
inequality of \cite[Theorem~3.2]{MWWX}, with the first
function equal to $w$ and all remaining functions equal to $\ell$.
Its equality characterization gives
\begin{eqnarray*}
w=a\ell+\sum\limits_{i=1}^n a_i\<\xi,E_i\>.
\end{eqnarray*}
The equation $\Delta w+nw=0$ implies $a=0$, proving
\eqref{h-uni}. Capillary evenness forces every $a_i$ to vanish.
\end{proof}

Next, we establish the uniqueness part of Theorem~\ref{thm-p-small} in the case $1<p<k-l+1$ with $l=0$. Specifically, we can show that the admissible solution to Eq. \eqref{quo-equ} is unique, which is formalized as the following theorem:
\begin{thm}\label{thm-p-small-uni}
    Let $1\leq k\leq n$, $l=0$ and $\theta\in (0, \pi)$. Suppose that $1<p<k+1$, then the positive strictly convex  solution of Eq. \eqref{quo-equ} is unique. 
\end{thm}

\begin{proof}
The proof follows essentially the same idea as that of Guan--Xia \cite[Section~4.3]{GX}, with the aid of the capillary Alexandrov--Fenchel inequality from  \cite[Theorem 1.1 or Corollary 3.3]{MWWX}.
For the reader’s convenience, we include the details here.

Suppose that there exist two positive strictly convex  solutions $h_{1}$ and $h_{2}$ satisfying Eq.  \eqref{quo-equ}. That is, for $s=1,2$, they satisfy
\begin{eqnarray}\label{equ-1}
    \sigma_{k}(W_{s})=fh^{p-1}_{s}, \quad {\text{in}}~~\C_{\theta}, 
\end{eqnarray}
together with the Robin boundary condition
 \begin{eqnarray}\label{eqn-unique-robin-bdry}
    \n_{\mu}h_{s}=\cot\theta h_{s}, ~~\text{ on } ~~\partial\C_{\theta}.
\end{eqnarray}
By recursively applying the capillary Alexandrov--Fenchel inequality (cf. \cite[Theorem 1.1]{MWWX} or Eq. (3.8)), we obtain
\begin{eqnarray}\label{p-uni-1}
    \int_{\C_{\theta}} h_{2}h_{1}^{p-1}f d\s&=&\int_{\C_{\theta}}h_{2}\sigma_{k}(W_{1})d\s \notag\\
    &\geq& \left(\int_{\C_{\theta}}h_{1}\sigma_{k}(W_{1})d\s\right)^{\frac{k}{k+1}}\left(\int_{\C_{\theta}}h_{2}\sigma_{k}(W_{2})d\s\right)^{\frac{1}{k+1}} \notag \\
    &=&\left(\int_{\C_{\theta}}h_{1}^{p}fd\s \right)^{\frac{k}{k+1}} \left(\int_{\C_{\theta}}h_{2}^{p}fd\s\right )^{\frac{1}{k+1}}.
\end{eqnarray}

On the other hand, by H\"older's inequality
\begin{eqnarray*}
    \int_{\C_{\theta}}h_{2}h_{1}^{p-1}fd\s\leq \left(\int_{\C_{\theta}}h_{2}^{p}fd\s\right)^{\frac{1}{p}}\left(\int_{\C_{\theta}}h_{1}^{p}fd\s\right)^{\frac{p-1}{p}}.
\end{eqnarray*}
Together with \eqref{p-uni-1} and the assumption that $1<p<k+1$, we conclude that
\begin{eqnarray*}
    \int_{\C_{\theta}}h_{1}^{p}fd\s\leq \int_{\C_{\theta}}h_{2}^{p}fd\s.
\end{eqnarray*}

By symmetry, the reverse inequality also holds:
    \begin{eqnarray*}
    \int_{\C_{\theta}}h_{2}^{p}fd\s\leq \int_{\C_{\theta}}h_{1}^{p}fd\s.
\end{eqnarray*}
Therefore, equality must hold in \eqref{p-uni-1}. By \cite[Theorem 1.1]{MWWX}, we obtain
\begin{eqnarray}\label{eqn-unique-h-1-2}
    h_{1}=b h_{2}+\sum_{i=1}^n a_i\<\cdot, E_i\>, ~~\text{ in } ~ \C_\theta
\end{eqnarray}  for some constants $b$ and $a_{i} \in \RR$, $1\leq i\leq n$. Substituting \eqref{eqn-unique-h-1-2} into \eqref{equ-1} and \eqref{eqn-unique-robin-bdry}, we deduce that $$b=1, ~~\text{ and } ~~a_i=0, \text{ for all } 1\leq i\le n.$$ Hence, we conclude that  $h_{1}\equiv h_{2}$.
    
\end{proof}

\begin{lem}\label{oper-inver}
    Let $0\leq l<k\leq n$ and $\theta \in (0, \pi/2)$. Suppose that  $p>k-l+1$ and  $h$ is a positive admissible solution to Eq. \eqref{quo-equ},  the corresponding linearized operator is given by
    \begin{eqnarray*}
        \mathcal{L}_{h}z\coloneqq F^{ij}(z_{ij}+z\delta_{ij})-\frac{p-1}{k-l}h^{a-1}\widehat{f} z, \quad z\in \mathcal{A}^{4, \gamma},
    \end{eqnarray*}
    where $$\mathcal{A}^{4, \gamma}\coloneqq  \left\{z\in C^{4, \gamma}(\C_{\theta}): \n_{\mu}z =\cot\theta z ~  \text{ on }~ \p \C_\theta \right\}.$$  Then $\mathcal L_h:\mathcal A^{4,\gamma}\to C^{2,\gamma}(\C_\theta)$ is invertible.

\end{lem}

\begin{proof}
By Schauder theory, $\mathcal{L}_{h}$ is Fredholm of index zero. Hence it suffices to prove that the kernel of $\mathcal L_h$ is trivial.
Let $z\in\mathcal A^{4,\gamma}$ satisfy
$$\mathcal L_h z=0.$$
Define
$$\tilde z\coloneqq\frac{z}{h}.$$
Then
\begin{eqnarray}\label{eqn:tilde-z-bdry}
\nabla_\mu\tilde z=0
\qquad\text{on }\partial\C_\theta,
\end{eqnarray}
and a direct computation yields
\begin{eqnarray}\label{ker-condi}
    0=\mathcal{L}_{h} z =\left(1-\frac{p-1}{k-l}\right)h^{a}\widehat{f}\tilde{z}+hF^{ij}\tilde{z}_{ij}+2F^{ij}\tilde{z}_{i}h_{j}. 
\end{eqnarray}
Suppose that $\tilde{z}$ attains its maximum value at some $\xi_{0}\in \C_\theta$. If $\xi_{0}\in\C_\theta\setminus\p \C_\theta$, then 
\begin{eqnarray}\label{v-deri}
    \n\tilde{z}=0, \quad {\rm{and}}\quad \n^{2}\tilde{z}\leq 0.
\end{eqnarray}
If $\xi_{0}\in\partial\C_\theta$, then \eqref{v-deri} still holds by virtue of the boundary condition \eqref{eqn:tilde-z-bdry}.

Inserting \eqref{v-deri} into \eqref{ker-condi}, at $\xi_0$, we derive 
\begin{eqnarray*}
   \left (1-\frac{p-1}{k-l} \right) h^{a}\widehat{f}\tilde{z}\geq 0.
\end{eqnarray*}
Since $p>k-l+1$, the coefficient on the left-hand side is negative. It follows that $$\max\limits_{\C_{\theta}}\tilde{z}=\tilde z(\xi_0)=0.$$ The same argument applied to the minimum point yields $$\min\limits_{\C_{\theta}}\tilde{z}=0.$$ 
Therefore, $$\tilde{z}\equiv 0, ~~~\text{ on } ~~ \C_{\theta},$$ and hence $z\equiv 0$ on $\C_{\theta}$. This completes the proof.
  \end{proof}

We are now in a position to prove Theorem \ref{thm-admissible-existence}.

\begin{proof}[\textbf{Proof of Theorem \ref{thm-admissible-existence}}]
The proof is divided into two cases, depending on the value of $p$.

\

\noindent\textbf{Case 1} $p>k-l+1$. We prove the existence by the method of continuity.
Consider the following one-parameter ($t\in[0,1]$) family of equations:
\begin{eqnarray}\label{t-equ}
\left\{
		\begin{array}{rcll}\vspace{2mm}\displaystyle
		  \frac{\sigma_{k}(W)}{\sigma_{l}(W)} &= &h^{p-1}f_{t},& \quad    \hbox{ in }~~ \C_{\theta},\\
			\n_\mu h&=&\cot\theta h, &\quad  \hbox{ on }~~ \p \C_\theta.
		\end{array}\right.
  \end{eqnarray} 

where \begin{eqnarray}\label{ft}
      f_{t}\coloneqq \left[(1-t)\left(\frac{\binom{n}{l}  }{\binom{n}{k}}\ell^{p-1}\right)^{\frac{1}{p+k-l-1}}+t f^{-\frac{1}{p+k-l-1}}\right]^{-(p+k-l-1)}.\end{eqnarray}
      
Define 
\begin{eqnarray*}
    \mathcal{I}\coloneqq  \left\{t~|~ 0\leq t\leq 1,\text{ and~Eq. } \eqref{t-equ}~\text{ admits~a~positive}, \text{admissible~solution} ~h_t\in \mathcal{A}^{4, \gamma} \right\}.
\end{eqnarray*}
Since $h\equiv\ell$ is a positive admissible solution of \eqref{t-equ} when
$t=0$, we have $0\in\mathcal I$. Hence $\mathcal I$ is nonempty.  The openness of $\mathcal I$ follows from Lemma~\ref{oper-inver} and the implicit function theorem, while the closedness follows from Theorem~\ref{thm priori est}. Hence, $\mathcal I=[0,1]$, which establishes the existence of a positive admissible solution to \eqref{quo-equ}. The uniqueness follows from Theorem~\ref{thm-unique}.

\ 

We next prove the strict convexity of the admissible solution under the additional assumption that $f$ satisfies \eqref{f-condition}. 
Denote
\begin{eqnarray*}
  \bar{c}_{n,k,l}\coloneqq \frac{\binom{n}{k}}{\binom{n}{l}},
 \qquad\lambda\coloneqq \frac{p-1}{p+k-l-1},
\end{eqnarray*}
for any $0\leq \lambda <1$, we set 
\begin{eqnarray*}
g_t\coloneqq f_t^{-1/(p+k-l-1)}
=(1-t)\bar{c}_{n,k,l}^{-1/(p+k-l-1)}\ell^\lambda+t f^{-1/(p+k-l-1)}.
\end{eqnarray*}
The identities $\n^2\ell=(1-\ell)\s$ and
$|\n\ell|^2\leq\ell$ on the acute spherical cap imply
\begin{eqnarray*}
\n^2\ell^\lambda+\ell^\lambda\s
&=&\ell^{\lambda-2}
\left[(\lambda\ell+(1-\lambda)\ell^2)\s
-\lambda(1-\lambda)\n\ell\otimes\n\ell\right]\\
&\geq&\ell^{\lambda-1}
\left[\lambda^2+(1-\lambda)\ell\right]\s>0.
\end{eqnarray*}
When $k<n$, from \eqref{f-condition}, we have
\begin{eqnarray}\label{ft-convex}
        \n^{2}f_{t}^{-\frac{1}{p+k-l-1}}+f_{t}^{-\frac{1}{p+k-l-1}}\sigma \geq 0,  \quad {\text{in}}\quad \C_{\theta},
    \end{eqnarray}
and 
\begin{eqnarray*}
&&\n_\mu g_t-\cot\theta g_t\\
&=&(1-t)\bar{c}_{n,k,l}^{-1/(p+k-l-1)}(\lambda-1)\cot\theta\ell^\lambda
+t\left(\n_\mu(f^{-1/(p+k-l-1)})-\cot\theta f^{-1/(p+k-l-1)}\right)
\leq0.
\end{eqnarray*}
Since $f_t=g_t^{-(p+k-l-1)}$, the above inequality is equivalent to
\begin{eqnarray}
\n_\mu f_t+(p+k-l-1)\cot\theta f_t\geq0.\label{ft-robin}
\end{eqnarray}
Thus, for $t\in [0,1]$, we have that $f_{t}$ satisfies \eqref{f-convex} and \eqref{f-bry}.
 
 When $t=0$, in view of Corollary \ref{cor:5.2-uniqueness}, the solution of \eqref{t-equ} is given by $$h\equiv\ell,$$ which is strictly convex. Suppose, for contradiction, that there exists $t_0\in(0,1]$ such that $h_t$ is strictly convex for all $0\le t<t_0$, whereas $h_{t_0}$ fails to be strictly convex. Since $f_{t_0}$ satisfies both \eqref{f-convex} and \eqref{f-bry}, Theorem~\ref{thm-convex} implies that $h_{t_0}$ is strictly convex, a contradiction. Therefore, $h_t$ is strictly convex for every $t\in[0,1]$. In particular, the solution $h=h_1$ is strictly convex.

 \

 \noindent \textbf{Case 2} $p=k-l+1$.
We employ an approximation argument following
Guan--Lin \cite[Section~3]{GL2000}
and Hu--Ma--Shen \cite[Section~4]{HMS2004}.

For each $\varepsilon\in(0,1)$, \textbf{ Case~1} yields a unique positive smooth admissible solution
$h_\varepsilon$ of
\begin{eqnarray}\label{eqn:eq of critical-h-varepsilon}
\left\{	\begin{array}{rcll}\vspace{2mm}\displaystyle
			 \frac{\sigma_{k}(W)}{\sigma_{l}(W)}&=&h^{k-l+\varepsilon} f, &\quad    \hbox{ in } \C_{\theta},\\
			\n_\mu h &=& \cot\theta  h, &\quad  \hbox{ on } \p \C_\theta. 
		\end{array}\right.
	\end{eqnarray} 
Define $$\tilde{h}_{\varepsilon}\coloneqq \frac{h_{\varepsilon}}{\min\limits_{\C_{\theta}}h_{\varepsilon}}.$$    By \eqref{eqn:eq of critical-h-varepsilon}, then $\tilde{h}_{\varepsilon}$ satisfies 
\begin{eqnarray*}
\left\{	\begin{array}{rcll}\vspace{2mm}\displaystyle
		  \frac{\sigma_{k}(W)}{\sigma_{l}(W)}&=&(\min\limits_{\C_{\theta}}h_{\varepsilon})^{\varepsilon}h^{k-l+\varepsilon} f, &\quad    \hbox{ in } \C_{\theta},\\
			\n_\mu h &=& \cot\theta  h, &\quad  \hbox{ on } \p \C_\theta.
		\end{array}\right.
	\end{eqnarray*} 
By Theorem \ref{thm priori est} (2),  we have
\begin{eqnarray*}
    \|\tilde{h}_{\varepsilon}\|_{C^{4,\alpha}(\C_{\theta})}\leq C,
\end{eqnarray*}
where $C$ depends only on $n$, $k$, $l$, and $f$, and is independent of $\varepsilon$. Passing to a subsequence if necessary, there exist $\varepsilon_j\to0$ and a function $\tilde h\in C^4(\C_\theta)$ such that
\begin{eqnarray*}
    \tilde h_{\varepsilon_j}\to\tilde h
\qquad\text{in }C^4(\C_\theta).
\end{eqnarray*} 
By Lemma~\ref{lem kC0} (1), $$(\min\limits_{\C_{\theta}}h_{\varepsilon_{j}})^{\varepsilon_{j}}\rightarrow \tau,$$ for some constant $\tau>0$, where 
\begin{eqnarray*}
\frac{\binom{n}{k} }{\binom{n}{l}}\frac{1}{\max\limits_{ \C_{\theta}}f \cdot (\sin\theta)^{2(k-l)}}\leq       \tau
    \leq\frac{\binom{n}{k} }{\binom{n}{l}}\frac{1}{\min\limits_{\C_{\theta}}f \cdot (1-\cos\theta)^{k-l}}.
    \end{eqnarray*}
Consequently, $\tilde h$ is a positive smooth admissible solution of
\begin{eqnarray*}
\left\{	\begin{array}{rcll}\vspace{2mm}\displaystyle
	\frac{\sigma_{k}(W)}{\sigma_{l}(W)}&=&\tau  f h^{k-l},	   &\quad    \hbox{ in } \C_{\theta},\\
			\n_\mu h &=& \cot\theta  h, &\quad  \hbox{ on } \p \C_\theta.
		\end{array}\right.
	\end{eqnarray*} 

Assume now that $k<n$ and that $f$ satisfies
\eqref{f-condition} with $p=k-l+1$. For $s\in(0,1)$, set
\begin{eqnarray*}
f_s\coloneqq\left[(1-s)
\left(\frac{\binom{n}{l}}{\binom{n}{k}}
\ell^{k-l}\right)^{\frac1{2(k-l)}}
+s f^{-\frac1{2(k-l)}}\right]^{-2(k-l)}.
\end{eqnarray*}
The calculation leading to \eqref{ft-convex} and
\eqref{ft-robin}, with $\lambda=1/2$, shows that
\begin{eqnarray*}
\n^2 f_s^{-\frac1{2(k-l)}}
+f_s^{-\frac1{2(k-l)}}\s>0,
\qquad
\n_\mu f_s+2(k-l)\cot\theta f_s\geq0.
\end{eqnarray*}
For each fixed $s\in(0,1)$, strict positivity and
continuity in $\varepsilon$ imply
\begin{eqnarray*}
\n^2 f_s^{-\frac1{2(k-l)+\varepsilon}}
+f_s^{-\frac1{2(k-l)+\varepsilon}}\s>0
\end{eqnarray*}
for all sufficiently small $\varepsilon>0$.
The boundary condition is preserved because
\begin{eqnarray*}
\n_\mu f_s+(2(k-l)+\varepsilon)\cot\theta f_s
\geq\varepsilon\cot\theta f_s>0.
\end{eqnarray*}
Thus Case~1, applied with $f_s$ and
$p=k-l+1+\varepsilon$, gives strictly convex solutions.
Normalize these solutions to have minimum one and pass
to the limit as $\varepsilon\downarrow0$, as above.
We obtain a positive admissible pair $(h_s,\tau_s)$ with
$\min_{\C_\theta}h_s=1$, $W[h_s]\geq0$, and
\begin{eqnarray}\label{pari-equ}
\left\{\begin{array}{rcll}\vspace{2mm}\displaystyle
\frac{\sigma_k(W[h_s])}{\sigma_l(W[h_s])}
&=&\tau_s f_s h_s^{k-l},&\quad\hbox{in }\C_\theta,\\
\n_\mu h_s&=&\cot\theta h_s,
&\quad\hbox{on }\partial\C_\theta.
\end{array}\right.
\end{eqnarray}
Multiplication of $f_s$ by the positive constant $\tau_s$
preserves the structural conditions. Therefore
Theorem~\ref{thm-convex} gives $W[h_s]>0$.
The normalized estimates in
Theorem~\ref{thm priori est}\eqref{thm-C2-est-p-middle}
pass to this limit. Since $f_s\to f$ smoothly and has
uniform positive lower and smooth upper bounds for
$s$ close to one, these estimates and the preceding
bounds for $\tau_s$ are uniform in $s$.
After passing to a subsequence,
\begin{eqnarray*}
\tau_{s_i}\to\tau>0,
\qquad h_{s_i}\to h\quad\hbox{in }C^4(\C_\theta),
\qquad\min_{\C_\theta}h=1.
\end{eqnarray*}
The limit is positive and admissible, satisfies
\eqref{homo-equ}, and has $W[h]\geq0$.
Theorem~\ref{thm-convex}, applied to $\tau f$, now yields
$W[h]>0$. For $k=n$, strict convexity follows directly
from admissibility.

 The uniqueness of the positive constant $\tau$ and the strictly convex solution $h$ (up to dilations) of \eqref{homo-equ} can be proved by the same argument as in \cite[Section~5]{Lee}, and we omit the details for conciseness. This completes the proof.
\end{proof}

\

\subsection{Proof of Theorem \ref{thm-p-small}}
We prove existence for $1\leq p<k-l+1$ by a degree argument
in the space of capillary even functions.

Fix $0<\epsilon<\alpha<1$ and write
$\widehat\xi=(-\xi_1,\ldots,-\xi_n,\xi_{n+1})$.
Define
\begin{eqnarray*}
\mathcal S^{4,\epsilon}
\coloneqq\left\{h\in C^{4,\epsilon}(\C_\theta):
h(\xi)=h(\widehat\xi),\quad
\n_\mu h=\cot\theta h\ \hbox{on }\p\C_\theta\right\},
\end{eqnarray*}
and
\begin{eqnarray*}
\mathcal Y^{2,\epsilon}
\coloneqq\left\{v\in C^{2,\epsilon}(\C_\theta):
v(\xi)=v(\widehat\xi)\right\}.
\end{eqnarray*}
Let
\begin{eqnarray*}
\mathcal U\coloneqq\left\{h\in\mathcal S^{4,\epsilon}:
h>0,\quad \n^2h+h\s>0\right\}.
\end{eqnarray*}
For $f_t$ defined by \eqref{ft}, set
$G:\mathcal U\times[0,1]\longrightarrow\mathcal Y^{2,\epsilon}$,
\begin{eqnarray*}
G(h,t)\coloneqq
\frac{\sigma_k(\n^2h+h\s)}{\sigma_l(\n^2h+h\s)}
-h^{p-1}f_t.
\end{eqnarray*}
The functions $f_t$ are positive, smooth, and capillary even,
with a uniform positive lower bound and uniform smooth norms. Moreover, from \eqref{ft-convex} and \eqref{ft-robin}, we see  that $f_{t}$ satisfies \eqref{f-convex} and \eqref{f-bry} whenever $1\leq p<k-l+1$ and $t\in [0,1]$.

\begin{lem}\label{kernel}
Let $1\leq p<k-l+1$ and $\theta\in(0,\pi/2)$.
The linearization of $G(\cdot,0)$ at $h=\ell$ is an
isomorphism from $\mathcal S^{4,\epsilon}$ to
$\mathcal Y^{2,\epsilon}$.
\end{lem}
\begin{proof}
The linearized operator is
\begin{eqnarray*}
\mathcal L\varphi
&\coloneqq&\frac{d}{ds}\bigg|_{s=0}G(\ell+s\varphi,0)\\
&=&a_0\left(\Delta\varphi+n\varphi-b_0\ell^{-1}\varphi\right),
\end{eqnarray*}
where
\begin{eqnarray*}
a_0=\frac{k-l}{n}\frac{\binom{n}{k}}{\binom{n}{l}},
\qquad b_0=\frac{n(p-1)}{k-l}.
\end{eqnarray*}
Suppose that $\varphi\in\mathcal S^{4,\epsilon}$ and
$\mathcal L\varphi=0$. Then
\begin{eqnarray}\label{kernel-zero}
\Delta\varphi+n\varphi-b_0\ell^{-1}\varphi=0.
\end{eqnarray}
Multiplying by $\ell$ and integrating by parts, we obtain
\begin{eqnarray*}
b_0\int_{\C_\theta}\varphi\,d\s
&=&\int_{\C_\theta}\ell(\Delta\varphi+n\varphi)d\s\notag\\
&=&\int_{\p\C_\theta}
(\ell\n_\mu\varphi-\varphi\n_\mu\ell)ds
+\int_{\C_\theta}\varphi(\Delta\ell+n\ell)d\s\notag\\
&=&n\int_{\C_\theta}\varphi\,d\s.
\end{eqnarray*}
Since $b_0<n$, it follows that
\begin{eqnarray*}
\int_{\C_\theta}\varphi\,d\s=0.
\end{eqnarray*}
Applying \cite[Theorem~3.2]{MWWX} with the first function
equal to $\varphi$ and all remaining functions equal to $\ell$,
we obtain
\begin{eqnarray*}
\left(\int_{\C_\theta}\varphi\,d\s\right)^2
\geq\frac1n
\left(\int_{\C_\theta}\varphi(\Delta\varphi+n\varphi)d\s\right)
\left(\int_{\C_\theta}\ell\,d\s\right).
\end{eqnarray*}
Therefore,
\begin{eqnarray*}
\int_{\C_\theta}\varphi(\Delta\varphi+n\varphi)d\s\leq0.
\end{eqnarray*}
Multiplying \eqref{kernel-zero} by $\varphi$ gives
\begin{eqnarray*}
b_0\int_{\C_\theta}\ell^{-1}\varphi^2d\s
=\int_{\C_\theta}\varphi(\Delta\varphi+n\varphi)d\s\leq0.
\end{eqnarray*}
If $p>1$, then $b_0>0$, so $\varphi\equiv0$.
If $p=1$, then $b_0=0$ and equality holds in the preceding
Alexandrov--Fenchel inequality. Its equality characterization
gives
\begin{eqnarray*}
\varphi=c\ell+\sum\limits_{i=1}^n c_i\<\xi,E_i\>.
\end{eqnarray*}
Equation \eqref{kernel-zero} gives $c=0$.
Since $\varphi$ is capillary even and each horizontal coordinate
function changes sign under $\xi\mapsto\widehat\xi$, every
$c_i$ vanishes. Thus the kernel is trivial also at $p=1$.

The Robin realization of $\Delta+n-b_0\ell^{-1}$ is
self-adjoint with compact resolvent and commutes with the
reflection $\xi\mapsto\widehat\xi$.
Its restriction to the capillary even subspace is Fredholm
of index zero. The triviality of the kernel and elliptic
regularity therefore imply that
$\mathcal L:\mathcal S^{4,\epsilon}\to\mathcal Y^{2,\epsilon}$
is an isomorphism.
\end{proof}

\ 

We now turn to the proof of Theorem~\ref{thm-p-small}. Together with the proof of Theorem~\ref{thm-admissible-existence}, this will complete the proofs of Theorems~\ref{thm-cap-Lp-CM-thm} and~\ref{thm:1.4}. 
\begin{proof}[\textbf{Proof of Theorem \ref{thm-p-small}}]
We consider the homotopy equations
\begin{eqnarray}\label{ft-equ}
\left\{\begin{array}{rcll}\vspace{2mm}\displaystyle
\frac{\sigma_k(\n^2h+h\s)}{\sigma_l(\n^2h+h\s)}
&=&h^{p-1}f_t,&\quad\hbox{in }\C_\theta,\\
\n_\mu h&=&\cot\theta h,&\quad\hbox{on }\p\C_\theta.
\end{array}\right.
\end{eqnarray}
Theorem~\ref{thm priori est} gives constants $c_0,C_0>0$,
independent of $t\in[0,1]$, such that every positive,
capillary even, strictly convex solution satisfies
\begin{eqnarray*}
h\geq c_0,\qquad \|h\|_{C^{4,\alpha}(\C_\theta)}\leq C_0.
\end{eqnarray*}
For $p=1$, these bounds follow from the separate endpoint
argument in that theorem.

We claim that there is also a uniform constant $c_1>0$ such that
\begin{eqnarray*}
W[h]=\n^2h+h\s\geq c_1\s.
\end{eqnarray*}
Otherwise, there would be solutions $h_j$, parameters $t_j$,
and points $\xi_j$ with
$\lambda_{\min}(W[h_j](\xi_j))\to0$.
Since $\epsilon<\alpha$, passing to a subsequence gives
\begin{eqnarray*}
t_j\to t_\infty,\qquad \xi_j\to\xi_\infty,\qquad
h_j\to h_\infty\quad\hbox{in }C^{4,\epsilon}(\C_\theta).
\end{eqnarray*}
The lower bound for $f_t h^{p-1}$ and the Newton--Maclaurin
inequalities give a uniform positive lower bound for
$\sigma_k(W[h_j])$.
Thus $h_\infty$ is positive and admissible,
$W[h_\infty]\geq0$, and $G(h_\infty,t_\infty)=0$.
If $k<n$, Theorem~\ref{thm-convex} implies $W[h_\infty]>0$.
If $k=n$, this follows directly from
$W[h_\infty]\geq0$ and $\sigma_n(W[h_\infty])>0$.
Both conclusions contradict the assumed degeneration.
This proves the claim.

Choose $\delta>0$ sufficiently small and $C_s>0$ sufficiently
large that all the preceding solutions and $\ell$ belong to
\begin{eqnarray*}
\mathcal S\coloneqq\left\{h\in\mathcal S^{4,\epsilon}:
h>\delta,\quad \n^2h+h\s>\delta\s,\quad
\|h\|_{C^{4,\epsilon}(\C_\theta)}<C_s\right\}.
\end{eqnarray*}
Then $G(h,t)\neq0$ on $\p\mathcal S\times[0,1]$.
We use the degree construction for elliptic operators
with oblique boundary conditions in \cite[Section~2]{LLN},
with the fixed homogeneous Robin condition incorporated
in the domain. The reflection preserves the equation,
the boundary condition, and the auxiliary elliptic
solution operators used in that construction.
The resulting compact maps therefore preserve the
capillary even subspace, on which we take the degree.
We denote this reduced degree by $\text{deg}$.
Homotopy invariance gives
\begin{eqnarray}\label{deg-1}
\text{deg}(G(\cdot,0),\mathcal S,0)
=\text{deg}(G(\cdot,1),\mathcal S,0).
\end{eqnarray}

When $p=1$, Theorem~\ref{thm-constant-p-one} and capillary
evenness imply that $\ell$ is the only solution of
$G(h,0)=0$ in $\mathcal S$.
For $1<p<k-l+1$, the rigidity result
\cite[Theorem~1.6, Theorem~1.10 (i)]{GLi} implies that the corresponding
hypersurface is a spherical cap. Capillary evenness centres
this cap on the symmetry axis, so $h=R\ell$.
Substitution into the model equation gives
$R^{k-l}=R^{p-1}$, and hence $R=1$.
Thus the only model solution in either case is $h=\ell$.

By Lemma~\ref{kernel}, the linearization $\mathcal L$ at
$\ell$ is an isomorphism.
Choose $\rho>0$ so that
$\ell+\overline{B_\rho(0)}\subset\mathcal S$,
where $B_\rho(0)$ is a ball in $\mathcal S^{4,\epsilon}$.

By excision and the local linearization and normalization
properties of the oblique elliptic degree
\cite[Section~2]{LLN} (see also 
\cite[Propositions~2.3 and~2.4]{LYY}),
\begin{eqnarray}\label{deg-2}
\text{deg}(G(\cdot,0),\mathcal S,0)
&=&\text{deg}(G(\ell+\cdot,0),B_\rho(0),0)\notag\\
&=&\text{deg}(\mathcal L,B_\rho(0),0)=\pm1.
\end{eqnarray}
Combining \eqref{deg-1} and \eqref{deg-2}, we obtain
$\text{deg}(G(\cdot,1),\mathcal S,0)\neq0$.
Hence \eqref{ft-equ} has a positive, capillary even,
strictly convex solution at $t=1$.
Elliptic regularity makes this solution smooth.
When $l=0$ and $1<p<k+1$, uniqueness follows from
Theorem~\ref{thm-p-small-uni}.
This completes the proof.
\end{proof}
\ 

\begin{proof}[\textbf{Proof of Theorems \ref{thm-cap-Lp-CM-thm} and \ref{thm:1.4}}]
By Definition~\ref{defn:2.10} and Proposition \ref{prop:cm-problem-lp}, 
Theorem~\ref{thm-cap-Lp-CM-thm} follows immediately from Theorems \ref{thm-admissible-existence} and \ref{thm-p-small} upon setting $l=0$ in Eq. \eqref{quo-equ}. Theorem~\ref{thm:1.4} is then an immediate consequence of Theorems \ref{thm-admissible-existence} and \ref{thm-p-small}.

\end{proof}

\bigskip 

\ 

\noindent\textit{Acknowledgment:} X.M. would like to express sincere gratitude to Prof. Yuguang Shi for his ongoing support and encouragement, and to Dr. Wei Wei for some useful discussions on the proof of Lemma \ref{lem kC1}. X.M. was supported by the National Key R $\&$ D Program of China (No. 2020YFA0712800) and the Postdoctoral Fellowship Program of CPSF under Grant Numbers 2025T180843 and 2025M773082.  L.W. was partially supported by CRM De Giorgi of Scuola Normale Superiore. L.W. is a member of GNAMPA as part of INdAM.

\bigskip

\printbibliography

\end{document}